\documentclass[reqno,11pt]{amsart}
\UseRawInputEncoding
\pdfoutput=1
\usepackage{amsfonts}
\usepackage{graphicx}
\usepackage{amsfonts,amsmath, amssymb}
\usepackage{bbm,mathrsfs,mathscinet}
\usepackage{hyperref}
\usepackage{marginnote}
\usepackage{color}
\usepackage{verbatim}
\usepackage{mathrsfs}
\usepackage{cancel}
\allowdisplaybreaks

\numberwithin{equation}{section}

\newcommand{\R}{\mathbb{R}}

\newtheorem{theorem}{Theorem}[section]
\newtheorem{corollary}[theorem]{Corollary}
\newtheorem{lemma}[theorem]{Lemma}
\newtheorem{proposition}[theorem]{Proposition}
\newtheorem{remark}[theorem]{Remark}

\def\v{\varepsilon}
\def\var{\varepsilon}

\begin{document}
	
	\title[Low Mach number Limit of Steady Thermally Driven Fluid]{Low Mach number Limit of Steady Thermally Driven Fluid}
	
	\author[F.-M. Huang]{Feimin Huang}
	\address[F.-M. Huang]{Academy of Mathematics and Systems Science, Chinese Academy of Sciences, Beijing 100190, China; School of Mathematical Sciences, University of Chinese Academy of Sciences, Beijing 100049, China.}
	\email{fhuang@amt.ac.cn}
	
	\author[W.-Q. Wang]{Weiqiang Wang}
	\address[W.-Q. Wang]{Academy of Mathematics and Systems Science, Chinese Academy of Sciences, Beijing 100190, China.}
	\email{wangweiqiang@amss.ac.cn}
	
	\author[Y. Wang]{Yong Wang}
	\address[Y. Wang]{Academy of Mathematics and Systems Science, Chinese Academy of Sciences, Beijing 100190, China; School of Mathematical Sciences, University of Chinese Academy of Sciences, Beijing 100049, China.}
	\email{yongwang@amss.ac.cn}
	
\begin{abstract}
In this paper, we establish the existence of strong solutions to the steady non-isentropic compressible Navier-Stokes system with Dirichlet boundary conditions in bounded domains where the fluid is driven by the wall temperature, and justify its low Mach number limit, i.e., $\v\to 0$, in $L^{\infty}$ sense with a rate of convergence. Notably, for the limiting system \eqref{fge} obtained in the low Mach number limit, the variation of the wall temperature is allowed to be independent of the Mach number. It is also worth pointing out that the
velocity field $u_{1}$ acts like a ghost since it appears at $\v$-order in the expansion, but still affects the density and temperature at $O(1)$-order. In the proof, we design a new expansion, in which the density, velocity and temperature have different expansion forms with respect to $\v$, so that the density at higher orders is well-defined under the Boussinesq relations and the constraint of zero average. We also introduce a new $\v$-dependent functional space, allowing us to obtain some uniform estimates for high-order normal derivatives near the boundary.
\end{abstract}

\date{\today}
%\thanks{}
\maketitle

\setcounter{tocdepth}{2}
\tableofcontents

\thispagestyle{empty}

\section{Introduction and Main Result}

\subsection{Introduction}
The steady non-isentropic compressible  Navier-Stokes system reads
\begin{equation}\label{sNS}
	\left\{
	\begin{aligned}
		&\operatorname{div}(\rho^{\v} \mathfrak{u}^{\v})=0,\\
		&\rho^{\v}( \mathfrak{u}^{\v}\cdot \nabla) \mathfrak{u}^{\v}+\nabla P^{\v}=\varepsilon\operatorname{div}\mathbb{S}(\nabla \mathfrak{u}^{\v}),\\
		&\rho^{\v}( \mathfrak{u}^{\v}\cdot \nabla)\theta^{\v}+P^{\v}\operatorname{div}\mathfrak{u}^{\v}=\varepsilon\operatorname{div}(\kappa\nabla\theta^{\v})+\varepsilon\Psi(\nabla \mathfrak{u}^{\v}),
	\end{aligned}\qquad x\in \Omega
	\right.
\end{equation}
where $\v>0$ is the viscosity. Here $\Omega\subset \mathbb{R}^{3}$ is a smooth bounded domain, $\rho^{\v}, \mathfrak{u}^{\v}$ and $\theta^{\v}$ are the density, velocity and temperature respectively. $P^{\v}=c_{v}\rho^{\v} \theta^{\v}$ denotes the pressure with $c_{v}$ being the gas constant. $\kappa>0$ is the heat conductivity. The stress tensor $\mathbb{S}(\nabla \mathfrak{u}^{\v})$ and the dissipation function $\Psi(\nabla \mathfrak{u}^{\v})$ are defined by
\begin{align}\label{1.1}
	\mathbb{S}(\nabla \mathfrak{u}^{\v})=2\mu D(\mathfrak{u}^{\v})+\lambda\operatorname{div}\mathfrak{u}^{\v}\mathbb{I},\qquad \Psi(\nabla \mathfrak{u}^{\v})=2\mu D(\mathfrak{u}^{\v}):D(\mathfrak{u}^{\v})+\lambda|\operatorname{div}\mathfrak{u}^{\v}|^2,
\end{align}
where $\mathbb{I}$ is the identical matrix, $D(\mathfrak{u}^{\v})=\frac{1}{2}(\nabla \mathfrak{u}^{\v}+(\nabla \mathfrak{u}^{\v})^{t})$ is the deformation tensor, $\mu>0$ and $\lambda$ are the viscosity coefficients satisfying the physical requirement: $2\mu+3\lambda\geq 0$.

We denote $\mathscr{M}_{a}$ as the Mach number. In the present paper, we consider the low Mach number limit, i.e., $\mathscr{M}_{a}\to 0$, of \eqref{sNS} under the scaling $\mathscr{M}_{a}=\v$ and  $|\mathfrak{u}^{\v}|\lesssim\varepsilon$ ({\it cf}. \cite{Lions}). Formally, assuming $\frac{\mathfrak{u}^{\v}}{\v}\to u$ as $\varepsilon\to 0$, one derives the limiting system:
\begin{align}\label{fge}
	\left\{
	\begin{aligned}
		&\nabla (\rho\theta)=0,\,\,\operatorname{div}(\rho u)=0,\\
		&\rho(u\cdot \nabla)u+\nabla P=\operatorname{div}\mathbb{S}(\nabla u),\\
		&\operatorname{div}(\kappa\nabla \theta)=(1+c_{v})\rho\theta\operatorname{div}u.
	\end{aligned}
	\right.
\end{align}

Many works have been made on the low Mach number limit of evolutionary isentropic /non-isentropic compressible Navier-Stokes system both with the well-prepared and ill-prepared initial data respectively, and we refer the interested readers to the survey papers \cite{Alazard-2008,Feireisl-2018,Jiang-Masmoudi} and the references therein. We also would like to mention recent works \cite{Ju-Ou-2022,Sun-2022}, in which they justified the low Mach number limit of evolutionary non-isentropic compressible Navier-Stokes equations with large temperature variations and Navier-slip boundary conditions in bounded domains.

For the steady case, to the best of our knowledge, there are only a few results \cite{Veiga-1987,Choe-Jin,Dou-Jiang-Jiang-Yang}.
In 1987, Beir\"{a}o da Veiga \cite{Veiga-1987} established the existence and uniqueness for stationary solution of the non-isentropic compressible Navier-Stokes with small external forces near a rest state in $L^p$-setting, and justified the low Mach number limit from isentropic compressible Navier-Stokes with small external forces to the incompressible Navier-Stokes system. Later, Choe-Jin \cite{Choe-Jin} established the existence of stationary isentropic compressible Navier-Stokes system with general large external forces when the Mach number is small, and then justified its incompressible limit. Recently, Dou-Jiang-Jiang-Yang  \cite{Dou-Jiang-Jiang-Yang} extended the result in \cite{Choe-Jin} to the non-isentropic Navier-Stokes system. We would like to point out that $O(1)$ temperature variations are not allowed in \cite{Veiga-1987,Choe-Jin,Dou-Jiang-Jiang-Yang} since they essentially need the wall temperature to be around a constant and $\nabla (\theta\vert_{\partial\Omega})\lesssim \varepsilon$. Consequently, the limiting system obtained in \cite{Veiga-1987,Choe-Jin,Dou-Jiang-Jiang-Yang} is the incompressible Navier-Stokes system rather than \eqref{fge}.

We also would like to mention a very important progress \cite{Esposito-Guo-Marra-Wu-2023-1,Esposito-Guo-Marra-Wu-2023-2} made by Esposito-Guo-Marra-Wu on the low Mach number limit of the steady Boltzmann equation in bounded domains.  They rigorously justified the hydrodynamic limit from the steady Boltzmann equation to the following ghost effect system:
%was made by Esposito {\it et al} \cite{Esposito-Guo-Marra-Wu-2023-1,Esposito-Guo-Marra-Wu-2023-2}, in which they rigorously justified the limit from the steady Boltzmann equation to the ghost effect system in a smooth bounded domain.
%Recall the Ghost effect system derived from the low Mach number limit of the Boltzmann equation ({\it cf}. \cite{Bardos-Levermore-Ukai-Yang,Esposito-Guo-Marra-Wu-2023-2}):
\begin{align}\label{ge}
	\left\{\begin{aligned}
	&\nabla (\rho\theta)=0,\quad \operatorname{div}(\rho u)=0,\\
	&\rho(u\cdot \nabla u)+\nabla\mathfrak{p}=\operatorname{div}\Big[\mu\big(2D(u)-\frac{2}{3}\operatorname{div}u\mathbb{I}\big)\Big]\\
	&\qquad\qquad\qquad\quad\,\,\,+ \operatorname{div}\Big[\frac{\mu^2}{\rho\theta}\big(K_{1}(\nabla^2\theta-\frac{1}{3}\Delta\theta\mathbb{I}\big)+\frac{K_{2}}{\theta}(\nabla\theta\otimes \nabla\theta-\frac{1}{3}|\nabla\theta|^2\mathbb{I})\big)\Big],\\
	&\operatorname{div}(\kappa\nabla\theta)=K_{3}\rho\theta\operatorname{div}u,
	\end{aligned}
	\right.
\end{align}
where $K_{i}\,\,(i=1,2,3)$ are some constants. The ghost effect system is used to describe an interesting phenomenon that the flow is only driven by the variation of the temperature ({\it cf. }\cite{Huang-2015}). Since it is quite different from the usual heat flow driven by the difference of pressure, Y. Sone \cite{Sone-2007} suggested the name of ``ghost effect" to such phenomenon because the vanishing small velocity field, playing like a ghost, still affects macroscopic equations of temperature. This phenomenon was analyzed in \cite{Chen-Chen-Liu-Sone} from the view of steady linearized Boltzmann equation between two parallel plates and in a circular pipe respectively. Later, by using the Hilbert expansion, Bardos-Levermore-Ukai-Yang \cite{Bardos-Levermore-Ukai-Yang} gave a formal limit from the Boltzmann equation to the ghost effect system.  We refer \cite{Huang-Wang-Wang-Yang} for the 1-D case with rigorous justification of the formal limit in \cite{Bardos-Levermore-Ukai-Yang}. The well-posedness of the ghost effect system was studied in \cite{Huang-Tan,Levermore-Sun-Trivisa}.

Compared with \eqref{ge}, the system \eqref{fge} only looses some terms involving with derivatives of temperature up to third order in the momentum equation $\eqref{fge}_{2}$. Such discrepancy arises due to the fact that the compressible Navier-Stokes equation is only the first-order approximation of the Boltzmann equation, and some microscopic terms are already ignored in the level of compressible Navier-Stokes equations. Both \eqref{fge} and \eqref{ge} indicate that the $O(1)$-temperature and density are affected by the $O(\varepsilon)$-velocity. Inspired by the work \cite{Esposito-Guo-Marra-Wu-2023-1,Esposito-Guo-Marra-Wu-2023-2} made by Esposito-Guo-Marra-Wu, it is an interesting problem to study the low Mach number limit from \eqref{sNS} to \eqref{fge} in smooth bounded domain.

\subsection{Formulation and main result} In present paper, for any given Mach number $\v>0$, we focus on the existence of strong solutions of \eqref{sNS} with $o(1)$ temperature variations on the boundary, and its rigorous low Mach number limit to \eqref{fge}. Without of lose of generality, we assume $|\Omega|=\int_{\Omega}1\,{\rm d}x=1$ and the gas constant $c_{v}=1$. To avoid technicalities, we will consider the case that the transport coefficients $\mu$, $\lambda$ and $\kappa$ are constants.
%The general case that the transport coefficients are smooth function of $\theta$ can also be operated similarly, whose proof will be sketched in the Appendix \ref{AppendixA}.
As in \cite{Esposito-Guo-Marra-Wu-2023-1}, we supplement the velocity $\mathfrak{u}$ and the temperature $\theta$ with following Dirichlet boundary conditions:
\begin{align}\label{bd}
	\theta^{\v}=T_{w},\quad(\mathfrak{u}^{\v}\cdot \vec{\iota}_{1}, \mathfrak{u}^{\v}\cdot \vec{\iota}_{2},\mathfrak{u}^{\v}\cdot \vec{n})=\varepsilon h(T_{w})(\partial_{\vec{\iota}_{1}}T_{w},\partial_{\vec{\iota}_{2}}T_{w},0)\quad \text{ on }\partial\Omega
\end{align}
where $T_{w}=1+O(|\nabla T_{w}|)$, $h(T_{w})$ is a given smooth function of $T_{w}$, $\vec{\iota}_{i},\,\,i=1,2,$ are the unit tangential vectors of $\partial\Omega$, and $\vec{n}$ is the unit normal vector of $\partial\Omega$.
%Since we consider the steady solution of \eqref{eNS}, we can rewrite \eqref{eNS} as
%\begin{equation}\label{sNS}
%	\left\{
%	\begin{aligned}
		%&\operatorname{div}(\rho \mathfrak{u})=0,\\
		%&\rho( \mathfrak{u}\cdot \nabla) \mathfrak{u}+\nabla P=\varepsilon\operatorname{div}\mathbb{S}(\nabla \mathfrak{u}),\\
		%&\rho( \mathfrak{u}\cdot %\nabla)\theta+P\operatorname{div}\mathfrak{u}=\varepsilon\operatorname{div}(\kappa(\theta)\nabla\theta)+%\varepsilon\Psi(\nabla \mathfrak{u}),
%	\end{aligned}\qquad x\in \Omega.
%	\right.
%\end{equation}
Moreover, the total mass is prescribed:
\begin{align}\label{1.2}
	\int_{\Omega}\rho^{\v}(x)\,{\rm d}x=M.	
\end{align}

We consider the solution of \eqref{sNS} in the following expansion:
\begin{align}\label{expansion}
\rho^{\v}=\rho_{0}+\varepsilon\rho_{1}+\varepsilon^2\rho_{2}+\varepsilon^3\rho_{3}+\varepsilon^2\rho_{R},\quad \mathfrak{u}^{\v}=\varepsilon u_{1}+\varepsilon^2u_{2}+\varepsilon^2u_{R},\quad \theta^{\v}=\theta_{0}+\varepsilon\theta_{1}+\varepsilon^2\theta_{R}.
\end{align}
Then substituting \eqref{expansion} into \eqref{sNS} and \eqref{bd}--\eqref{1.2}, and comparing the order of $\varepsilon$, one has
\begin{align}\label{ge2}
	\left\{
	\begin{aligned}
		&\nabla(\rho_{0}\theta_{0})=0,\quad \int_{\Omega}\rho_{0}\,{\rm d}x=M,\\
		&\operatorname{div}(\rho_{0}u_{1})=0,\\
		&\rho_{0}(u_{1}\cdot \nabla)u_{1}+\nabla P_{2}=\mu\Delta u_{1}+\zeta\nabla\operatorname{div}u_{1},\\
		&\kappa\Delta\theta_{0}=2(\rho_{0}\theta_{0})\operatorname{div}u_{1},\\
		&\theta_{0}=T_{w},\quad(u_{1}\cdot \vec{\iota}_{1}, u_{1}\cdot \vec{\iota}_{2},u_{1}\cdot \vec{n})=h(T_{w})(\partial_{\vec{\iota}_{1}}T_{w},\partial_{\vec{\iota}_{2}}T_{w},0)\quad\text{ on }\partial\Omega,
	\end{aligned}
	\right.
\end{align}
where we have denoted $\zeta=\mu+\lambda$ and $P_{2}:=\rho_{2}\theta_{0}+\rho_{1}\theta_{1}+C_{1}$. Here the constant $C_{1}$ is used to guarantee
\begin{align}\label{C1-1}
\int_{\Omega}\rho_{2}\,{\rm d}x=0\quad \Leftrightarrow\quad
C_{1}=\frac{\int_{\Omega}\frac{1}{\theta_{0}}(P_{2}-\rho_{1}\theta_{1})\,{\rm d}x}{\int_{\Omega}\frac{1}{\theta_{0}}\,{\rm d}x}.
\end{align}
From $\eqref{ge2}_{1}$, we see $\rho_{0}\theta_{0}\equiv P_{0}$ for some constant $P_{0}$. Observing \eqref{1.2}, we chose $P_{0}$ as
$$
P_{0}=\frac{M}{\int_{\Omega}\theta_{0}^{-1}\,{\rm d}x}.
$$
The existence and uniqueness of the strong solution $(\rho_{0},u_{1},\theta_{0},P_{2})$ of \eqref{ge2} is established in Lemma \ref{lem1} below, see also \cite{Esposito-Guo-Marra-Wu-2023-2}.
%$$
%P_{0}=\frac{M}{\int_{\Omega}\frac{1}{\theta_{0}}\,{\rm d}x}.
%$$

Motivated by \cite{Dou-Jiang-Jiang-Yang}, we require $(\rho_{1},u_{2},\theta_{1})$ to satisfy
\begin{align}\label{lge2}
	\left\{
	\begin{aligned}
		&\nabla (\rho_{0}\theta_{1}+\rho_{1}\theta_{0})=0,\quad \int_{\Omega}\rho_{1}\,{\rm d}x=0,\\
		&\operatorname{div}(\rho_{0}u_{2})=-\operatorname{div}(\rho_{1}u_{1}),\\
		&\rho_{0}(u_{1}\cdot \nabla)u_{2}+\rho_{0}(u_{2}\cdot \nabla)u_{1}+\nabla P_{3}=-\rho_{1}(u_{1}\cdot \nabla)u_{1}+\mu\Delta u_{2}+\zeta\nabla\operatorname{div}u_{2},\\
		&\kappa\Delta\theta_{1}=-\theta_{0}(u_{2}\cdot \nabla)\rho_{0}+\rho_{1}(u_{1}\cdot \nabla)\theta_{0}+\rho_{0}(u_{1}\cdot \nabla)\theta_{1}+(\rho_{0}\theta_{1}+\rho_{1}\theta_{0})\operatorname{div}u_{1}\\
		&\qquad\quad\,\, +(\rho_{0}\theta_{0})\operatorname{div}u_{2}-2\theta_{0}(u_{R}\cdot \nabla)\rho_{0},\\
		&u_{2}=0,\quad \theta_{1}=0,\quad \text{on }\partial\Omega,
	\end{aligned}
	\right.
\end{align}
where we have denoted $P_{3}:=\rho_{3}\theta_{0}+\rho_{2}\theta_{1}+C_{2}$. Here the constant $C_{2}$ is used to guarantee
\begin{align}\label{C2-1}
\int_{\Omega}\rho_{3}\,{\rm d}x=0\quad\Leftrightarrow\quad C_{2}=\frac{\int_{\Omega}\frac{1}{\theta_{0}}(P_{3}-\rho_{2}\theta_{1})\,{\rm d}x}{\int_{\Omega}\frac{1}{\theta_{0}}\,{\rm d}x}.
\end{align}
From $\eqref{lge2}_{1}$, we see $\rho_{1}\theta_{0}+\rho_{0}\theta_{1}\equiv P_{1}$ for some constant $P_{1}$. Since $\int_{\Omega}\rho_{1}\,{\rm d}x=0$, the constant $P_{1}$ should be chosen as
$$
P_{1}=\displaystyle\frac{\int_{\Omega}\rho_{0}\theta_{1}\theta_{0}^{-1}\,{\rm d}x}{\int_{\Omega}\theta_{0}^{-1}\,{\rm d}x}.
$$

For the remainder $(\rho_{R},u_{R},\theta_{R})$, it satisfies
\begin{equation}\label{r2}
	\left\{
	\begin{aligned}
		&\operatorname{div}[\rho_{R}(u_{1}+\v u_{2}+\v u_{R})]+\frac{1}{\varepsilon}\operatorname{div}(\rho_{0}u_{R})=-\operatorname{div}[u_{R}(\rho_{1}+\v\rho_{2}+\v^2\rho_{3})]+r_{1},\\
		&\mu\Delta u_{R}+\zeta\nabla \operatorname{div}u_{R}=\frac{1}{\varepsilon}\nabla(\rho_{0}\theta_{R}+\rho_{R}\theta_{0})+\rho_{0}(u_{1}\cdot \nabla u_{R}+u_{R}\cdot \nabla u_{1})+\nabla(\rho_{R}\theta_{1}+\rho_{1}\theta_{R})\\
		&\qquad\qquad\qquad\qquad\quad  +F^{\varepsilon}(\rho_{R},u_{R},\theta_{R})+r_{2},\\
		&\frac{\kappa}{\theta_{0}} \Delta \theta_{R}=\frac{1}{\varepsilon}\operatorname{div}(\rho_{0}u_{R})-\frac{1}{\theta_{0}}\Psi(\nabla (u_{1}+\varepsilon(u_{2}+u_{R})))+\frac{1}{\theta_{0}}(\rho_{0}\theta_{1}+\rho_{1}\theta_{0})\operatorname{div}u_{R}\\
		&\qquad\quad\,\,\,\, +\frac{1}{\theta_{0}}\big[\rho_{0}(u_{R}\cdot \nabla \theta_{1}+u_{1}\cdot \nabla \theta_{R})+\rho_{1}u_{R}\cdot\nabla\theta_{0}+\rho_{R}u_{1}\cdot \nabla \theta_{0}\big]\\
		&\qquad\quad\,\,\,\, +\frac{1}{\theta_{0}}G^{\varepsilon}(\rho_{R},u_{R},\theta_{R})+\frac{1}{\theta_{0}}r_{3}+\frac{1}{\theta_{0}}(\rho_{0}\theta_{R}+\rho_{R}\theta_{0})\operatorname{div}u_{1},\\
		&u_{R}=0,\quad \theta_{R}=0,\quad \text{on }\partial \Omega,\qquad \int_{\Omega}\rho_{R}\,{\rm d}x=0,
	\end{aligned}
	\right.
\end{equation}
where
\begin{align}
	&F^{\varepsilon}(\rho_{R},u_{R},\theta_{R})\nonumber\\
	&=\varepsilon\big[\rho_{0}(u_{2}\cdot \nabla u_{R}+u_{R}\cdot \nabla u_{2}+u_{R}\cdot \nabla u_{R})+\rho_{1}(u_{R}\cdot \nabla u_{1}+ u_{1}\cdot \nabla u_{R})+\rho_{R}u_{1}\cdot \nabla u_{1}\big]\nonumber\\
	&\quad +\varepsilon^2\big[\rho_{1}(u_{R}\cdot \nabla u_{2}+u_{2}\cdot \nabla u_{R}+u_{R}\cdot \nabla u_{R})+\rho_{R}(u_{1}\cdot \nabla u_{R}+u_{R}\cdot \nabla u_{1}+u_{2}\cdot \nabla u_{1}+u_{1}\cdot \nabla u_{2})\nonumber\\
	&\quad +\rho_{2}(u_{R}\cdot \nabla u_{1}+u_{1}\cdot \nabla u_R)\big]+\varepsilon^2\nabla (\rho_{3}\theta_{R}) +\varepsilon^3\big[\rho_{2}(u_{R}\cdot \nabla u_{2}+u_{2}\cdot \nabla u_{R}+u_{R}\cdot \nabla u_{R})\nonumber\\
	&\quad +\rho_{3}(u_{R}\cdot \nabla u_{1}+u_{1}\cdot \nabla u_{R})+\rho_{R}(u_{R}\cdot \nabla u_{2}+u_{2}\cdot \nabla u_{R}+u_{2}\cdot \nabla u_{2}+u_{R}\cdot \nabla u_{R})\big]\nonumber\\
	&\quad+\varepsilon^{4}\big[\rho_{3}(u_{R}\cdot \nabla u_{2}+u_2\cdot u_{R}+u_{R}\cdot \nabla u_{R})\big]+\varepsilon\nabla (\rho_{2}\theta_{R}+\rho_{R}\theta_{R}),\label{F}\\
	&G^{\varepsilon}(\rho_{R},u_{R},\theta_{R})\nonumber\\
	&= \varepsilon\big[\rho_{0}(u_{2}\cdot\nabla \theta_{R}+u_{R}\cdot \nabla \theta_{R})+\rho_{2}u_{R}\cdot \nabla\theta_{0}+\rho_{R}(u_{2}\cdot \nabla\theta_{0}+u_{R}\cdot \nabla \theta_{0}+u_{1}\cdot \nabla \theta_{1})\nonumber\\
	&\quad +\rho_{1}(u_{R}\cdot \nabla \theta_{1}+u_{1}\cdot \nabla \theta_{R})\big]
	+\varepsilon^2\big[\rho_{1}(u_{2}\cdot\nabla \theta_{R}+u_{R}\cdot \nabla \theta_{R})+\rho_{R}(u_{2}\cdot \nabla\theta_{1}+u_{R}\cdot \nabla \theta_{1}+u_{1}\cdot \nabla \theta_{R})\nonumber\\
	&\quad +\rho_{2}(u_{R}\cdot \nabla\theta_{1}+u_{1}\cdot \nabla \theta_{R})+\rho_{3}u_{R}\cdot \nabla \theta_{0}\big]+\varepsilon^3\big[\rho_{2}(u_{2}\cdot\nabla \theta_{R}+u_{R}\cdot \nabla \theta_{R})+\rho_{3}(u_{1}\cdot \nabla \theta_{R}+u_{R}\cdot \nabla \theta_{1})\nonumber\\
	&\quad+\rho_{R}(u_{2}\cdot\nabla \theta_{R}+u_{R}\cdot \nabla \theta_{R})\big]+\varepsilon^{4}\big[\rho_{3}(u_{2}\cdot \nabla \theta_{R}+u_{R}\cdot \nabla \theta_{R})\big] +\varepsilon\big[(\rho_{0}\theta_{R}+\rho_{R}\theta_{0})\operatorname{div}u_{2}\nonumber\\
	&\quad+(\rho_{0}\theta_{R}+\rho_{R}\theta_{0}+\rho_{1}\theta_{1}+\rho_{2}\theta_{0})\operatorname{div}u_{R} +(\rho_{1}\theta_{R}+\rho_{R}\theta_{1})\operatorname{div}u_{1}\big] +\varepsilon^2\big[(\rho_{1}\theta_{R}+\rho_{R}\theta_{1})\operatorname{div}u_{2}\nonumber\\
	&\quad+(\rho_{2}\theta_{R}+\rho_{R}\theta_{R})\operatorname{div}u_{1}+(\rho_{2}\theta_{1}+\rho_{1}\theta_{R}+\rho_{R}\theta_{1}+\rho_{3}\theta_{0})\operatorname{div}u_{R}  \big]+\varepsilon^3\big[(\rho_{2}\theta_{R}+\rho_{R}\theta_{R})\operatorname{div}u_{2}\nonumber\\
	&\quad +\rho_{3}\theta_{R}\operatorname{div}u_{1}+(\rho_{2}\theta_{R}+\rho_{R}\theta_{R}+\rho_{3}\theta_{1})\operatorname{div}u_{R}\big ]+\varepsilon^{4}\big[\rho_{3}\theta_{R}\operatorname{div}u_2+\rho_{3}\theta_{R}\operatorname{div}u_R\big],\label{G}
\end{align}
and
\begin{align}
	&r_{1}=-\operatorname{div}(\rho_{1}u_{2}+\rho_{2}u_{1})-\varepsilon\operatorname{div}(\rho_{2}u_{2}+\rho_{3}u_{1})-\varepsilon^2\operatorname{div}(\rho_{3}u_{2}),\label{r1-1}\\
	&r_{2}=\rho_{0}(u_{2}\cdot \nabla u_{1}+u_{1}\cdot \nabla u_{2})+\varepsilon\big[\rho_{0}u_{2}\cdot \nabla u_{2}+\rho_{1}(u_{2}\cdot \nabla u_{1}+u_{1}\cdot \nabla u_{2})+\rho_{2}u_{1}\cdot \nabla u_{1}\big]\nonumber\\
	&\qquad +\varepsilon^2\big[\rho_{1}u_{2}\cdot \nabla u_{2}+\rho_{2}(u_{1}\cdot \nabla u_{2}+u_{2}\cdot \nabla u_{1})+\rho_{3}u_{1}\cdot \nabla u_{1}\big]+\varepsilon^{4}\rho_{3}u_{2}\cdot \nabla u_{2}\nonumber\\
	&\qquad+\varepsilon^3\big[\rho_{2}u_{2}\cdot \nabla u_{2} +\rho_{3}(u_{1}\cdot \nabla u_{2}+u_{2}\cdot \nabla u_{1})\big]+\varepsilon\nabla (\rho_{3}\theta_{1}),\label{r2-1}\\
	&r_{3}=\big[\rho_{0}(u_{2}\cdot \nabla \theta_{1})+\rho_{1}(u_{2}\cdot \nabla\theta_{0}+u_{1}\cdot \nabla\theta_{1})+\rho_{2}u_{1}\cdot \nabla\theta_{0}\big]\nonumber\\
	&\qquad +\varepsilon\big[\rho_{1}(u_{2}\cdot \nabla\theta_{1})+\rho_{2}(u_{2}\cdot \nabla\theta_{0}+u_{1}\cdot \nabla\theta_{1})+\rho_{3}u_{1}\cdot \nabla \theta_{0}\big]\nonumber\\
	&\qquad +\varepsilon^2\big[\rho_{2}u_2\cdot \nabla \theta_{1}+\rho_{3}(u_{2}\cdot \nabla \theta_{0}+u_{1}\cdot \nabla \theta_{1})\big]+\varepsilon^3\rho_{3}u_{2}\cdot \nabla \theta_{1}\nonumber\\
	&\qquad +(\rho_{0}\theta_{1}+\rho_{1}\theta_{0})\operatorname{div}u_{2}+(\rho_{2}\theta_{0}+\rho_{1}\theta_{1})\operatorname{div}u_{1}\nonumber\\
	&\qquad +\varepsilon\big[(\rho_{2}\theta_{0}+\rho_{1}\theta_{1})\operatorname{div}u_{2}+(\rho_{2}\theta_{1}+\rho_{3}\theta_{0})\operatorname{div}u_{1}\big]\nonumber\\
	&\qquad +\varepsilon^2[(\rho_{2}\theta_{1}+\rho_{3}\theta_{0})\operatorname{div}u_{2}+\rho_{3}\theta_{1}\operatorname{div}u_{1}]+\varepsilon^3\rho_{3}\theta_{1}\operatorname{div}u_{2}.\label{r3}
\end{align}

Since we put an extra term $-2\theta_{0}(u_{R}\cdot \nabla)\rho_{0}$ into \eqref{lge2} as a source term, the singular terms (the terms involving $\var^{-1}$) are skew-symmetric in the remainder equations \eqref{r2}. Though it will make \eqref{lge2} and \eqref{r2} slightly coupled, the skew-symmetric property can  eliminate the singular terms through integrating by parts in the energy estimates. This is a key point to derive the {\it a priori} uniform estimates for $(\rho_{R},u_{R},\theta_{R})$.

\medskip

Now, we state the main result of present paper:
\begin{theorem}\label{thm1.1}
	Let $\Omega$ be a $C^{4}-$smooth bounded domain and $\nabla T_{w}\in H^{4}(\partial \Omega)$. Let $(\rho_{0},u_{1},\theta_{0})\in H^5\times H^4\times H^5$ be the strong solution of \eqref{ge2} established in Lemma \ref{lem1} below. There exists a constant $c_{0}>0$
	%, $c_{1}>0$
	%and $\varepsilon_{0}>0$
	such that if $|\nabla T_{w}|_{H^4}\leq c_{0}$, then for any given Mach number $0<\v<1$,
	%and $\varepsilon\in (0,\varepsilon_{0})$.
	there exists a unique strong solution $(\rho^{\v}, \mathfrak{u}^{\v},\theta^{\v})\in H^2\times H^3\times H^3$ in the form of expansion \eqref{expansion} to the boundary problem \eqref{sNS}, \eqref{bd} and \eqref{1.2}, satisfying
%\begin{align}\label{c0}
%	\|\rho^{\v}\|_{H^2}+\|\mathfrak{u}^{\v}\|_{H^3}+\|\theta^{\v}\|_{H^3}\leq c_{1},
%\end{align}
%and
\begin{align}\label{c}
	 \|\frac{\mathfrak{u}^{\v}}{\varepsilon}-u_{1}\|_{H^2}
	 %+\|\nabla^3(\mathfrak{u}^{\v}-\v u_{1})\|_{L^2}
	 +\|\rho^{\v}-\rho_{0}\|_{H^2}+\|\theta^{\v}-\theta_{0}\|_{H^3}\leq C\varepsilon,\quad \|\nabla P^{\v}\|_{H^1}\leq C\varepsilon^2,
\end{align}
	where the positive constants $c_{0}$ and $C$ are  independent of $\varepsilon$.
\end{theorem}

%\begin{remark}
%	The uniqueness of the strong solution $(\rho^{\v}, \mathfrak{u}^{\v},\theta^{\v})\in H^2\times H^3\times H^3$ for \eqref{sNS}, \eqref{bd} and \eqref{1.2} is not clear. However, the uniqueness of strong solution is true in the form of expansion \eqref{expansion}.
%\end{remark}

\begin{remark}
	%It is worth pointing out that $c_{0}$ is independent of $\varepsilon$ in Theorem \ref{thm1.1}, which
In Theorem \ref{thm1.1},  the variation of the wall temperature can be large compared with the Mach number $\v$.
\end{remark}

%\begin{remark}
%There is no requirements on the smallness of the Mach number $\v$ for the existence of strong solution to the boundary problem \eqref{sNS}, \eqref{bd} and \eqref{1.2}, in the form of expansion \eqref{expansion}. %In fact, the Mach number $\v$ can be any finite number in this aspect.
%\end{remark}

\begin{remark}
	The Dirichlet boundary condition for velocity in \eqref{bd} can be generalized to
	$$
	(\mathfrak{u}^{\v}\cdot \vec{\iota}_{1},\mathfrak{u}^{\v}\cdot \vec{\iota}_{2},\mathfrak{u}^{\v}\cdot \vec{n})=\varepsilon (h_{1}(x),h_{2}(x),0),
	$$
	where $h_{i}(x)$ can be any $H^4-$ function on $\partial\Omega$ with $o(1)-$norm.
\end{remark}

\begin{remark}
It follows from \eqref{c} and the Sobolev embedding that
\begin{equation}\label{c1}
	\|\frac{\mathfrak{u}^{\v}}{\varepsilon}-u_{1}\|_{L^{\infty}}+\|\rho^{\v}-\rho_{0}\|_{L^{\infty}}+\|\theta^{\v}-\theta_{0}\|_{L^{\infty}}\leq C\varepsilon\rightarrow 0.
\end{equation}
Compared with \cite{Esposito-Guo-Marra-Wu-2023-1}, in which the convergence is proved in $L^2$ for general smooth bounded domain, \eqref{c1} provides a convergence in $L^{\infty}$.
\end{remark}

\subsection{The main ingredients in proof}\label{kp} Now, we briefly introduce the main ingredients in our proof. First, the expansion \eqref{expansion} is carefully designed such that $\rho_{2}$, $\rho_{3}$ can be well-defined under the Boussinesq relation and the constraints: $\int_{\Omega}\rho_{i}\,{\rm d}x=0\,\,(i=2,3)$. In fact, if one makes a symmetric expansion for density and temperature, that is, let
$$
\rho^{\v}=\rho_{0}+\varepsilon\rho_{1}+\varepsilon^2\rho_{2}+\varepsilon^2\rho_{R},\,\,\theta^{\v}=\theta_{0}+\varepsilon\theta_{1}+\varepsilon^2\theta_{2}+\varepsilon\theta_{R}.
$$
Then the Boussinesq relations in \eqref{ge2} and \eqref{lge2} may lead to
$$
P_{2}=\rho_{0}\theta_{2}+\theta_{0}\rho_{2}+\rho_{1}\theta_{1},\quad P_{3}=\rho_{1}\theta_{2}+\rho_{2}\theta_{1}.
$$
To define $\rho_{2}$ and $\theta_{2}$, one has to require $\rho_{0}\theta_{1}-\theta_{0}\rho_{1}\neq 0$. However, this constraint is hard to verify since $\rho_{1}$ and $\theta_{1}$ are already solved from \eqref{lge2}. Based on this observation, we design a non-symmetric expansion \eqref{expansion} for density and temperature such that there is a no extra restriction to define $\rho_{2}$ and $\rho_{3}$ satisfying the constraints of zero average under the Boussinesq relation once we solved $P_{2}$ and $P_{3}$.

Next, noting that the systems \eqref{lge2} and \eqref{r2} are slightly coupled, we shall take the following five steps to establish the existence of strong solutions of \eqref{lge2} and \eqref{r2}.
%which is more complicated than those in (see also \cite{Choe-Jin,Dou-Jiang-Jiang-Yang,Li-Liao-2019}):

1. Motivated by \cite{Choe-Jin,Dou-Jiang-Jiang-Yang,Li-Liao-2019}, we first replace $u_{R}$ by a given function $\tilde{u}_{R}\in \mathcal{K}$
in \eqref{lge2} to get a new linearized system \eqref{3.1}, where $\mathcal{K}$ is defined as
$$
\mathcal{K}:=\{f\in H^3(\Omega)\cap H_{0}^{1}(\Omega):\,\|f\|_{\mathcal{K}}=:\|f\|_{H^2(\Omega)}+\varepsilon\|\nabla^3 f\|_{L^2(\Omega)}<\infty\}.
$$
We point out that the $\v-$dependent functional space $\mathcal{K}$ is introduced to deal with the uniform estimates of the highest normal derivatives of $u_{R}$ since its boundary conditions are not available. By using a similar fixed-point argument as in the proof of the existence of \eqref{ge2}, we prove that there exists a unique strong solution $(\rho_{1},u_{2},\theta_{1},P_{3})$ of \eqref{3.1}, see Lemma \ref{lem2.1} for details.
%By using a similar fixed-point argument as in the proof of the existence of \eqref{ge2}, we first prove there exists a unique strong solution $(\rho_{1},u_{2},\theta_{1},P_{3})$ to \eqref{lge2} with

2. For any given pair $(\tilde{u}_{R},\tilde{\theta}_{R})$ with $\tilde{u}_{R}\in \mathcal{K}$ and $\tilde{\theta}_{R}\in H_{0}^1(\Omega)\cap H^3(\Omega)$, we consider a linearized system of \eqref{r2} (that is, \eqref{3.10} below), where $(\rho_{0},u_{1},\theta_{0},P_{2})$ is the strong solution of \eqref{ge2} established in Lemma \ref{lem1},  $(\rho_{1},u_{2},\theta_{1},P_{3})$ is the strong solution of \eqref{3.1} established in Lemma \ref{lem2.1}, and $(\rho_{2},\rho_{3})$ is defined in \eqref{3.7-4} through the Boussinesq relations. Motivated by \cite{Novotny-Padula-1994,Novotny-Straskraba,Valli-1987}, with the help of the effective viscous flux and the Helmholtz decomposition, we decompose the system \eqref{3.10} into two Laplace equations with Dirichlet and Neumann boundary conditions respectively, a classical Stokes problem and a steady transport equation. Then based on this decomposition, we design an elaborate approximate scheme to establish the existence the strong solution $(\rho_{R},u_{R},\theta_{R})\in H^2\times \mathcal{K}\times (H_{0}^1\cap H^3)$ of \eqref{3.10} for a specific viscous coefficients pair $(\mu_{0},\zeta_{0})$ with $\frac{\zeta_{0}}{\mu_{0}}\gg 1$ by making good use of the classical elliptic estimates, the Stokes estimates and the properties of steady transport equation, see Lemma \ref{lem7.1} for details.

%we design an elaborate approximate scheme to establish the existence the strong solution $(\rho_{R},u_{R},\theta_{R})\in H^2\times \mathcal{K}\times (H_{0}^1\cap H^3)$ of \eqref{3.10} with a specific viscous coefficients pair $(\mu_{0},\zeta_{0})$ with $\frac{\zeta_{0}}{\mu_{0}}\gg 1$ by making good use of the Helmholtz decomposition, the effective viscous flux and the property of steady transport equation, see Lemma \ref{lem7.1} for details.
%Motivated by \cite{Dou-Jiang-Jiang-Yang}, by using a vanishing artificial viscosity approach and Leray-Schauder fixed point theorem, we can get the existence and uniqueness of the weak solution $(\rho_{R},u_{R},\theta_{R})\in L^2\times H_{0}^{1}(\Omega)\times H_{0}^1(\Omega)$ of \eqref{3.10} under the assumption $\|\tilde{u}_{R}\|_{\mathcal{K}}+\|\tilde{\theta}_{R}\|_{H^3}$ is small; see Lemma \ref{lem2.2}.
%Then we improve the regularity of $(\rho_{R},u_{R},\theta_{R})$ to $H^2\times H^3\times H^3$ via Hlemohotz decomposition and the regularity theories of elliptic equations, Stokes problem and steady transport equations such that $(\rho_{R},u_{R},\theta_{R})$ is a unique strong solution of \eqref{3.10}, see Lemma \ref{lem2.3}.

3. To generalize the existence of strong solution of \eqref{3.10} in Lemma \ref{lem7.1} to any given viscous coefficient $(\mu, \zeta)$, as in \cite{Valli-1987}, we need to show that the strong solution $(\rho_{R},u_{R},\theta_{R})$ of
\eqref{3.10} enjoys {\it a priori} uniform estimate with respect to $\var$ and continuous dependence on $\mu$ and $\zeta$. The uniform lower order estimates of $(\rho_{R},u_{R},\theta_{R})$ can be directly obtained by the classical energy method due to the skew-symmetry of the singular terms (the terms involving $\var^{-1}$) in \eqref{3.10}.
%As indicated in \cite{Valli-1987}, the regularity of the weak solution of the linearized system \eqref{3.10} cannot be improved to $(\rho_{R},u_{R},\theta_{R})\in H^2\times \mathcal{K}\times (H^3\cap H_{0}^1)$ through the approximate scheme ({\it i.e.} \eqref{3.14} below) used in Step 2. This is because the term $u_{R}\cdot \nabla \rho_{0}$ in $\eqref{3.14}_{1}$ makes \eqref{3.14} lack of elliptic regularity. To overcome this difficulty, we make full use of the effective viscous flux, good property of steady transport equation and then follow similar approximate scheme as in \cite{Novotny-Padula-1994,Novotny-Straskraba,Valli-1987} to construct strong solution of \eqref{3.14}. Before this, we have to show that
%To establish a fixed point framework, we further show that the strong
%the strong solution $(\rho_{R},u_{R},\theta_{R})$ of
%\eqref{3.10} enjoys a uniform higher order {\it a priori} estimate with respect to $\var$.

To get the uniform higher order estimates, we regard $\eqref{3.10}_2$ as an inhomogeneous Stokes problem, and by using the classical Stokes estimate \cite[Theorems IV 5.2 and IV 5.8]{Boyer-Fabrie}, we get
\begin{align}\label{S1}
\|u_{R}\|_{H^2}+\|\frac{\nabla (\rho_{0}\theta_{R}+\theta_{0}\rho_{R})}{\varepsilon}\|_{L^2}&\leq C(\|\rho_{R}\|_{H^1}\|\theta_{1}\|_{H^2}+\|\rho_{1}\|_{H^2}\|\theta_{R}\|_{H^1})+C\|\nabla \operatorname{div}u_{R}\|_{L^2}\nonumber\\
&\quad +C\|r_{2}\|_{L^2}+\text{lower order derivatives of $(\rho_{R},u_{R},\theta_{R})$},
\end{align}
and
\begin{align}\label{S2}
	\|u_{R}\|_{H^3}+\|\frac{\nabla (\rho_{0}\theta_{R}+\theta_{0}\rho_{R})}{\varepsilon}\|_{H^1}&\leq C(\|\rho_{R}\|_{H^2}\|\theta_{1}\|_{H^2}+\|\rho_{1}\|_{H^2}\|\theta_{R}\|_{H^2})+C\|\nabla^2 \operatorname{div}u_{R}\|_{L^2}\nonumber\\
	&\quad +C\|r_{2}\|_{H^1}+\text{lower order derivatives of $(\rho_{R},u_{R},\theta_{R})$}.
\end{align}
Therefore, to close these estimates, we have to control $\|\theta_{R}\|_{H^2}$, $\|\nabla \operatorname{div}u_{R}\|_{L^2}$ and $\|\nabla^2\operatorname{div}u_{R}\|_{L^2}$. We divide these estimates into two parts : the interior part and the part near the boundary.

The interior part can be bounded by the classical energy estimate since the singular terms (the terms involving $\var^{-1}$) in \eqref{3.10} are skew-symmetric, and will cancel each other out through integrating by parts. For the part near the boundary, we shall flatten the boundary. Noting the zero Dirichlet boundary condition of $(u_{R},\theta_{R})$, we can get the estimates of tangential derivatives for the part near the boundary via energy estimates as in the interior part.

The main difficulty is to estimate the terms involving the normal derivatives due to the lack of boundary conditions of the higher normal derivatives of $(u_{R},\theta_{R})$. To overcome this difficulty, we fully exploit the structure of $\eqref{3.10}$ to see
\begin{align}\label{D3r}
	\frac{1}{\v}\operatorname{div}(\rho_{0}u_{R})=-\operatorname{div}[\rho_{R}(u_{1}+\v(u_{2}+\tilde{u}_{R}))]+\cdots,
\end{align}
\begin{align}\label{D3u}
(\mu+\zeta)D_{3}\operatorname{div}u_{R}-\frac{D_{3}(\theta_{0}\rho_{R}+\rho_{0}\theta_{R})}{\varepsilon}=-\mu[\Delta u_{R}^{3}-D_{3}\operatorname{div}u_{R}]+\cdots,
\end{align}
and
\begin{align}\label{D3t}
\frac{1}{\theta_{0}}D_{33}^2\theta_{R}-\frac{\operatorname{div}(\rho_{0}u_{R})}{\varepsilon}=-\frac{1}{\theta_{0}}(D_{11}^2\theta_{R}+D_{22}^2\theta_{R})+\cdots.
\end{align}
Then, in order to get $\|D_{3}\operatorname{div}u_{R}\|_{L^2}$ and $\|D_{33}^2\theta_{R}\|_{L^2}$ in the vicinity of boundary, denoting $\chi_{1}$ by the smooth cut-off function satisfying $\chi_{1}\equiv 1$ in the vicinity of boundary, we multiply \eqref{D3r}, \eqref{D3u} and \eqref{D3t} by
$-\chi_{1}^2D_{33}^2(\theta_{0}\rho_{R})$, $\chi_{1}^2D_{3}\operatorname{div}(\rho_{0}u_{R})$ and $\chi_{1}^2D_{33}^2(\rho_{0}\theta_{R})$ respectively to get
\begin{align*}
&\int_{\Omega}\chi_{1}^2\rho_{0}|D_{3}\operatorname{div}u_{R}|^2\,{\rm d}x+\int_{\Omega}\chi_{1}^2\frac{\rho_{0}}{\theta_{0}}|D_{33}^2\theta_{R}|^2\,{\rm d}x\\
&\leq -\mu\int_{\Omega}\chi_{1}^2\big\{[\Delta u_{R}^{3}-D_{3}\operatorname{div}u_{R}]D_{3}\operatorname{div}(\rho_{0}u_{R})+\frac{1}{\theta_{0}}[D_{11}^2\theta_{R}+D_{22}^2\theta_{R}]D_{33}^2\theta_{R}\big\}\,{\rm d}x\\
&\quad %+\frac{1}{\varepsilon}\int_{\Omega}D_{3}(\theta_{0}\rho_{R})D_{3}\operatorname{div}(\rho_{0}u_{R})\,{\rm d}x
+\frac{1}{\varepsilon}\int_{\Omega}\chi_{1}^2\big[D_{3}(\rho_{0}\theta_{R}+\theta_{0}\rho_{R})D_{3}\operatorname{div}(\rho_{0}u_{R})+D_{33}^2(\rho_{0}\theta_{R}+\theta_{0}\rho_{R})\operatorname{div}(\rho_{0}u_{R})\big]\,{\rm d}x\\
&\quad +\int_{\Omega}\chi_{1}^2D_{33}^2(\theta_{0}\rho_{R})\operatorname{div}[\rho_{R}(u_{1}+\v(u_{2}+\tilde{u}_{R}))]\,{\rm d}x+\cdots\\
&=:J_{1}+J_{2}+J_{3}+\cdots.
\end{align*}
For $J_{1}$, noting that $\Delta u_{R}^{3}-D_{3}\operatorname{div}u_{R}$ does not include the term $D_{33}^2u_{R}^{3}$, by the H\"{o}lder inequality, it can be bounded as
$$
|J_{1}|\leq \frac{1}{8}\int_{\Omega}\chi_{1}^2\rho_{0}|D_{3}\operatorname{div}u_{R}|^2\,{\rm d}x+\frac{1}{8}\int_{\Omega}\chi_{1}^2\frac{\rho_{0}}{\theta_{0}}|D_{33}^2\theta_{R}|^2\,{\rm d}x+\text{tangential derivatives of $(\nabla u_{R},\nabla\theta_{R})$}.
$$
For $J_{3}$, by the Sobolev embedding, one directly has
$$
|J_{3}|\leq C[\|u_{1}\|_{H^3}+\v(\|u_{2}\|_{H^3}+\|\tilde{u}_{R}\|_{H^3})]\|\rho_{R}\|_{H^2}.
$$
Notably, $[\|u_{1}\|_{H^3}+\v(\|u_{2}\|_{H^3}+\|\tilde{u}_{R}\|_{H^3})]$ is a small coefficient.
%by the tangential derivatives of $(\nabla u_{R},\nabla\theta_{R})$.
%For $J_{2}$, we can apply $D_{3}$ to $\eqref{3.10}_{1}$ and multiply the resultant equation by $D_{3}(\rho_{R}\theta_{0})$ to cancel it.
For $J_{2}$, by integrating by parts, it leads to
\begin{align*}
J_{2}&=-\frac{1}{\v}\int_{\Omega}\chi_{1}D_{3}\chi_{1}D_{3}(\rho_{0}\theta_{R}+\theta_{0}\rho_{R})\operatorname{div}(\rho_{0}u_{R})\,{\rm d}x+\frac{1}{\varepsilon}\chi_{1}^2\int_{\partial\Omega}D_{3}(\rho_{0}\theta_{R}+\theta_{0}\rho_{R})\operatorname{div}(\rho_{0}u_{R})\,{\rm d}x\\
&=J_{2,1}+J_{2,2}.
\end{align*}
For $J_{2,1}$, using the H\"{o}lder inequality, for any small $\tau>0$, one has
\begin{align*}
	|J_{2,1}|\leq \tau\|\frac{\nabla (\rho_{0}\theta_{R}+\theta_{0}\rho_{R})}{\v}\|_{L^2}+C_{\tau}\|u_{R}\|_{H^1}^2.
\end{align*}
For the boundary term $J_{2,2}$, noting \eqref{D3r}, we can represent $\frac{1}{\varepsilon}\operatorname{div}(\rho_{0}u_{R})$ through other non-singular terms, and then we can control $J_{2,2}$ in terms of $\|(\rho_{R},\theta_{R})\|_{H^2}$ with some small coefficients. Therefore, combining the above estimates with \eqref{S1}, we can control $\|(u_{R},\theta_{R})\|_{H^2}$ in terms of $\|\rho_{R}\|_{H^2}$ with some small coefficients; see Lemma \ref{lem5.0-1} for details.

However, it is still hard to apply the above arguments to control $\|\nabla^2\operatorname{div}u_{R}\|_{L^2}$ since it requires $\|\rho_{R}\|_{H^3}$ to control the singular terms induced by the boundary terms in the process of integrating by parts, which cannot be closed by using \eqref{S1}--\eqref{S2}. To circumvent the difficulty, we resort to controlling $\varepsilon\|\nabla^2\operatorname{div}u_{R}\|_{L^2}$. That is the reason why we introduce the $\v-$dependent functional space $\mathcal{K}$. We rewrite \eqref{S2} as
\begin{align}\label{S2-1}
	\varepsilon\|u_{R}\|_{H^3}+\|\rho_{R}\|_{H^2}\leq C(\|\theta_{R}\|_{H^2}+\|\rho_{R}\|_{H^1})+C\varepsilon(\|u_{R}\|_{H^1}+\|r_{2}\|_{H^1}+\|\nabla^2\operatorname{div}u_{R}\|_{L^2}).
\end{align}
Then we apply an elaborate $\varepsilon-$dependent energy method to control $\varepsilon\|(\nabla^2\operatorname{div}u_{R},\nabla^3\theta_{R})\|_{L^2}$ in terms of $\|\rho_{R}\|_{H^2}$, and hence close the estimates of $\|\rho_{R}\|_{H^2}+\|(u_{R},\theta_{R})\|_{\mathcal{K}}$ via \eqref{S1} and \eqref{S2-1}. Precisely, for the singular terms involving $\theta_{R}$ and $u_{R}$:
$$
\var\int_{\Omega}D_{33}^2(\rho_{0}\theta_{R})D_{33}^2\operatorname{div}(\rho_{0}u_{R})\,{\rm d}x\quad\text{and}\quad  \var\int_{\Omega}D_{3}\operatorname{div}(\rho_{0}u_{R})D_{333}^3(\rho_{0}\theta_{R})\,{\rm d}x,
$$
we control them via $\|(u_{R},\theta_{R})\|_{H^2}$ by using the H\"{o}lder inequality directly, while for the singular terms involving $\rho_{R}$ and $u_{R}$:
$$
\var\int_{\Omega}D_{33}^2(\theta_{0}\rho_{R})D_{33}^2(\operatorname{div}u_{R})\,{\rm d}x,
$$
we apply $D_{33}^2$ to $\eqref{3.10}_{1}$ and multiply the resultant equation by $D_{33}^2(\rho_{R}\theta_{0})$ to find a cancellation between bad terms; see Lemmas \ref{lem5.1}--\ref{lem5.3} for details.

4. With the help of the above uniform estimates of $\|\rho_{R}\|_{H^2}+\|(u_{R},\theta_{R})\|_{\mathcal{K}}$, we can further get the uniform estimates of $\|\theta_{R}\|_{H^3}$ by using $\eqref{3.10}_{1}$ to remove the singular terms $\frac{1}{\v}\operatorname{div}(\rho_{0}u_{R})$ in $\eqref{3.10}_{3}$ and the standard elliptic regularity. Noting the {\it a priori} uniform estimates depend continuously on the viscous coefficients $(\mu,\zeta)$, we can apply an elaborate continuity method as in \cite{Valli-1987} to generalize this specific pair of $(\mu_{0},\zeta_{0})$ in Lemma \ref{lem7.1} to any given viscous coefficients pair; see Lemma \ref{lem7.2} for details.

5. For the existence of strong solution $(\rho_{R},u_{R},\theta_{R})\in H^2\times \mathcal{K}\times (H_{0}^1\times H^3)$ to the nonlinear coupled problem \eqref{lge2}--\eqref{r2}, we shall apply the Tikhononv's fixed point theorem to show the mapping $(\tilde{u}_{R},\tilde{\theta}_{R})\mapsto (u_{R},\theta_{R})$ has a unique fixed point in the space $\mathbf{K}_{A}=:\{(u_{R},\theta_{R})\in \mathcal{K}\times (H_{0}^1\cap H^3):\,\,\|u_{R}\|_{\mathcal{K}}+\|\theta_{R}\|_{H^3}\leq A\}$ provided that $A$ is small enough, see Lemma \ref{lem6.1} for details.
%Then we obtain a unique solution $(\rho_{1},u_{2},\theta_{1},P_{3},\rho_{R},u_{R},\theta_{R})$ of the nonlinear coupled systems \eqref{lge2} and \eqref{r2}.

%Next, to make the singular terms (the terms involved with $\frac{1}{\varepsilon}$) skew-symmetric in the remainder equation \ref{r2}, we put an additional term $-2\theta_{0}(u_{R}\cdot \nabla)\rho_{0}$ into \eqref{lge2} as an additional source term. This makes \eqref{lge2} and \eqref{r2} slightly coupled, and the proof of existence theory of \eqref{lge2} and \eqref{r2} more complicated. However, the skew-symmetric property such that we can remove the singular terms through integrating by parts during the process of the energy estimates. This is very to important to derive the lower order estimates for $(\rho_{R},u_{R},\theta_{R})$.

%To derive the estimates of third order derivatives of $(u_{R},\theta_{R})$, the integrating by parts failed due to the lack of the trace of second order normal derivatives of $(u_{R},\theta_{R})$, and then the singular terms does not vanish by applying the energy estimates directly. To overcome this difficulty, we turn to estimate $\varepsilon (\|\nabla^3 u_{R}\|_{L^2}+\|\nabla^3 \theta_{R}\|_{H^3})$, which is closed with the help of $\|\rho_{R}\|_{H^2}$
\smallskip

The rest of present paper is organized as follows. In Section 2, we  establish the existence of the strong solution to \eqref{ge2} via the Banach fixed point theorem. In Section 3, the existence of strong solution to \eqref{lge2} with replacing $u_{R}$ by $\tilde{u}_{R}$ is similarly proved. Then we also establish the existence of the strong solution $(\rho_{R},u_{R},\theta_{R})\in H^2\times \mathcal{K}\times (H_{0}^1\cap H^3)$ of \eqref{3.10} with a specific viscous coefficients pair $(\mu_{0},\zeta_{0})$.
%we establish the existence and uniqueness of the weak solution to a linearized system of \eqref{r2} (the system \eqref{3.10}) through a vanishing artificial viscosity approach and Leray-Schauder fixed point theorem.
In Section 4, we will establish {\it a priori} uniform estimates of the strong solution of \eqref{3.10}, including the lower and higher order estimates. In Section 5, with the help of {\it a priori} uniform estimates established in Section 4, we apply an elaborate continuity argument to extend the existence result in the case of specific viscous coefficients pair $(\mu_{0},\zeta_{0})$ to the case of any given viscous coefficients pair $(\mu,\zeta)$. In Section 6, we shall use the Tikhonov fixed point theorem to establish the existence and uniqueness of the strong solution of the coupled nonlinear systems \eqref{lge2} and \eqref{r2}, and then complete the proof of Theorem \ref{thm1.1}.

\smallskip
{\bf Notations:} Throughout this paper, we denote $ L^{p}(\Omega), W^{k, p}(\Omega)$ and $W_{0}^{k,p}(\Omega)$ as the standard Sobolev spaces on domain $\Omega\subset \R^3$ and $p\in [1,\infty]$. In particular, we denote $H^{k}(\Omega)=:W^{k,2}(\Omega)$ and $H_{0}^{k}(\Omega)=:W_{0}^{k,2}(\Omega)$. $C_{c}^{k}(\R^3)\,\,(0\leq k\leq \infty)$ represents the space of
continuously differentiable functions up to the $k$th order
with compact support over $\R^3$. $a\lesssim b$ means there exists a constant $C$ independent of $\varepsilon$ such that $a\leq Cb$.

\section{Well-posedness of the Limiting System \eqref{ge2}}
Motivated by \cite[Theorem 3.3]{Esposito-Guo-Marra-Wu-2023-2}, we have following well-posedness theory for the limiting system \eqref{ge2}.
\begin{lemma}\label{lem1}
	Let $\nabla T_{w}\in H^{m}(\partial \Omega)$ with $m\geq 4$. For any given positive constant $\delta_{0}\ll 1$, there exists a positive constant $\delta_{1}$, depending on $\delta_{0}$, such that if $|\nabla T_{w}|_{H^m}\leq \delta_{1}$, then there exists a unique solution $(\rho_{0},u_{1},\theta_{0},P_{2})$ 
	${\rm{(}}\int_{\Omega}P_{2}(x)\,{\rm d}x=0\rm{)}$ to \eqref{ge2} satisfying
	$$
 \|u_{1}\|_{H^m}+\|P_{2}\|_{H^{m-1}}+\|\theta_{0}-1\|_{H^{m+1}}+\|\rho_{0}-P_{0}\|_{H^{m+1}}\lesssim \delta_{0}.
	$$
\end{lemma}

\noindent\textbf{Proof.} Noting $\eqref{ge2}_{1}$, one has that $\rho_{0}\theta_{0}=P_{0}$ with $P_{0}=M(\int_{\Omega}\theta_{0}^{-1}\,{\rm }dx)^{-1}$. Denoting $v_{1}=\rho_{0}u_{1}$, one obtains directly from $\eqref{ge2}_{2}$ and $\eqref{ge2}_{5}$ that
\begin{align}\label{2.1}
	\operatorname{div}v_{1}=0,\quad (v_{1}\cdot \vec{\iota}_{1},v_{1}\cdot \vec{\iota}_2,v_{1}\cdot \vec{n})=P_{0}\frac{h(T_{w})}{T_{w}}(\partial_{\vec{\iota}_{1}}T_{w},\partial_{\vec{\iota}_{2}}T_{w},0).
\end{align}
Furthermore, a direct calculation shows that
\begin{align}\label{2.2}
&\operatorname{div}u_{1}=\frac{1}{P_{0}}v_{1}\cdot \nabla\theta_{0},\quad \rho_{0}(u_{1}\cdot \nabla)u_{1}=\frac{1}{P_{0}}(v_{1}\cdot \nabla \theta_{0})v_{1}+\frac{\theta_{0}}{P_{0}}(v_{1}\cdot \nabla)v_{1}\nonumber\\
& \Delta u_{1}=\frac{1}{P_{0}}v_{1}\Delta\theta_{0}+\frac{\theta_{0}}{P_{0}}\Delta v_{1}+\frac{2}{P_{0}}\nabla \theta_{0}\cdot (\nabla v_{1})^{t}.
\end{align}
Combining $\eqref{ge2}_{3}$--$\eqref{ge2}_{4}$ and \eqref{2.1}--\eqref{2.2}, we get
\begin{equation}\label{2.3}
\left\{
\begin{aligned}
	&\operatorname{div}v_{1}=0,\\
	&-\mu\theta_{0}\Delta v_{1}+\nabla (P_{0} P_{2})=(\mu-\frac{\kappa}{2})v_{1}\Delta\theta_{0}-\theta_{0}(v_{1}\cdot \nabla)v_{1}\\
	&\qquad \qquad \qquad\qquad\qquad\,\,  +\zeta\nabla(v_{1}\cdot \nabla \theta_{0})+2\mu\nabla\theta_{0}\cdot (\nabla v_{1})^{t},\\
	&\kappa\Delta\theta_{0}=2v_{1}\cdot \nabla\theta_{0},
\end{aligned}
\right.
\end{equation}
%Denote
%$$
%\begin{aligned}
%&\mathcal{H}_{1}:=\{v_{1}\in H^{m}(\Omega)|(v_{1}\cdot \vec{\iota}_{1},v_{1}\cdot \vec{\iota}_2,v_{1}\cdot \vec{n})\vert_{\partial \Omega}=P_{0}\frac{\beta(T_{w})}{T_{w}}(\partial_{\vec{\iota}_{1}}T_{w},\partial_{\vec{\iota}_{2}}T_{w},0)\},\\
%&\mathcal{H}_{2}:=\{\theta_{0}\in H^{m+1}(\Omega)|\theta_{0}\vert_{\partial \Omega}=T_{w}\}.
%\end{aligned}
%$$
Then we can design a linear mapping $H^{m}\times H^{m+1}\to H^{m}\times H^{m+1}$: $(\tilde{v}_{1},\tilde{\theta}_{0})\mapsto (v_{1},\theta_{0})$
\begin{equation}\label{2.4}
\left\{
\begin{aligned}
	&\operatorname{div}v_{1}=0,\\
	&-\mu\Delta v_{1}+\nabla (\tilde{P}_{0} P_{2})=\mu(\tilde{\theta}_{0}-1)\Delta \tilde{v}_{1}+(\mu-  v\frac{\kappa}{2})\tilde{v}_{1}\Delta\tilde{\theta}_{0}-\tilde{\theta}_{0}(\tilde{v}_{1}\cdot \nabla)\tilde{v}_{1}\\
	&\qquad \qquad \qquad\qquad\quad\,\,\, +\zeta\nabla(\tilde{vv}_{1}\cdot \nabla \tilde{\theta}_{0})+2\mu\nabla\tilde{\theta}_{0}\cdot (\nabla \tilde{v}_{1})^{t}=:\tilde{Z}_{1},\\
	&\kappa\Delta\theta_{0}=2\tilde{v}_{1}\cdot \nabla\tilde{\theta}_{0}=:\tilde{Z}_{2},\\
	&(v_{1}\cdot \vec{\iota}_{1},v_{1}\cdot \vec{\iota}_2,v_{1}\cdot \vec{n})\vert_{\partial \Omega}=\tilde{P}_{0}\frac{h(T_{w})}{T_{w}}(\partial_{\vec{\iota}_{1}}T_{w},\partial_{\vec{\iota}_{2}}T_{w},0),\quad \theta_{0}\vert_{\partial\Omega}=T_{w},
\end{aligned}
\right.
\end{equation}
where we have denoted $\tilde{P}_{0}=M(\int_{\Omega}\tilde{\theta}_{0}^{-1}\,{\rm d}x)^{-1}$. We assume
$$
\|\tilde{v}_{1}\|_{H^{m}}+\|\tilde{\theta}_{0}-1\|_{H^{m+1}}\leq \delta_{0},
$$
where $\delta_{0}\ll 1$ is a small constant and will be chosen later. Then a direct calculation shows that
$$
\tilde{P}_{0}=\frac{M}{\int_{\Omega}\tilde{\theta}_{0}^{-1}\,{\rm d}x}\leq 2M.
$$
Furthermore, it follows from \cite[Theorems IV 5.2 and IV 5.8]{Boyer-Fabrie} and the classical elliptic theory that there exists a unique solution $(v_{1},\theta_{0},P_{2})$ ($\int_{\Omega}P_{2}\,{\rm d}x=0$) of the linear system \eqref{2.4} satisfying
\begin{align}\label{2.5}
	&\|v_{1}\|_{H^{m}}+\|P_{2}\|_{H^{m-1}}\leq K_{1}(\|\tilde{Z}_{1}\|_{H^{m-2}}+|\nabla T_{w}|_{H^{m-\frac{1}{2}}}),\nonumber\\
	&\|\theta_{0}-1\|_{H^{m+1}}\leq K_{1}(\|\tilde{Z}_{2}\|_{H^{m-1}}+|T_{w}-1|_{H^{m+\frac{1}{2}}}),
\end{align}
where $K_{1}>0$ depends only on $M$, $\mu$, $\kappa$, $m$ and $\Omega$. By Sobolev embedding and direct calculations, one has
\begin{align}\label{2.6}
	\|\tilde{Z}_{1}\|_{H^{m-2}}+\|\tilde{Z}_{2}\|_{H^{m-1}}\leq K_{2}(\|\tilde{\theta}_{0}-1\|_{H^{m}}\|\tilde{v}_{1}\|_{H^{m+1}}+\|\nabla\tilde{\theta}\|_{H^{m}}\|\tilde{v}_{1}\|_{H^{m}})\leq 2K_{2}\delta_{0}^2,
\end{align}
where $K_{2}>0$ is constant depends on $m$ and $\Omega$.
%Substituting \eqref{2.6} into $\eqref{2.5}_{2}$, one has
%\begin{equation}\label{2.7-1}
	%\|\theta_{0}-1\|_{H^{m+1}}\leq 2C_{1}C_{2}\delta_{0}^2+C_{1}|\nabla T_{w}|_{H^{m}}
%\end{equation}
%Then letting $2C_{2}C_{1}\delta_{0}\leq \frac{1}{2}$ and $2C_{1}|\nabla T_{w}|_{H^{m}}\leq \frac{1}{2}\delta_{0}$
Substituting \eqref{2.6} into \eqref{2.5}, one
has
\begin{align}\label{2.7}
	\|v_{1}\|_{H^{m}}+\|P_{2}\|_{H^{m-1}}+\|\theta_{0}-1\|_{H^{m+1}}\leq 2K_{1}K_{2}\delta_{0}^2+2K_{1}|\nabla T_{w}|_{H^{m}}.
\end{align}
Then taking $\delta_{0}$ small enough such that $2K_{2}K_{1}\delta_{0}\leq \frac{1}{2}$, and then choosing $|\nabla T_{w}|$ small enough such that  $2K_{1}|\nabla T_{w}|_{H^{m}}\leq \frac{1}{2}\delta_{0}$, we obtain from \eqref{2.7} that
\begin{align*}
	\|v_{1}\|_{H^{m}}+\|P_{2}\|_{H^{m-1}}+\|\theta_{0}-1\|_{H^{m+1}}\leq \delta_{0},
\end{align*}
which implies the linear mapping in \eqref{2.4} is bounded.

By a similar argument, for $(\tilde{v}_{1}^{[k]},\tilde{\theta}_{0}^{[k]})\mapsto (v_{1}^{[k]},\theta_{0}^{[k]})$ with $k=1,2$, one has
\begin{align}\label{2.9}
	&\|v_{1}^{[1]}-v_{1}^{[2]}\|_{H^{m}}+\|P_{2}^{[1]}-P_{2}^{[2]}\|_{H^{m-1}}+\|\theta_{0}^{[1]}-\theta_{0}^{[2]}\|_{H^{m+1}}\nonumber\\
	&\leqslant K_{3}\delta_{0}(\|\tilde{v}_{1}^{[1]}-\tilde{v}_{1}^{[2]}\|_{H^{m}}+\|\nabla\theta_{0}^{[1]}-\nabla\theta_{0}^{[2]}\|_{H^{m}})+K_{3}|\nabla T_{w}|_{H^{m}}\|\theta_{0}^{[1]}-\theta_{0}^{[2]}\|_{H^{m+1}}.
\end{align}
where $K_{3}>0$ depends only on $m$ and $\Omega$. Letting $\delta_{0}$ and $|\nabla T_{w}|_{H^{m}}$ small enough such that $K_{3}(\delta_{0}+|\nabla T_{w}|_{H^{m}})\leq \frac{1}{2}$, we gets from \eqref{2.9} that the linear mapping in \eqref{2.4} is contractive. Hence, the Banach fixed point theorem implies that there exists a unique strong solution $(v_{1},P_{2},\theta_{0})\in H^{m}\times H^{m-1}\times H^{m+1}$ of \eqref{2.3} satisfying
$$
\|v_{1}\|_{H^{m}}+\|P_{2}\|_{H^{m-1}}+\|\theta_{0}-1\|_{H^{m+1}}\lesssim \delta_{0},
$$
Moreover, let $\rho_{0}=P_{0}\theta_{0}^{-1}$ and $u_{1}=\rho_{0}^{-1}v_{1}$, it is clear that  $(u_{1},\rho_{0},\theta_{0},P_{2})$ is the unique strong solution of \eqref{ge2}. Therefore the proof of Lemma \ref{lem2.1} is complete. $\hfill\square$

%%\begin{remark}\label{rem2.1}
%%The introduction of $v_{1}$ is very helpful for us, since the compatibly condition
%%$$
%%\int_{\Omega}\operatorname{div}v\,{\rm d}x=\int_{\partial\Omega}v_{1}\cdot\vec{n}\,{\rm d}x=0
%%$$
%%holds automatically when we apply \cite[Theorem IV 6.6]{Boyer-Fabrie}. If we directly apply a similar argument to establish the existence of $u_{1}$, that is, if we directly study the following linear system for $(u_{1},\theta_{0})$:
%%\begin{align}\label{2.10}
%%\left\{
%%\begin{aligned}
	%%&\operatorname{div}u_{1}=\frac{1}{\tilde{\theta}_{0}}\tilde{u}_{1}\cdot \nabla \tilde{\theta_{0}},\\
	%%&-\mu\Delta u_{1}+\nabla P_{2}=-\frac{\tilde{P}_{0}}{\tilde{\theta}_{0}}(\tilde{u}_{1}\cdot \nabla)\tilde{u}_{1}+\zeta\nabla(\frac{1}{\tilde{\theta}_{0}}\tilde{u}_{1}\cdot \nabla \tilde{\theta}_{0}),\\
	%%&\kappa \Delta\theta_{0}=-\frac{2\tilde{P}_{0}}{\tilde{\theta}_{0}}\tilde{u}_{1}\cdot \nabla \tilde{\theta}_{0},\\
	%%&(u_{1}\cdot \vec{\iota}_{1}, u_{1}\cdot \vec{\iota}_{2},u_{1}\cdot \vec{n})=\beta(T_{w})(\partial_{\vec{\iota}_{1}}T_{w},\partial_{\vec{\iota}_{2}}T_{w},0),\quad \theta_{0}=T_{w}\quad \text{ on }\partial\Omega,
%%\end{aligned}
%%\right.
%%\end{align}
%%where $\tilde{P}_{0}=M(\int_{\Omega}\theta_{0}^{-1}\,{\rm d}x)^{-1}$, then the compatibly condition
%%$$
%%0=\int_{\partial\Omega}u_{1}\cdot \vec{n}\,{\rm d}x=\int_{\Omega}\operatorname{div}u_{1}\,{\rm d}x=\int_{\Omega}\frac{1}{\tilde{\theta}_{0}}\tilde{u}_{1}\cdot \nabla \tilde{\theta_{0}}\,{\rm d}x
%%$$
%%is difficult to be verified.
%%\end{remark}

\section{Existence of the Linearized Problem}

\subsection{Existence of linearized system \eqref{3.1}}
Noting that \eqref{lge2} and \eqref{r2} are slightly decoupled, we shall borrow the idea of \cite{Dou-Jiang-Jiang-Yang} to establish the well-posedness theory of \eqref{lge2} and \eqref{r2}. Here, we would like to point out the systems \eqref{r2} is more complicated than the ones in \cite{Dou-Jiang-Jiang-Yang}, it is hard to apply the arguments in \cite{Dou-Jiang-Jiang-Yang} directly. For later use, we introduce following $\v-$dependent space
\begin{align}\label{K}
	&\mathcal{K}=\{f\in H^3\cap H_{0}^1:\,\,\|f\|_{\mathcal{K}}=:\|u\|_{H^2}+\varepsilon\|\nabla^3 u\|_{L^2}<\infty\}.
	%\qquad \text{with }\|u\|_{\mathcal{K}}=:\|u\|_{H^2}+\varepsilon\|\nabla^3 u\|_{L^2}.
	%\|u\|_{\mathcal{K}}=:\|u\|_{H^2(\Omega)}+\|\nabla (\partial_{ij}^2u)\|_{L^2}+\varepsilon\|\nabla (\partial_{33}^2u)\|_{L^2}\quad (i,j\in \{1,2\}).	
\end{align}
As mentioned in subsection \ref{kp}, the $\v-$dependent functional space $\mathcal{K}$ will be mainly used to deal with the highest normal derivatives of $u_{R}$. Let $\tilde{u}_{R}\in \mathcal{K}$ be a given function, we consider the existence of \eqref{lge2} with replacing $u_{R}$ by $\tilde{u}_{R}$, that is,
\begin{align}\label{3.1}
	\left\{
	\begin{aligned}
		&\rho_{0}\theta_{1}+\rho_{1}\theta_{0}=P_{1},\quad \int_{\Omega}\rho_{1}\,{\rm d}x=0,\\
		&\operatorname{div}(\rho_{0}u_{2})=-\operatorname{div}(\rho_{1}u_{1}),\\
		&\rho_{0}(u_{1}\cdot \nabla)u_{2}+\rho_{0}(u_{2}\cdot \nabla)u_{1}+\nabla P_{3}=-\rho_{1}(u_{1}\cdot \nabla)u_{1}+\mu\Delta u_{2}+\zeta\nabla\operatorname{div}u_{2},\\
		&\kappa\Delta\theta_{1}=-\theta_{0}(u_{2}\cdot \nabla)\rho_{0}+\rho_{1}(u_{1}\cdot \nabla)\theta_{0}+\rho_{0}(u_{1}\cdot \nabla)\theta_{1}+(\rho_{0}\theta_{1}+\rho_{1}\theta_{0})\operatorname{div}u_{1}\\
		&\qquad\quad\,\, +\rho_{0}\theta_{0}\operatorname{div}u_{2}-2\theta_{0}(\tilde{u}_{R}\cdot \nabla)\rho_{0},\\
		&u_{2}=0,\quad \theta_{1}=0,\quad \text{on }\partial\Omega.
	\end{aligned}
	\right.
\end{align}

\begin{lemma}\label{lem2.1}
Let $(\rho_{0},u_{1},\theta_{0})\in H^{5}\times H^{4}\times H^{5}$ be the unique strong solution of \eqref{ge2} established in Lemma \ref{lem1} with $\delta_{0}$ small enough. Then for any given $\tilde{u}_{R}\in \mathcal{K}$,
%constant $C_{0}>0$, there exists a constant $\delta_{2}>0$ depending on $C_{0}$ such that if $\|\tilde{u}_{R}\|_{H^3(\Omega)}\leq \delta_{2}$,
there exists a unique strong solution $(\rho_{1}, u_{2}, \theta_{2}, P_{3})$ ${\rm {(}} \int_{\Omega}P_{3}\,{\rm d}x=0 {\rm{)}}$ of \eqref{3.1} satisfying
\begin{align}\label{3.1-1}
\|(\rho_{1},\theta_{1},u_{2})\|_{H^4}+\|P_{3}\|_{H^3}\lesssim \|\tilde{u}_{R}\|_{H^2}\|\nabla \theta_{0}\|_{H^2}.
%\|u_{2}\|_{H^4(\Omega)}+\|P_{3}\|_{H^3(\Omega)}+\|\theta_{1}\|_{H^4(\Omega)}+\|\rho_{1}\|_{H^4(\Omega)}\lesssim \|\tilde{u}_{R}\|_{H^2}\|\nabla \theta_{0}\|_{H^2}.
\end{align}
\end{lemma}

\noindent\textbf{Proof.} Similar to the proof of Lemma \ref{lem1}, we shall use the Banach fixed point theorem to solve \eqref{3.1}. By $\eqref{3.1}_{1}$, one has $\rho_{1}=\frac{1}{\theta_{0}}(P_{1}-\rho_{0}\theta_{1})$ with $P_{1}=(\int_{\Omega}\rho_{0}\theta_{1}\theta_{0}^{-1}\,{\rm d}x)/(\int_{\Omega}\theta_{0}^{-1}\,{\rm d}x)$.
%$$
%P_{1}=\frac{\int_{\Omega}\frac{\rho_{0}\theta_{1}}{\theta_{0}}\,{\rm d}x}{\int_{\Omega}\theta_{0}^{-1}\,{\rm d}x}.
%$$
Let $v_{2}=\rho_{0}u_{2}$, then $v_{2}\vert_{\partial\Omega}=0$ and
\begin{align}\label{3.3}
\operatorname{div}v_{2}=-\operatorname{div}(\rho_{1}u_{1})&=-(\frac{P_{1}}{\theta_{0}}-\frac{\rho_{0}\theta_{1}}{\theta_{0}})\operatorname{div}u_{1}-u_{1}\cdot \nabla (\frac{P_{1}}{\theta_{0}}-\frac{\rho_{0}\theta_{1}}{\theta_{0}})\nonumber\\
&=\operatorname{div}(\frac{\rho_{0}\theta_{1}}{\theta_{0}}u_{1})-P_{1}\operatorname{div}(\frac{u_{1}}{\theta_{0}}).
\end{align}
Furthermore, it follows from $\eqref{3.1}_{3}$ that
\begin{align}\label{3.4}
  -\frac{\theta_{0}}{P_{0}}\mu\Delta v_{2}+\nabla P_{3}=
&-\frac{\rho_{0}}{P_{0}}(u_{1}\cdot \nabla\theta_{0})v_{2}-(u_{1}\cdot \nabla)v_{2}-(v_{2}\cdot \nabla)u_{1}\nonumber\\
& +\frac{\mu}{P_{0}}v_{2}\Delta\theta_{0}+\frac{2\mu}{P_{0}}\nabla\theta_{0}\cdot( \nabla v_{2})^{t}-\frac{1}{\theta_{0}}(P_{1}-\rho_{0}\theta_{1})(u_{1}\cdot \nabla)u_{1}\nonumber\\
& +\frac{\zeta}{P_{0}}\nabla\big[\theta_{0}\big(\operatorname{div}(\frac{\rho_{0}\theta_{1}}{\theta_{0}}u_{1})-P_{1}\operatorname{div}(\frac{u_{1}}{\theta_{0}})\big)\big]+\frac{\zeta}{P_{0}}\nabla (v_{2}\cdot \nabla \theta_{0}).
%&\qquad\,\,+\frac{\zeta}{P_{0}}\nabla \big[\rho_{0}\theta_{1}\operatorname{div}u_{1}+\theta_{0}u_{1}\cdot \nabla(\frac{\rho_{0}\theta_{1}}{\theta_{0}})-P_{1}\operatorname{div}u_{1}-P_{1}\theta_{0}u_{1}\cdot\nabla(\frac{1}{\theta_{0}})\big]+\frac{\zeta}{P_{0}}\nabla (v_{2}\cdot \nabla \theta_{0}).
\end{align}
Combining \eqref{3.3}, \eqref{3.4} and $\eqref{3.1}_{4}-\eqref{3.1}_{5}$, we can design a linear mapping $H^4\times H^4\to H^{4}\times H^{4}$: $(\tilde{v}_{2},\tilde{\theta}_{1})\mapsto (v_{2},\theta_{1})$
\begin{align}\label{3.2-1}
	\left\{
	\begin{aligned}
	&\operatorname{div}v_{2}=\operatorname{div}(\frac{\rho_{0}\tilde{\theta}_{1}}{\theta_{0}}u_{1})-\tilde{P}_{1}\operatorname{div}(\frac{u_{1}}{\theta_{0}})=:\tilde{Z}_{3},\\
	&-\mu\Delta v_{2}+\nabla (P_{0}P_{3})=(\theta_{0}-1)\mu\Delta \tilde{v}_{2}-\rho_{0}(u_{1}\cdot \nabla \theta_{0})\tilde{v}_{2}-P_{0}(u_{1}\cdot \nabla)\tilde{v}_{2}-P_{0}(\tilde{v}_{2}\cdot \nabla)u_{1}\\
	&\qquad \qquad\qquad\qquad\quad\,\,+\mu\tilde{v}_{2}\Delta\theta_{0}+2\mu\nabla\theta_{0}\cdot (\nabla \tilde{v}_{2})^{t}-\rho_{0}(\tilde{P}_{1}-\rho_{0}\tilde{\theta}_{1})(u_{1}\cdot \nabla)u_{1}\\
	&\qquad \qquad\qquad\qquad\quad\,\,+\zeta\nabla\big[\theta_{0}\big(\operatorname{div}(\frac{\rho_{0}\tilde{\theta}_{1}}{\theta_{0}}u_{1})-\tilde{P}_{1}\operatorname{div}(\frac{u_{1}}{\theta_{0}})\big)\big]+\zeta \nabla (\tilde{v}_{2}\cdot \nabla \theta_{0})=:\tilde{Z}_{4}\\
	%&\qquad \qquad\qquad\quad+\zeta \nabla [\rho_{0}\tilde{\theta}_{1}\operatorname{div}u_{1}+\theta_{0}u_{1}\cdot \nabla (\frac{\rho_{0}\tilde{\theta}_{1}}{\theta_{0}})-\tilde{P}_{1}\operatorname{div}u_{1}-\tilde{P}_{1}\theta_{0}u_{1}\cdot \nabla (\frac{1}{\theta_{0}})]=:\tilde{Z}_{4},\\
	& \kappa \Delta\theta_{1}=-\frac{\theta_{0}}{\rho_{0}}\tilde{v}_{2}\cdot \nabla\rho_{0}+\frac{\tilde{P}_{1}-\rho_{0}\tilde{\theta}_{1}}{\theta_{0}}(u_{1}\cdot \nabla)\theta_{0}+\rho_{0}(u_{1}\cdot \nabla)\tilde{\theta}_{1}+\tilde{v}_{2}\cdot \nabla \theta_{0}\\
	&\qquad \quad\,\, +\rho_{0}\tilde{\theta}_{1}\operatorname{div}u_{1}+\theta_{0}u_{1}\cdot \nabla(\frac{\rho_{0}\tilde{\theta}_{1}}{\theta_{0}})-\tilde{P}_{1}\theta_{0}u_{1}\cdot \nabla (\frac{1}{\theta_{0}})-2\theta_{0}(\tilde{u}_{R}\cdot \nabla )\rho_{0}=:\tilde{Z}_{5},\\
	&v_{2}=0,\quad \theta_{1}=0\quad \text{on }\partial \Omega,
	\end{aligned}
	\right.
\end{align}
where $\tilde{P}_{1}=(\int_{\Omega}\rho_{0}\theta_{0}^{-1}\tilde{\theta}_{1}\,{\rm d}x)/(\int_{\Omega}\theta_{0}^{-1}\,{\rm d}x)$.

In view of \eqref{3.3} and $u_{1}\cdot \vec{n}=0$ on $\partial\Omega$, it easy to check that the terms on the right hand side (RHS) of $\eqref{3.2-1}_{1}$ satisfies the compatible conditions. Therefore, it follows from \cite[Theorems IV 5.2 and IV 5.8]{Boyer-Fabrie} and the classical elliptic theory that there exists a unique solution $(v_{2},\theta_{1},P_{3})$ ($\int_{\Omega}P_{3}\,{\rm d}x=0$) of the linear system \eqref{3.2-1} satisfying
\begin{align}\label{3.5}
&\|v_{2}\|_{H^4}+\|P_{3}\|_{H^3}\leq K_{1}(\|\tilde{Z}_{3}\|_{H^3}+\|\tilde{Z}_{4}\|_{H^2}),\quad \|\theta_{1}\|_{H^4}\leq K_{1}\|\tilde{Z}_{5}\|_{H^2},
\end{align}
where $K_{1}>0$ depends only on $\Omega$.
%For any given $K_{0}>0$, under the assumption
%$$
%\|\tilde{v}_{2}\|_{H^{4}}+\|\tilde{\theta}_{1}\|_{H^5}\leq K_{0},
%$$
Noting that $|\tilde{P}_{1}|\lesssim \|\tilde{\theta}_{1}\|_{L^2}$, by Sobolev embedding and a direct calculation, one obtains that there exists a constant $K_{2}>0$ such that
\begin{align}\label{3.6}
&\|\tilde{Z}_{3}\|_{H^3}\leq K_{2}\|u_{1}\|_{H^4}\|\tilde{\theta}_{1}\|_{H^4},\nonumber\\
&\|\tilde{Z}_{4}\|_{H^2}\leq K_{2}(\|\theta_{0}-1\|_{H^4}+\|u_{1}\|_{H^{4}})(\|\tilde{v}_{2}\|_{H^4}+\|\tilde{\theta}_{1}\|_{H^{4}}),\nonumber\\
&\|\tilde{Z}_{5}\|_{H^2}\leq K_{2}(\|\nabla \theta_{0}\|_{H^4}+\|u_{1}\|_{H^4})(\|\tilde{\theta}_{1}\|_{H^3}+\|\tilde{v}_{2}\|_{H^{2}})+K_{2}\|\tilde{u}_{R}\|_{H^2}\|\nabla\theta_{0}\|_{H^{2}}.
\end{align}
Noting $\delta_{0}\ll 1$ in Lemma \ref{lem1}, we get from \eqref{3.5}--\eqref{3.6} that
\begin{align}\label{3.7}
\|v_{2}\|_{H^4}+\|P_{3}\|_{H^3}+\|\theta_{1}\|_{H^4}\leq \frac{1}{8}(\|\tilde{v}_{2}\|_{H^4}+\|\tilde{\theta}_{1}\|_{H^4})+K_{2}\|\tilde{u}_{R}\|_{H^2}\|\nabla \theta_{0}\|_{H^2}.
\end{align}
Under assumption $\|\tilde{v}_{2}\|_{H^4}+\|\tilde{\theta}_{1}\|_{H^4}\leq 2K_{2}\|\tilde{u}_{R}\|_{H^2}\|\nabla \theta_{0}\|_{H^2}$, it follows from \eqref{3.7} that
%%For any given $C_{0}>0$, letting $\|\tilde{u}_{R}\|_{H^3}\leq \frac{1}{2}C_{0}$, we get from \eqref{3.7} that
\begin{align}\label{3.7-1}
\|v_{2}\|_{H^4}+\|P_{3}\|_{H^3}+\|\theta_{1}\|_{H^5}\leq\frac{5}{4}K_{2}\|\tilde{u}_{R}\|_{H^2}\|\nabla \theta_{0}\|_{H^2},
\end{align}
which implies the linear mapping in \eqref{3.2-1} is bounded. Similar to the proof Lemma \ref{lem1}, we can show the linear mapping in \eqref{3.2-1} is a contraction mapping if $\delta_{0}\ll 1$ in Lemma \ref{lem1}. Hence, it follows from the Banach fixed point theorem that there exists a unique a solution $(u_{2},P_{3},\rho_{1}, \theta_{1})\in H^{4}\times H^3\times H^{4}\times H^4$ of \eqref{3.1} with $u_{2}=\rho_{0}^{-1}v_{2}$ and $\rho_{1}=\theta_{0}^{-1}(P_{1}-\rho_{0}\theta_{1})$. Moreover, \eqref{3.1-1} follows directly from \eqref{3.7-1}.
%%it holds
%%$$
%%\|P_{3}\|_{H^3}+\|\theta_{1}\|_{H^5}\leq C_{0},\quad \|u_{2}\|_{H^4}+\|\rho_{1}\|_{H^2}\lesssim C_{0}.
%%$$
Therefore, the proof of Lemma \ref{lem2.1} is complete. $\hfill\square$

\begin{remark}
	Since $\tilde{u}_{R}\in \mathcal{K}$, then $\tilde{u}_{R}\in H^3(\Omega)$, and one may improve the regularity of $u_{2}$ and $\theta_{1}$ into $H^5\times H^5$. However, $\|\nabla^5(u_{R},\theta_{R})\|_{L^2}$ will depend on $\varepsilon$ due to the dependence of $\|\tilde{u}_{R}\|_{H^3}$ on $\varepsilon$. Fortunately, the uniform bound of $\|(u_{2},\theta_{1})\|_{H^4}$ established in Lemma \ref{lem2.1} will be enough for us. %the following estimate.
\end{remark}

%Noting the definition of $\rho_{2}$ and $\rho_{3}$:
Let $P_{2}$ and $P_{3}$ be the ones obtained in Lemmas \ref{lem1} and \ref{lem2.1}, we define $\rho_{2}$ and $\rho_{3}$ by
\begin{align}\label{3.7-4}
\rho_{2}=\frac{1}{\theta_{0}}(P_{2}-\rho_{1}\theta_{1}-C_{1}),\quad \rho_{3}=\frac{1}{\theta_{0}}(P_{3}-\rho_{2}\theta_{1}-C_{2}),
\end{align}
where $C_{1}$ and $C_{2}$ are the ones in \eqref{C1-1} and \eqref{C2-1} respectively.
Applying Lemmas \ref{lem1} and \ref{lem2.1}, we have following estimates:
\begin{align}
\|\rho_{2}\|_{H^3}&\leq C(\|P_{2}\|_{H^3}+\|\rho_{1}\|_{H^3}\|\theta_{1}\|_{H^3})\lesssim (\delta_{0}+\|\tilde{u}_{R}\|_{H^2}^2\delta_{0}^2),\label{3.7-2}\\
\|\rho_{3}\|_{H^3}&\leq C(\|P_{3}\|_{H^3}+\|\rho_{2}\|_{H^3}\|\theta_{1}\|_{H^3})
\lesssim \delta_{0}\|\tilde{u}_{R}\|_{H^2}(1+\delta_{0}+\|\tilde{u}_{R}\|_{H^2}^2\delta_{0}^2).\label{3.7-3}
\end{align}

\subsection{Existence of linearized system \eqref{3.10} for $(\mu_{0},\zeta_{0})$}
Let $\tilde{u}_{R}\in \mathcal{K}$ and $\tilde{\theta}_{R}\in (H_{0}^1\cap H^3)$,
$(\rho_{0},u_{1},\theta_{1},P_{2})$ be the unique solution of \eqref{ge2} established in Lemma \ref{lem1}, $(\rho_{2},u_{2},\theta_{2},P_{3})$ be the unique solution of \eqref{lge2} established in Lemma \ref{lem2.1}, and $\rho_{2}, \rho_{3}$ are given in \eqref{3.7-4}. Now we consider following linearized system of \eqref{r2}:
\begin{equation}\label{3.10}
\left\{
\begin{aligned}
	&\operatorname{div}\big[\rho_{R}(u_{1}+\v(u_{2}+\tilde{u}_{R}))\big]+\frac{1}{\v}\operatorname{div}(\rho_{0}u_{R})
	=\mathfrak{R}_{1},\\
	&\mu\Delta u_{R}+\zeta\nabla \operatorname{div}u_{R}-\frac{1}{\varepsilon}\nabla(\rho_{0}\theta_{R}+\rho_{R}\theta_{0})=\rho_{0}(u_{1}\cdot \nabla u_{R}+u_{R}\cdot \nabla u_{1})+\nabla (\rho_{R}\theta_{1}+\rho_{1}\theta_{R})\\
	&\qquad\qquad\qquad\qquad\qquad\qquad\qquad\qquad \qquad\, +\tilde{F}^{\varepsilon}(\rho_{R},u_{R},\theta_{R})+\mathfrak{R}_{2},\\
	&\frac{\kappa}{\theta_{0}}\Delta \theta_{R}-\frac{1}{\varepsilon}\operatorname{div}(\rho_{0}u_{R})=\frac{P_{1}}{\theta_{0}}\operatorname{div}u_{R}+\frac{1}{\theta_{0}}(\rho_{0}\theta_{R}+\rho_{R}\theta_{0})\operatorname{div}u_{1}+\frac{1}{\theta_{0}}\tilde{G}^{\varepsilon}(\rho_{R},u_{R},\theta_{R})+\mathfrak{R}_{3}\\
	&\qquad\qquad\qquad\qquad\qquad\,\,+\frac{1}{\theta_{0}}[\rho_{0}(u_{R}\cdot \nabla \theta_{1}+u_{1}\cdot \nabla \theta_{R})+\rho_{1}u_{R}\cdot \nabla\theta_{0}+\rho_{R}u_{1}\cdot \nabla \theta_{0}],\\
	&u_{R}\vert_{\partial\Omega}=\theta_{R}\vert_{\Omega}=0,\qquad \int_{\Omega}\rho_{R}\,{\rm d}x=0,
\end{aligned}
\right.
\end{equation}
where
\begin{align}\label{3.10-1}
&\mathfrak{R}_{1}=:-\operatorname{div}\big[(\rho_{1}+\v\rho_{2}+\v^2\rho_{3})\tilde{u}_{R}\big]+r_{1},\quad \mathfrak{R}_{2}=:r_{2},\nonumber\\
&\mathfrak{R}_{3}=:\frac{1}{\theta_{0}}r_{3}-\frac{1}{\theta_{0}}\Psi(\nabla (u_{1}+\varepsilon(u_{2}+\tilde{u}_{R}))),
\end{align} $\tilde{F}^{\varepsilon}(\rho_{R},u_{R},\theta_{R})$ and $\tilde{G}^{\varepsilon}(\rho_{R},u_{R},\theta_{R})$ are the linearized forms of $F^{\varepsilon}(\rho_{R},u_{R},\theta_{R})$ and $G^{\varepsilon}(\rho_{R},u_{R},\theta_{R})$ respectively, i.e.,
$$
\tilde{F}^{\varepsilon}(\rho_{R},u_{R},\theta_{R})=\tilde{f}_{1}^{\varepsilon}(u_{R},\theta_{R})+\tilde{f}_{2}^{\varepsilon}(\rho_{R})\,\,\text{ and }\,\, \tilde{G}^{\varepsilon}(\rho_{R},u_{R},\theta_{R})=\tilde{g}_{1}^{\varepsilon}(u_{R},\theta_{R})+\tilde{g}_{2}^{\varepsilon}(\rho_{R}),$$
where
\begin{align}
	&\tilde{f}_{1}^{\varepsilon}(u_{R},\theta_{R})=\varepsilon\big[\rho_{0}(u_{2}\cdot \nabla u_{R}+u_{R}\cdot \nabla u_{2}+\tilde{u}_{R}\cdot \nabla u_{R})+\rho_{1}(u_{R}\cdot \nabla u_{1}+ u_{1}\cdot \nabla u_{R})\big]\nonumber\\
	&\qquad\qquad \quad\,\, +\varepsilon^2\big[\rho_{1}(u_{R}\cdot \nabla u_{2}+u_{2}\cdot \nabla u_{R}+\tilde{u}_{R}\cdot \nabla u_{R})+\rho_{2}(u_{R}\cdot \nabla u_{1}+u_{1}\cdot \nabla u_R)\big]\nonumber\\
	&\qquad\qquad \quad\,\, +\varepsilon^3\big[\rho_{2}(u_{R}\cdot \nabla u_{2}+u_{2}\cdot \nabla u_{R}+\tilde{u}_{R}\cdot \nabla u_{R})+\rho_{3}(u_{R}\cdot \nabla u_{1}+u_{1}\cdot \nabla u_{R})\big]\nonumber\\
	&\qquad\qquad \quad\,\, +\varepsilon^{4}\rho_{3}(u_{R}\cdot \nabla u_{2}+u_2\cdot u_{R}+\tilde{u}_{R}\cdot \nabla u_{R})+\varepsilon\nabla (\rho_{2}\theta_{R})+\varepsilon^2\nabla (\rho_{3}\theta_{R}),\label{4.5}\\
	&\tilde{f}_{2}^{\varepsilon}(\rho_{R})=\varepsilon\rho_{R}u_{1}\cdot \nabla u_{1}+\varepsilon^2\rho_{R}(u_{1}\cdot \nabla \tilde{u}_{R}+\tilde{u}_{R}\cdot \nabla u_{1}+u_{2}\cdot \nabla u_{1}+u_{1}\cdot \nabla u_{2})\nonumber\\
	&\qquad \qquad+\varepsilon^3\rho_{R}(\tilde{u}_{R}\cdot \nabla u_{2}+u_{2}\cdot \nabla \tilde{u}_{R}+u_{2}\cdot \nabla u_{2}+\tilde{u}_{R}\cdot \nabla \tilde{u}_{R})+\varepsilon\nabla(\rho_{R}\tilde{\theta}_{R}),\label{4.6}
\end{align}
%and
\begin{align}
	&\tilde{g}_{1}^{\varepsilon}(u_{R},\theta_{R})=\varepsilon\big[\rho_{0}(u_{2}\cdot\nabla \theta_{R}+\tilde{u}_{R}\cdot \nabla \theta_{R})+\rho_{2}u_{R}\cdot \nabla\theta_{0}+\rho_{1}(u_{R}\cdot \nabla \theta_{1}+u_{1}\cdot \nabla \theta_{R})\big]\nonumber\\
	&\qquad\qquad \quad\, +\varepsilon^2\big[\rho_{1}(u_{2}\cdot\nabla \theta_{R}+\tilde{u}_{R}\cdot \nabla \theta_{R})+\rho_{2}(u_{R}\cdot \nabla\theta_{1}+u_{1}\cdot \nabla \theta_{R})+\rho_{3}u_{R}\cdot \nabla \theta_{0}\big]\nonumber\\
	&\qquad\qquad \quad\,  +\varepsilon^3\big[\rho_{2}(u_{2}\cdot\nabla \theta_{R}+\tilde{u}_{R}\cdot \nabla \theta_{R})+\rho_{3}(u_{1}\cdot \nabla \theta_{R}+u_{R}\cdot \nabla \theta_{1})\big]\nonumber\\
	&\qquad\qquad \quad\,  +\varepsilon\big[\rho_{0}\theta_{R}\operatorname{div}u_{2}+(\rho_{1}\theta_{1}+\rho_{2}\theta_{0})\operatorname{div}u_{R}+\rho_{0}\theta_{R}\operatorname{div}\tilde{u}_{R}+\rho_{1}\theta_{R}\operatorname{div}u_{1}\big]\nonumber\\
	&\qquad\qquad \quad\,  +\varepsilon^2\big[\rho_{1}\theta_{R}\operatorname{div}u_{2}+(\rho_{2}\theta_{1}+\rho_{3}\theta_{0})\operatorname{div}u_{R}+\rho_{1}\theta_{R}\operatorname{div}\tilde{u}_{R}\nonumber +\rho_{2}\theta_{R}\operatorname{div}u_{1}\big]\nonumber\\
	&\qquad\qquad \quad\,  +\varepsilon\big[\rho_{3}\theta_{R}\operatorname{div}u_{1}+\rho_{2}\theta_{R}\operatorname{div}u_{2}+\rho_{3}\theta_{1}\operatorname{div}u_{R}+\rho_{2}\theta_{R}\operatorname{div}\tilde{u}_{R}\big]\nonumber\\
	&\qquad\qquad \quad\, +\varepsilon^{4}\rho_{3}(u_{2}\cdot \nabla \theta_{R}+\tilde{u}_{R}\cdot \nabla \theta_{R})+\varepsilon^{4}[\rho_{3}\theta_{R}\operatorname{div}u_2+\rho_{3}\theta_{R}\operatorname{div}\tilde{u}_R],\label{4.7}\\
	&\tilde{g}_{2}^{\varepsilon}(\rho_{R})=\varepsilon\rho_{R}(u_{2}\cdot \nabla\theta_{0}+\tilde{u}_{R}\cdot \nabla \theta_{0}+u_{1}\cdot \nabla \theta_{1})+\varepsilon^2\rho_{R}(u_{2}\cdot \nabla\theta_{1}+\tilde{u}_{R}\cdot \nabla \theta_{1}+u_{1}\cdot \nabla \tilde{\theta}_{R})\nonumber\\
	&\qquad\quad\,\,\,\, +\varepsilon^3\rho_{R}(u_{2}\cdot\nabla \tilde{\theta}_{R}+\tilde{u}_{R}\cdot \nabla \tilde{\theta}_{R})+\varepsilon\rho_{R} \big[\theta_{0}\operatorname{div}u_{2}+\theta_{0}\operatorname{div}\tilde{u}_{R}+\theta_{1}\operatorname{div}u_{1}\big]\nonumber\\
	&\qquad \quad\,\,\,\,
	+\varepsilon^2\rho_{R}\big[\theta_{1}\operatorname{div}u_{2}+\theta_{1}\operatorname{div}\tilde{u}_{R}+\tilde{\theta}_{R}\operatorname{div}u_{1}\big]
	+\varepsilon^3\rho_{R}\big[\tilde{\theta}_{R}\operatorname{div}\tilde{u}_{R}+\tilde{\theta}_{R}\operatorname{div}u_{2}\big].\label{4.8}
\end{align}

Motivated by \cite{Novotny-Padula-1994,Novotny-Straskraba,Valli-1987}, to establish the strong solution of \eqref{3.10}, we shall decompose \eqref{3.10} into the system \eqref{3.10} into two Laplace equations with Dirichlet and Neumann boundary conditions respectively, a classical Stokes problem and a steady transport equation with the help of effective viscous flux and the Helmholtz decomposition. This decomposition is helpful for us to construct approximate solutions. It follows from $\eqref{3.10}_{1}$ and $\eqref{3.10}_{3}$ that
\begin{align}\label{5.20-1}
	\frac{\kappa}{\theta_{0}}\Delta\theta_{R}&=-\operatorname{div}[\rho_{R}(u_{1}+\v(u_{2}+\tilde{u}_{R}))]+\mathfrak{R}_{1}+\frac{P_{1}}{\theta_{0}}\operatorname{div}u_{R}+\frac{1}{\theta_{0}}(\rho_{0}\theta_{R}+\rho_{R}\theta_{0})\operatorname{div}u_{1}\nonumber\\
	&\quad +\frac{1}{\theta_{0}}[\rho_{0}(u_{R}\cdot \nabla \theta_{1}+u_{1}\cdot \nabla \theta_{R})+\rho_{1}u_{R}\cdot \nabla\theta_{0}+\rho_{R}u_{1}\cdot \nabla \theta_{0}]\nonumber\\
	&\quad +\frac{1}{\theta_{0}}\tilde{G}^{\v}(\rho_{R},u_{R},\theta_{R})+\mathfrak{R}_{3}.
\end{align}
Denoting $w_{R}=\rho_{0}u_{R}$, one has from $P_{0}=\rho_{0}\theta_{0}$ and direct calculations that
\begin{align}
	&\operatorname{div}u_{R}=\frac{1}{P_{0}}\operatorname{div}(w_{R}\theta_{0})=\frac{1}{P_{0}}w_{R}\cdot \nabla \theta_{0}+\frac{1}{P_{0}}\theta_{0}\operatorname{div}w_{R},\label{F1}\\
	&\Delta u_{R}=\frac{1}{P_{0}}w_{R}\Delta \theta_{0}+\frac{\theta_{0}}{P_{0}}\Delta w_{R}+\frac{2}{P_{0}}\nabla \theta_{0}\cdot (\nabla w_{R})^{t}\label{F1-1}.
\end{align}
Then it follows from $\eqref{3.10}_{2}$ that
\begin{align}\label{F2}
	\mu \theta_{0}\Delta w_{R}+\nabla \tilde{P}_{R}&=\rho_{0}\big[w_{R}(u_{1}\cdot \nabla)\theta_{0}+\theta_{0}(u_{1}\cdot \nabla)w_{R}+\theta_{0}(w_{R}\cdot \nabla)u_{1}\big]-\mu w_{R}\Delta\theta_{0}\nonumber\\
	&\quad -2\mu \nabla \theta_{0}\cdot (\nabla w_{R})^{t}-\zeta w_{R}\cdot \nabla \theta_{0}+\mathscr{F}^{\v}(\rho_{R},w_{R},\theta_{R})+P_{0}\mathfrak{R}_{2}
\end{align}
where
\begin{align}\label{F3}
	\tilde{P}_{R}:=\zeta \theta_{0}\operatorname{div}w_{R}-\frac{P_{0}}{\v}(\rho_{0}\theta_{R}+\rho_{R}\theta_{0})-P_{0}(\theta_{1}+\v\tilde{\theta}_{R})\rho_{R}-P_{0}(\rho_{1}+\v\rho_{2}+\v^2\rho_{3})\theta_{R},
\end{align}
and
\begin{align}\label{F4}
	\mathscr{F}^{\v}(\rho_{R},w_{R},\theta_{R})
	:=&P_{0}\tilde{F}^{\v}(\rho_{R},\frac{w_{R}}{\rho_{0}},\theta_{R})-P_{0}\v\nabla(\rho_{R}\tilde{\theta}_{R}+\rho_{2}\theta_{R})-P_{0}\v^2\nabla (\rho_{3}\theta_{R}).
\end{align}
As in \cite{Novotny-Padula-1994,Novotny-Straskraba}, we use the Helmholtz decomposition to set
$$
w_{R}=v_{R}+\nabla q_{R}\text{ with }\operatorname{div}v_{R}=0\text{ and }v_{R}\vert_{\partial\Omega}=-\nabla q_{R},\,\,\frac{\partial q_{R}}{\partial \vec{n}}\vert_{\partial\Omega}=0.
$$
Then it follows from $\eqref{3.10}_{1}$ and $\eqref{3.10}_{2}$ that
\begin{align}\label{E1-1}
	\Delta q_{R}&=\operatorname{div}w_{R}=-\v\operatorname{div}[\rho_{R}(u_{1}+\varepsilon(u_{2}+\tilde{u}_{R}))]+\v\mathfrak{R}_{1},
\end{align}
and
\begin{align}\label{E2}
	\mu\theta_{0} \Delta v_{R}+\nabla P_{R}&=\rho_{0}\big[w_{R}(u_{1}\cdot \nabla)\theta_{0}+\theta_{0}(u_{1}\cdot \nabla)w_{R}+\theta_{0}(w_{R}\cdot \nabla)u_{1}\big]+\mu\operatorname{div}w_{R} \nabla \theta_{0}\nonumber\\
	&\quad -\mu w_{R}\Delta\theta_{0}-2\mu \nabla \theta_{0}\cdot (\nabla w_{R})^{t}-\zeta w_{R}\cdot \nabla \theta_{0}+\mathscr{F}^{\v}(\rho_{R},w_{R},\theta_{R})+P_{0}\mathfrak{R}_{2},
\end{align}
where
\begin{align}\label{E3}
	P_{R}&:=\tilde{P}_{R}+\mu\theta_{0} \operatorname{div} w_{R}\nonumber\\
	&=\theta_{0}(\zeta +\mu)\operatorname{div}w_{R}-\frac{P_{0}}{\v}(\rho_{0}\theta_{R}+\rho_{R}\theta_{0})-P_{0}(\theta_{1}+\v\tilde{\theta}_{R})\rho_{R}-P_{0}(\rho_{1}+\v\rho_{2}+\v^2\rho_{3})\theta_{R}.
\end{align}
%$$
%\mathscr{F}(\rho_{R},u_{R},\theta_{R})=:\tilde{F}^{\v}(\rho_{R},u_{R},\theta_{R})-\nabla(\rho_{R}\theta_{1}+\rho_{1}\theta_{R})-\v\nabla(\rho_{R}\tilde{\theta}_{R}+\rho_{2}\theta_{R})-\v^2\nabla(\rho_{3}\theta_{R}),
%$$
%and
%\begin{align}\label{E3}
%	P_{R}=(\mu+\zeta)\operatorname{div}u_{R}-\frac{1}{\v}(\rho_{R}\theta_{0}+\rho_{0}\theta_{R})-(\theta_{1}+\v\tilde{\theta}_{R})\rho_{R}-(\rho_{1}+\v\rho_{2}+\v^2\rho_{3})\theta_{R}.
%\end{align}
It follows from $\eqref{3.10}_{1}$ and \eqref{E3} that
%\begin{align}\label{E4-0}
%&P_{0}\rho_{R}+\varepsilon^2(\mu+\zeta)\theta_{0}\operatorname{div}[\rho_{R}(u_{1}+\varepsilon(u_{2}+\tilde{u}_{R}))]\nonumber\\
%&=-\var P_{R}-\v^2(\mu+\zeta)\theta_{0}\{\operatorname{div}[\tilde{u}_{R}(\rho_{1}+\varepsilon\rho_{2}+\var^2\rho_{3})]-r_{1}\}\nonumber\\
%&\quad -P_{0}(\rho_{0}+\v\rho_{1}+\v^2\rho_{2}+\v^3\rho_{3})\theta_{R}-P_{0}(\theta_{0}-1+\var\theta_{1}+\var^2\tilde{\theta}_{R})\rho_{R}
%&\quad -(\rho_{0}+\v\rho_{1}+\v^2\rho_{2}+\v^3\rho_{3})\theta_{R}-(\theta_{0}-1+\var\theta_{1}+\var^2\tilde{\theta}_{R})\rho_{R}.
%\end{align}
%which is equivalent to
\begin{align}\label{E4}
	&\rho_{R}+\frac{\varepsilon^2(\mu+\zeta)}{P_{0}}\operatorname{div}[\rho_{R}\theta_{0}(u_{1}+\varepsilon(u_{2}+\tilde{u}_{R}))]\nonumber\\
	&=- \frac{\var}{P_{0}}P_{R}+\frac{\v^2(\mu+\zeta)}{\rho_{0}}\mathfrak{R}_{1}
	%-\frac{\v^2(\mu+\zeta)}{\rho_{0}}\{\operatorname{div}[\tilde{u}_{R}(\rho_{1}+\varepsilon\rho_{2}+\var^2\rho_{3})]-r_{1}\}
	-(\rho_{0}+\v\rho_{1}+\v^2\rho_{2}+\v^3\rho_{3})\theta_{R}\nonumber\\
	&\quad -(\theta_{0}-1+\var\theta_{1}+\var^2\tilde{\theta}_{R})\rho_{R}+\frac{\varepsilon^2(\mu+\zeta)}{P_{0}}\rho_{R}(u_{1}+\varepsilon(u_{2}+\tilde{u}_{R}))\cdot \nabla\theta_{0}.
\end{align}

For later use, we denote
\begin{align*}
	\mathcal{H}^2&=\{f\in H^2:\v^2\|\nabla^2 \operatorname{div}[f(u_{1}+\v(u_{2}+\tilde{u}_{R}))]\|_{L^2}<\infty\},
\end{align*}
\begin{align}
	\mathcal{I}_{0}&:=\|(\rho_{0}-P_{0},u_{1},P_{2},\theta_{0}-1)\|_{H^3}+\|(\rho_{1},u_{2},P_{3},\theta_{1})\|_{H^3} +\|(\rho_{2},\rho_{3})\|_{H^3}\nonumber\\
	&\lesssim \delta_{0}(1+\|\tilde{u}_{R}\|_{H^2}),\label{I0}
\end{align}
and
\begin{align}\label{I1}
	\mathcal{I}_{1}:=\mathcal{I}_{0}+\|\tilde{u}_{R}\|_{\mathcal{K}}+\|\tilde{\theta}_{R}\|_{H^3}\lesssim \delta_{0}+\|\tilde{u}_{R}\|_{\mathcal{K}}+\|\tilde{\theta}_{R}\|_{H^3}.
\end{align}
Hereafter, we denote $C(\mathcal{I}_{0})$ and $C(\mathcal{I}_{1})$ as two small positive constants depending on $\mathcal{I}_{0}$ and $\mathcal{I}_{1}$ respectively, and satisfying $C(\mathcal{I}_{0}), C(\mathcal{I}_{1})\to 0$ as $\mathcal{I}_{0},\mathcal{I}_{1}\to 0$ respectively.

\begin{lemma}\label{lem7.1}
	For any fixed $0<\v<1$, there exist a positive constant $K_{0}\gg 1$, and two positive viscous coefficients $\mu_{0}$ and $\zeta_{0}$
	satisfying $0<\mu_{0},\zeta_{0}<1$ and $\frac{\zeta_{0}}{\mu_{0}}\geq K_{0}$ such that if $\delta_{0}$ and $\|\tilde{u}_{R}\|_{\mathcal{K}}+\|\tilde{\theta}_{R}\|_{H^3}$ are small enough, there exists a strong solution $(\rho_{R},u_{R},\theta_{R})\in H^2\times \mathcal{K}\times (H_{0}^1\cap H^3)$ of \eqref{3.10} with $(\mu,\zeta)=(\mu_{0},\zeta_{0})$. Moreover, it holds
	\begin{align}\label{g1}
		\|\rho_{R}\|_{H^2}+\|u_{R}\|_{\mathcal{K}}+\|\theta_{R}\|_{H^3}\leq C\big[\|\mathfrak{R}_{2}\|_{L^2}+\|(\mathfrak{R}_{1},\mathfrak{R}_{3})\|_{H^1}\big]+C\v\big[\|\mathfrak{R}_{2}\|_{H^1}+\v\|\mathfrak{R}_{1}\|_{H^2}\big].
	\end{align}
\end{lemma}

\begin{remark}\label{rem3.1}
We would like to point out that the positive constant $C$ in \eqref{g1} may blow up as $\frac{\zeta_{0}}{\mu_{0}}\to K_{0}$. Therefore, to establish the existent of the strong solution of \eqref{3.10} for any given $(\mu,\zeta)$, we still need to establish {\it a priori} uniform estimates depending continuously on $(\mu,\zeta)$ in next section, see also \cite{Valli-1987}.
\end{remark}

\noindent\textbf{Proof.} In the following, we shall use \eqref{5.20-1}, \eqref{E1-1}, \eqref{E2} and \eqref{E4} to construct solution. Since the proof is very long, we divide it into three steps.

{\it Step 1}. Setting
\begin{align}\label{E5}
	(\rho_{R}^{(0)},\theta_{R}^{(0)},v_{R}^{(0)},q_{R}^{(0)})=(0,0,0,0,0,0),
\end{align}
and $w_{R}^{(0)}=v_{R}^{(0)}+\nabla q_{R}^{(0)}=0$.
%It is clear that $(\rho_{R}^{(0)},\theta_{R}^{(0)},v_{R}^{(0)},q_{R}^{(0)})\in \mathcal{H}^2\times H_{0}^3\times H_{0}^{3}\times H^{4}$.
Once $(\rho_{R}^{(k)},\theta_{R}^{(k)},v_{R}^{(k)},q_{R}^{(k)})\in \mathcal{H}^2\times (H^3\cap H_{0}^1)\times (H^3\cap H_{0}^1)\times H^{4}$ are given, we first construct $\theta_{R}^{(k+1)}$ and $q_{R}^{(k+1)}$ by
\begin{align}\label{E6}
	\left\{
	\begin{aligned}
		\frac{\kappa}{\theta_{0}}\Delta\theta_{R}^{(k+1)}&=-\operatorname{div}[\rho_{R}^{(k)}(u_{1}+\v u_{2}+\v\tilde{u}_{R})]+\mathfrak{R}_{1}+\frac{1}{\theta_{0}}\tilde{G}^{\v}\big(\rho_{R}^{(k)},\frac{w_{R}^{(k)}}{\rho_{0}},\theta_{R}^{(k)}\big)+\mathfrak{R}_{3}\\
		&\quad +\frac{P_{1}}{P_{0}\theta_{0}}\big(w_{R}^{(k)}\cdot \nabla \theta_{0}+\theta_{0}\operatorname{div}w_{R}^{(k)}\big)+\frac{1}{\theta_{0}}(\rho_{0}\theta_{R}^{(k)}+\rho_{R}^{(k)}\theta_{0})\operatorname{div}u_{1}\\
		&\quad +\frac{1}{\theta_{0}}\big[(w_{R}^{(k)}\cdot \nabla \theta_{1}+\rho_{0}u_{1}\cdot \nabla \theta_{R}^{(k)})+\frac{\rho_{1}}{\rho_{0}}w_{R}^{(k)}\cdot \nabla \theta_{0}+\rho_{R}^{(k)}u_{1}\cdot \nabla \theta_{0}\big]\\
		&=:\mathcal{R}_{1}(\rho_{R}^{(k)},w_{R}^{(k)},\theta_{R}^{(k)})\\
		\theta_{R}^{(k+1)}\vert_{\partial \Omega}&=0,
	\end{aligned}
	\right.
\end{align}
and
\begin{align}\label{E7}
	\left\{
	\begin{aligned}
		&\Delta q_{R}^{(k+1)}=-\v\operatorname{div}[\rho_{R}^{(k)}(u_{1}+\varepsilon(u_{2}+\tilde{u}_{R}))]+\v\mathfrak{R}_{1},=:\mathcal{R}_{2}(\rho_{R}^{(k)}),\\
		&\frac{\partial q_{R}^{(k+1)}}{\partial \vec{n}}\vert_{\partial\Omega}=0,
	\end{aligned}
	\right.
\end{align}
where $w_{R}^{(k)}=v_{R}^{(k)}+\nabla q_{R}^{(k)}$. It follows from standard theory of the Laplace equation with Dirichlet boundary condition ({\it cf}. \cite[Theorems III.4.2]{Boyer-Fabrie}) that there exist a unique solution $\theta_{R}^{(k+1)}\in H^3\cap  H_{0}^1$ of \eqref{E6} satisfying
%\begin{align}\label{F5}
%	\|\theta_{R}^{(k+1)}\|_{H^2}&\leq C\|R_{1}(\rho_{R}^{(k)},w_{R}^{(k)},\theta_{R}^{(k)})\|_{L^2}\nonumber\\
%	&\leq C(\mathcal{I}_{1})\|(\rho_{R}^{(k)},w_{R}^{(k)},\theta_{R}^{(k)})\|_{H^1}+\|(\mathfrak{R}_{1},\mathfrak{R}_{3})\|_{L^2}.
%\end{align}
%and
\begin{align}\label{F6}
	\|\theta_{R}^{(k+1)}\|_{H^3}&\leq C\|\mathcal{R}_{1}(\rho_{R}^{(k)},w_{R}^{(k)},\theta_{R}^{(k)})\|_{H^1}\leq C(\mathcal{I}_{1})\|(\rho_{R}^{(k)},w_{R}^{(k)},\theta_{R}^{(k)})\|_{H^2}+C\|(\mathfrak{R}_{1},\mathfrak{R}_{3})\|_{H^1},
\end{align}
where we have used that fact
\begin{equation}\label{F6-1}
\begin{aligned}
&\|\tilde{G}^{\v}(\rho_{R},u_{R},\theta_{R})\|_{L^2}\leq \|\tilde{g}_{1}^{\v}(u_{R},\theta_{R})\|_{L^2}+\|\tilde{g}_{2}^{\v}(\rho_{R})\|_{L^2}\leq C(\mathcal{I}_{1})\|(\rho_{R},u_{R},\theta_{R})\|_{H^1},\\
&\|\tilde{G}^{\v}(\rho_{R},u_{R},\theta_{R})\|_{H^1}\leq \|\tilde{g}_{1}^{\v}(u_{R},\theta_{R})\|_{H^1}+\|\tilde{g}_{2}^{\v}(\rho_{R})\|_{H^1}\leq C(\mathcal{I}_{1})\|(\rho_{R},u_{R},\theta_{R})\|_{H^2}.
\end{aligned}
\end{equation}

Noting \eqref{3.10-1} and \eqref{r1-1}, the solvability condition of \eqref{E7} holds automatically. It follows from standard theory of the Laplace equation with Neumann boundary condition ({\it cf}. \cite[Theorems III.4.3]{Boyer-Fabrie}) that there exist a unique solution $q_{R}^{(k+1)}\in H^4$ with $\int_{\Omega}q_{R}^{(k+1)}\,{\rm d}x=0$ of \eqref{E7} satisfying
\begin{align}\label{F7}
	\|q_{R}^{(k+1)}\|_{H^3}\leq C\|\mathcal{R}_{2}(\rho_{R}^{(k)})\|_{H^1}\leq C(\mathcal{I}_{1})\v\|\rho_{R}^{(k)}\|_{H^2}+C\v\|\mathfrak{R}_{1}\|_{H^1},
\end{align}
and
\begin{align}\label{F7-1}
	\v\|q_{R}^{(k+1)}\|_{H^4}\leq \bar{C}\v\|\mathcal{R}_{2}(\rho_{R}^{(k)})\|_{H^2}&\leq
	\bar{C}\v^2\|\nabla^2\operatorname{div}[\rho_{R}^{(k)}(u_{1}+\v(u_{2}+\tilde{u}_{R}))]\|_{L^2}\nonumber\\
	&\quad +C(\mathcal{I}_{1})\v^2\|\rho_{R}^{(k)}\|_{H^2}+\bar{C}\v^2\|\mathfrak{R}_{1}\|_{H^2}.
\end{align}
We would like to point out that the constant $\bar{C}>0$ in \eqref{F7-1} depends only on $\Omega$.

{\it Step 2}. Next, we construct $v_{R}^{(k+1)}$ and $P_{R}^{(k+1)}$ by
\begin{align}\label{E8}
	\left\{
	\begin{aligned}
		&\mu \Delta v_{R}^{(k+1)}+\nabla P_{R}^{(k+1)}=\rho_{0}\big[w_{R}^{(k)}(u_{1}\cdot \nabla)\theta_{0}+\theta_{0}(u_{1}\cdot \nabla)w_{R}^{(k)}+\theta_{0}(w_{R}^{(k)}\cdot \nabla)u_{1}\big]\\
		&\qquad\qquad\qquad\qquad\qquad-\mu(\theta_{0}-1)\Delta v_{R}^{(k)}-\mu w_{R}^{(k)}\Delta\theta_{0}+\mu\operatorname{div}w_{R}^{(k)}\nabla\theta_{0}-2\mu \nabla \theta_{0}\cdot (\nabla w_{R}^{(k)})^{t}\\
		&\qquad\qquad\qquad\qquad\qquad-\zeta w_{R}^{(k)}\cdot \nabla \theta_{0}+\mathscr{F}^{\v}(\rho_{R}^{(k)},w_{R}^{(k)},\theta_{R}^{(k)})+P_{0}\mathfrak{R}_{2}\\
		&\qquad\qquad\qquad\qquad\quad:=\mathcal{R}_{3}(\rho_{R}^{(k)},w_{R}^{(k)},\theta_{R}^{(k)}),\\
		&\operatorname{div}v_{R}^{(k+1)}=0,\\
		&v_{R}^{(k+1)}\vert_{\partial \Omega}=-\nabla q_{R}^{(k+1)}\vert_{\partial\Omega},
		%=-\int_{\Omega}\frac{\v(\mu+\zeta)}{\rho_{0}}\big\{\frac{1}{\varepsilon}u_{R}^{(k)}\cdot \nabla \rho_{0}+\operatorname{div}[\tilde{u}_{R}(\rho_{1}+\varepsilon\rho_{2}+\var^2\rho_{3})]+r_{1}\}\,{\rm d}x\\
		%&\qquad \qquad \qquad\,\,\,\,\, -\frac{1}{\v}\int_{\Omega}\big[(\rho_{0}+\v\rho_{1}+\v^2\rho_{2}+\v^3\rho_{3})\theta_{R}^{(k)}-(\theta_{0}-1+\var\theta_{1}+\var^2\tilde{\theta}_{R})\rho_{R}^{(k)}\big]\,{\rm d}x
	\end{aligned}
	\right.
\end{align}
with the constraint
\begin{align}\label{E8-1}
	&\int_{\Omega}P_{R}^{(k+1)}\,{\rm d}x=P_{0}(\mu+\zeta)\v\int_{\Omega}\frac{1}{\rho_{0}}\mathfrak{R}_{1}\,{\rm d}x
	%\big\{\operatorname{div}[\tilde{u}_{R}(\rho_{1}+\v\rho_{2}+\v^2\rho_{3})]-r_{1}\big\}\,{\rm d}x\\
	-\frac{P_{0}}{\v}\int_{\Omega}(\rho_{0}+\v\rho_{1}+\v^2\rho_{2}+\v^3\rho_{3})\theta_{R}^{(k+1)}\,{\rm d}x\nonumber\\
	&\qquad \qquad \qquad\,\,\,\,\,+\v(\mu+\zeta)\int_{\Omega}\rho_{R}^{(k)}(u_{1}+\v(u_{2}+\tilde{u}_{R}))\cdot \nabla\theta_{0}\,{\rm d}x\nonumber\\
	&\qquad \qquad \qquad\,\,\,\,\,
	-\frac{P_{0}}{\v}\int_{\Omega}(\theta_{0}-1+\v\theta_{1}+\v^2\tilde{\theta}_{R})\rho_{R}^{(k)}\,{\rm d}x:=\mathcal{R}_{4}(\rho_{R}^{(k)},\theta_{R}^{(k+1)}).	
\end{align}
It follows from the classical Stokes problem \cite[Theorems IV.5.2 and IV.5.8]{Boyer-Fabrie} that there exists a unique solution $(v_{R}^{(k+1)},P_{R}^{(k+1)})\in (H^3\cap H_{0}^1)\times H^2$ of \eqref{E8}--\eqref{E8-1} satisfying
\begin{align}\label{F8}
	&\|v_{R}^{(k+1)}\|_{H^2}+\|\nabla P_{R}^{(k+1)}\|_{L^2}\leq C\|\mathcal{R}_{3}(\rho_{R}^{(k)},w_{R}^{(k)},\theta_{R}^{(k)})\|_{L^2}+C|\nabla q_{R}^{(k+1)}|_{H^{\frac{3}{2}}(\partial \Omega)}\nonumber\\
	&\leq C(\mathcal{I}_{1})\|(\rho_{R}^{(k)},w_{R}^{(k)},\theta_{R}^{(k)})\|_{H^1}+C\|\mathfrak{R}_{2}\|_{L^2}+C(\mathcal{I}_{0})\|v_{R}^{(k)}\|_{H^2}+C\|q_{R}^{(k+1)}\|_{H^3}\nonumber\\
	&\leq C(\mathcal{I}_{1})\|(\rho_{R}^{(k)},w_{R}^{(k)},\theta_{R}^{(k)})\|_{H^1}+C(\mathcal{I}_{0})\|v_{R}^{(k)}\|_{H^2} +C(\mathcal{I}_{1})\v\|\rho_{R}^{(k)}\|_{H^2}\nonumber\\
	&\quad +C\big[\v\|\mathfrak{R}_{1}\|_{H^1}+\|\mathfrak{R}_{2}\|_{L^2}\big],
\end{align}
and
\begin{align}\label{F8-1}
	&\v\mu\|v_{R}^{(k+1)}\|_{H^3}+\v\|\nabla P_{R}^{(k+1)}\|_{H^1}\nonumber\\
	&\leq \bar{C}\v\|\mathcal{R}_{3}(\rho_{R}^{(k)},w_{R}^{(k)},\theta_{R}^{(k)})\|_{H^1}+\bar{C}\mu\v|\nabla q_{R}^{(k+1)}|_{H^{\frac{5}{2}}(\partial \Omega)}\nonumber\\
	&\leq C(\mathcal{I}_{1})\|(\rho_{R}^{(k)},w_{R}^{(k)},\theta_{R}^{(k)})\|_{H^2}+C\v\|\mathfrak{R}_{2}\|_{H^1}+C(\mathcal{I}_{0})\mu\v\|v_{R}^{(k)}\|_{H^3}+\bar{C}\mu\v\|q_{R}^{(k+1)}\|_{H^4}\nonumber\\
	& \leq C(\mathcal{I}_{1})\|(\rho_{R}^{(k)},w_{R}^{(k)},\theta_{R}^{(k)})\|_{H^2}+C(\mathcal{I}_{0})\mu\v\|v_{R}^{(k)}\|_{H^3}+C\v^2\|\mathfrak{R}_{1}\|_{H^2}\nonumber\\
	&\quad +\bar{C}\mu\v^2\|\nabla^2 \operatorname{div}[\rho_{R}^{(k)}(u_{1}+\v(u_{2}+\tilde{u}_{R}))]\|_{L^2}+C\v\|\mathfrak{R}_{2}\|_{H^1},
\end{align}
where we have used \eqref{F7}--\eqref{F7-1} and the fact
where we have used that fact
\begin{equation}\label{F8-3}
	\begin{aligned}
		&\|\tilde{F}^{\v}(\rho_{R},u_{R},\theta_{R})\|_{L^2}\leq \|\tilde{f}_{1}^{\v}(u_{R},\theta_{R})\|_{L^2}+\|\tilde{f}_{2}^{\v}(\rho_{R})\|_{L^2}\leq C(\mathcal{I}_{1})\|(\rho_{R},u_{R},\theta_{R})\|_{H^1},\\
		&\|\tilde{F}^{\v}(\rho_{R},u_{R},\theta_{R})\|_{H^1}\leq \|\tilde{f}_{1}^{\v}(u_{R},\theta_{R})\|_{H^1}+\|\tilde{f}_{2}^{\v}(\rho_{R})\|_{H^1}\leq C(\mathcal{I}_{1})\|(\rho_{R},u_{R},\theta_{R})\|_{H^2}.
	\end{aligned}
\end{equation}
Here,
$\bar{C}>0$ in \eqref{F8-1} is independent of $\mu$ and $\zeta$. Moreover, it follows from \eqref{F8} and the Poincar\'{e} inequality that
\begin{align}\label{F8-2}
	\v\|P_{R}^{(k+1)}\|_{L^2}&\leq C\v\|\nabla P_{R}^{(k+1)}\|_{L^2}+C\v|\mathcal{R}_{4}(\rho_{R}^{(k)},\theta_{R}^{(k+1)})|\nonumber\\
	&\leq C(\mathcal{I}_{1})\|(\rho_{R}^{(k)},w_{R}^{(k)},\theta_{R}^{(k)})\|_{H^1}+C(\mathcal{I}_{0})\v\|v_{R}^{(k)}\|_{H^2}\nonumber\\
	&\quad +C(\mathcal{I}_{1})\v^2\|\rho_{R}^{(k)}\|_{H^2}+C\big[\v^2\|\mathfrak{R}_{1}\|_{H^1}+\v\|\mathfrak{R}_{2}\|_{L^2}\big]+C\|\theta_{R}^{(k+1)}\|_{L^2}.
\end{align}
Define $w_{R}^{(k+1)}=:v_{R}^{(k+1)}+\nabla q_{R}^{(k+1)}$, one directly has $w_{R}^{(k+1)}\vert_{\partial \Omega}=0$.

{\it Step 3}. Finally, we construct $\rho_{R}^{(k+1)}$ by
\begin{align}\label{E9}
	&\rho_{R}^{(k+1)}+\frac{\varepsilon^2(\mu+\zeta)}{P_{0}}\operatorname{div}[\rho_{R}^{(k+1)}\theta_{0}(u_{1}+\varepsilon(u_{2}+\tilde{u}_{R}))]\nonumber\\
	&=- \frac{\var}{P_{0}}P_{R}^{(k+1)}+\frac{\v^2(\mu+\zeta)}{\rho_{0}}\mathfrak{R}_{1}-(\rho_{0}+\v\rho_{1}+\v^2\rho_{2}+\v^3\rho_{3})\theta_{R}^{(k+1)}\nonumber\\
	&\quad -(\theta_{0}-1+\var\theta_{1}+\var^2\tilde{\theta}_{R})\rho_{R}^{(k)}+\frac{\varepsilon^2(\mu+\zeta)}{P_{0}}\rho_{R}^{(k)}(u_{1}+\varepsilon(u_{2}+\tilde{u}_{R}))\cdot \nabla\theta_{0}\nonumber\\
	&=:-\frac{\v}{P_{0}}P_{R}^{(k+1)}-(\rho_{0}+\v\rho_{1}+\v^2\rho_{2}+\v^3\rho_{3})\theta_{R}^{(k+1)}+\mathcal{R}_{5}(\rho_{R}^{(k)}).
	%&\rho_{R}^{(k+1)}+\varepsilon^2(\mu+\zeta)\operatorname{div}[\rho_{R}^{(k+1)}(u_{1}+\varepsilon(u_{2}+\tilde{u}_{R}))]\nonumber\\
	%&=-\var P_{R}^{(k+1)}-\frac{\var^2(\mu+\zeta)}{\rho_{0}}\big\{\frac{1}{\varepsilon}u_{R}^{(k)}\cdot \nabla \rho_{0}+\operatorname{div}[\tilde{u}_{R}(\rho_{1}+\varepsilon\rho_{2}+\var^2\rho_{3})]+r_{1}\}\nonumber\\
	%&\quad -(\rho_{0}+\v\rho_{1}+\v^2\rho_{2}+\v^3\rho_{3})\theta_{R}^{(k)}-(\theta_{0}-1+\var\theta_{1}+\var^2\tilde{\theta}_{R})\rho_{R}^{(k)}.
\end{align}

From the regularity theory of the steady transport equation ({\it cf.} \cite[Lemma 3.5]{Choe-Jin}) and \eqref{E8-1}, we know if
\begin{align}\label{H1}
	\v^2(\mu+\zeta)\big[\|u_{1}\|_{H^3}+\v\|u_{2}\|_{H^3}+\v\|\tilde{u}_{R}\|_{H^3}\big]\ll 1,
\end{align}
then there exists a unique solution $\rho_{R}^{(k+1)}\in H^2$ with $\int_{\Omega}\rho_{R}^{(k+1)}\,{\rm d}x=0$ of \eqref{E9} satisfying
\begin{align}\label{F9}
	\|\rho_{R}^{(k+1)}\|_{H^2}&\leq C\v\|P_{R}^{(k+1)}\|_{H^2}+C\|\theta_{R}^{(k+1)}\|_{H^2}+C\v\|\mathcal{R}_{5}(\rho_{R}^{(k)})\|_{H^2}\nonumber\\
	&\leq  C(\mathcal{I}_{1})\|(\rho_{R}^{(k)},w_{R}^{(k)},\theta_{R}^{(k)})\|_{H^2}+C(\mathcal{I}_{0})\v\|v_{R}^{(k)}\|_{H^3}+C\big[\v^2\|\mathfrak{R}_{1}\|_{H^2}+\v\|\mathfrak{R}_{2}\|_{H^1}\big]\nonumber\\
	&\quad +\bar{C}\mu\v^2\|\nabla^2 \operatorname{div}[\rho_{R}^{(k)}(u_{1}+\v(u_{2}+\tilde{u}_{R}))]\|_{L^2}+C\|\theta_{R}^{(k+1)}\|_{H^2},
\end{align}
where we have used \eqref{F8}--\eqref{F8-2}.
%It is worth pointing out that the constant in \eqref{F9} depends only on $\mu$ and $\Omega$.
One also has from \eqref{E9} that
\begin{align}\label{F9-1}
	&\v^2(\mu+\zeta)\|\operatorname{div}[\rho_{R}^{(k+1)}\theta_{0}(u_{1}+\v(u_{2}+\tilde{u}_{R}))]\|_{H^2}\nonumber\\&\leq C\|\rho_{R}^{(k+1)}\|_{H^2}+ C\v\|P_{R}^{(k+1)}\|_{H^2}+C\|\theta_{R}^{(k+1)}\|_{H^2}+C\|\mathcal{R}_{5}(\rho_{R}^{(k)},\theta_{R}^{(k)})\|_{H^2}\nonumber\\
	&\leq C(\mathcal{I}_{1})\|(\rho_{R}^{(k)},w_{R}^{(k)},\theta_{R}^{(k)})\|_{H^2}+C(\mathcal{I}_{0})\mu\v\|v_{R}^{(k)}\|_{H^3}+C\big[\v^2\|\mathfrak{R}_{1}\|_{H^2}+\v\|\mathfrak{R}_{2}\|_{H^1}\big]\nonumber\\
	&\quad +\bar{C}\mu\v^2\|\nabla^2 \operatorname{div}[\rho_{R}^{(k)}(u_{1}+\v(u_{2}+\tilde{u}_{R}))]\|_{L^2}+C\|\theta_{R}^{(k+1)}\|_{H^2},
	%&\leq C(\mathcal{I}_{1})\|(\rho_{R}^{(k)},w_{R}^{(k)},\theta_{R}^{(k)})\|_{H^2}+C(\mathcal{I}_{0})\v\|v_{R}^{(k)}\|_{H^3}+C\mu\v^2\|\mathfrak{R}_{1}\|_{H^2}\nonumber\\
	%&\quad +\bar{C}\mu\v^2\|\nabla^2 \operatorname{div}[\rho_{R}^{(k)}(u_{1}+\v(u_{2}+\tilde{u}_{R}))]\|_{L^2}+C\|\theta_{R}^{(k)}\|_{H^2},
\end{align}
which leads to
\begin{align}\label{F9-2}
	&\v^2(\mu+\zeta)\|\nabla^2\operatorname{div}[\rho_{R}^{(k+1)}(u_{1}+\v(u_{2}+\tilde{u}_{R}))]\|_{L^2}\nonumber\\
	&\leq \v^2(\mu+\zeta)\|\nabla^2\operatorname{div}[\rho_{R}^{(k+1)}\theta_{0}(u_{1}+\v(u_{2}+\tilde{u}_{R}))]\|_{L^2}+C(\mathcal{I}_{1})\v^2\||\rho_{R}^{(k+1)}\|_{H^2}
	\nonumber\\
	&\leq C(\mathcal{I}_{1})\|(\rho_{R}^{(k)},w_{R}^{(k)},\theta_{R}^{(k)})\|_{H^2}+C(\mathcal{I}_{0})\mu\v\|v_{R}^{(k)}\|_{H^3}+C\big[\v^2\|\mathfrak{R}_{1}\|_{H^2}+\v\|\mathfrak{R}_{2}\|_{H^1}\big]\nonumber\\
	&\quad +\bar{C}\mu\v^2\|\nabla^2 \operatorname{div}[\rho_{R}^{(k)}(u_{1}+\v(u_{2}+\tilde{u}_{R}))]\|_{L^2}+C\|\theta_{R}^{(k+1)}\|_{H^2}.
\end{align}

%\begin{align}\label{F9-3}
%		&\v^2(\mu+\zeta)\|\nabla^2\operatorname{div}[\rho_{R}^{(k+1)}(u_{1}+\v(u_{2}+\tilde{u}_{R}))]\|_{L^2}\nonumber\\
%		&\leq C(\mathcal{I}_{1})\|(\rho_{R}^{(k)},w_{R}^{(k)},\theta_{R}^{(k)})\|_{H^2}+C(\mathcal{I}_{0})\v\|v_{R}^{(k)}\|_{H^3}+C\|\theta_{R}^{(k+1)}\|_{H^2}\nonumber\\
%		&\quad +C\big[\v^2\|\mathfrak{R}_{1}\|_{H^2}+\v\|\mathfrak{R}_{2}\|_{H^1}\big]
%		%&\leq \v^2\|\nabla^2\operatorname{div}[\rho_{R}^{(k+1)}\theta_{0}(u_{1}+\v(u_{2}+\tilde{u}_{R}))]\|_{L^2}+C(\mathcal{I}_{1})\||\rho_{R}^{(k+1)}\|_{H^2}\nonumber\\
%		%&\leq C\|\rho_{R}^{(k+1)}\|_{H^2}+\v C\|P_{R}^{(k+1)}\|_{H^2}+C\|R_{5}(\rho_{R}^{(k)},\theta_{R}^{(k)})\|_{H^2}\nonumber\\
%		%&\leq C(\mathcal{I}_{1})\|(\rho_{R}^{(k)},w_{R}^{(k)},\theta_{R}^{(k)})\|_{H^2}+C(\mathcal{I}_{0})\|v_{R}^{(k)}\|_{H^3}+C(\mathcal{I}_{0})(1+\v\|\tilde{u}_{R}\|_{H^3})+C\|\theta_{R}^{(k)}\|_{H^2}.
%\end{align}
Combining \eqref{F6}--\eqref{F7-1}, \eqref{F8}--\eqref{F8-2}, \eqref{F9}--\eqref{F9-2}, one concludes
\begin{align}\label{F10}
	&\|(\rho_{R}^{(k+1)},v_{R}^{(k+1)},w_{R}^{(k+1)})\|_{H^2}+\|(\theta_{R}^{(k+1)},q_{R}^{(k+1)})\|_{H^3}+\v\mu\|(w_{R}^{(k+1)},v_{R}^{(k+1)})\|_{H^3}\nonumber\\
	&\quad +\v\mu\|q_{R}^{(k+1)}\|_{H^4}+\v^2\|\nabla^2 \operatorname{div}[\rho_{R}^{(k+1)}(u_{1}+\v(u_{2}+\tilde{u}_{R}))]\|_{L^2}\nonumber\\
	&\leq C(\mathcal{I}_{1})\|(\rho_{R}^{(k)},w_{R}^{(k)},\theta_{R}^{(k})\|_{H^2}+C(\mathcal{I}_{0})\v\mu\|v_{R}^{(k)}\|_{H^3}\nonumber\\
	&\quad +\bar{C}\frac{\mu}{\mu+\zeta}\v^2\|\nabla^2 \operatorname{div}[\rho_{R}^{(k)}(u_{1}+\v(u_{2}+\tilde{u}_{R}))]\|_{L^2}\nonumber\\
	&\quad +C\big[\|\mathfrak{R}_{2}\|_{L^2}+\|(\mathfrak{R}_{1},\mathfrak{R}_{3})\|_{H^1}\big]+C\v\big[\|\mathfrak{R}_{2}\|_{H^1}+\v\|\mathfrak{R}_{1}\|_{H^2}\big].
\end{align}
provided that $\mu+\zeta<1$ and \eqref{H1} hold. Noting $\bar{C}$ is independent of $\mu$ and $\zeta$, we can choose two positive constants $\mu_{0}$ and $\zeta_{0}$ satisfying $\mu_{0}+\zeta_{0}<1$ and  $\frac{\zeta_{0}}{\mu_{0}}\gg 1$ such that $\bar{C}\frac{\mu_{0}}{\mu_{0}+\zeta_{0}}\leq \frac{1}{2}$ in \eqref{F10}. Furthermore, noting \eqref{I0}--\eqref{I1}, we can take $\delta_{0}$ and $\|\tilde{u}_{R}\|_{\mathcal{K}}+\|\tilde{\theta}_{R}\|_{H^3}$ small enough such that \eqref{H1} holds and $C(\mathcal{I}_{0}),C(\mathcal{I}_{1})$ in \eqref{F10} satisfy $C(\mathcal{I}_{i})\leq \frac{1}{2}\,\,(i=0,1)$. Then, by iterative argument, we can derive from \eqref{F10} that
\begin{align}\label{F10-1}
	&\|(\rho_{R}^{(k)},v_{R}^{(k)},w_{R}^{(k)})\|_{H^2}+\|(\theta_{R}^{(k)},q_{R}^{(k)})\|_{H^3}+\v\mu_{0}\|(w_{R}^{(k)},v_{R}^{(k)})\|_{H^3}\nonumber\\
	&\quad +\v\mu_{0}\|q_{R}^{(k)}\|_{H^4}+\v^2\|\nabla^2 \operatorname{div}[\rho_{R}^{(k)}(u_{1}+\v(u_{2}+\tilde{u}_{R}))]\|_{L^2}\nonumber\\
	&\leq C\big[\|\mathfrak{R}_{2}\|_{L^2}+\|(\mathfrak{R}_{1},\mathfrak{R}_{3})\|_{H^1}\big]+C\v\big[\|\mathfrak{R}_{2}\|_{H^1}+\v\|\mathfrak{R}_{1}\|_{H^2}\big]
\end{align}
where $C>0$ is a constant depending only on $\mu_{0}$ and $\zeta_{0}$.

Noting \eqref{E6}--\eqref{E7}, \eqref{E8} and \eqref{E9}, then by similar arguments, we can prove
%$$
%(\rho_{R}^{(k+1)}-\rho_{R}^{(k)},v_{R}^{(k+1)}-v_{R}^{(k)},q_{R}^{(k+1)}-q_{R}^{(k)},\theta_{R}^{(k+1)}-\theta_{R}^{(k)})
%$$
$(\rho_{R}^{(k)},w_{R}^{(k)},\theta_{R}^{(k)})$ is a Cauchy sequence in $H^2\times \mathcal{K}\times (H_{0}^1\cap H^3)$ for $(\mu_{0},\zeta_{0})$. It is easy to prove that
$$
(\rho_{R},u_{R},\theta_{R}):=\lim\limits_{k\to \infty}(\rho_{R}^{(k)},u_{R}^{(k)},\theta_{R}^{(k)})=\lim\limits_{k\to \infty}(\rho_{R}^{(k)},\frac{1}{\rho_{0}}w_{R}^{(k)},\theta_{R}^{(k)})
$$
is a strong solution of \eqref{3.10} satisfying \eqref{g1} with $(\mu,\zeta)=(\mu_{0},\zeta_{0})$. Therefore, the proof of Lemma \ref{lem7.1} is complete. $\hfill\square$

\section{{\it A Priori} Uniform Estimates in Mach Number for \eqref{3.10}}
As mentioned in Remark \ref{rem3.1}, in this section, we shall show any strong solution $(\rho_{R},u_{R},\theta_{R})\in H^2\times \mathcal{K}\times (H_{0}^1\cap H^3)$ of \eqref{3.10} enjoys {\it a priori} uniform estimates in $\v$, and is of continuous dependence on $(\mu,\zeta)$, which will play an essential role in establishing the existence of strong solution of \eqref{3.10} for any given viscous coefficients pair $(\mu,\zeta)$ in Section 5. Hereafter, we assume that $(\rho_{R},u_{R},\theta_{R})\in H^2\times \mathcal{K}\times (H_{0}^1\cap H^3)$ is a strong solution of \eqref{3.10}.

\subsection{Lower order estimates}
Noting that Dirichlet boundary conditions of $(u_{R},\theta_{R})$ the singular terms (the terms involving $\v^{-1}$) is skew-symmetric in \eqref{3.10}, we can apply the standard the energy method to get following lemma:
\begin{lemma}\label{lem2.2}
For any $0<\v<1$, there exist a small constant $\tilde{\delta}_{0}$ such that if $\delta_{0}\leq \tilde{\delta}_{0}$ and  $\|\tilde{u}_{R}\|_{\mathcal{K}}+\|\tilde{\theta}_{R}\|_{H^3}\leq \tilde{\delta}_{0}$, there exists a positive constant $C$ independent of $\varepsilon$ such that
	\begin{align}\label{4.1}
		\|\rho_{R}\|_{L^2}^2+\|(u_{R},\theta_{R})\|_{H^1}^2
		&\leq C\|(\mathfrak{R}_{1},\mathfrak{R}_{2},\mathfrak{R}_{3})\|_{L^2}.
	\end{align}
\end{lemma}

\noindent\textbf{Proof.} Multiplying $\eqref{3.10}_{1}$ by $\theta_{0}\rho_{R}$, $\eqref{3.10}_{2}$ by $\rho_{0}u_{R}$ and $\eqref{3.10}_{3}$ by $\rho_{0}\theta_{R}$, then summing the resultant equations, one has from the H\"{o}lder and Poincar\'{e} inequalities that
\begin{align}\label{3.15-1}
	&\mu\int_{\Omega}\rho_{0}|\nabla u_{R}|^2\,{\rm d}x+\zeta\int_{\Omega}\rho_{0}|\operatorname{div}u_{R}|^2\,{\rm d}x+\kappa\int_{\Omega}\frac{\rho_{0}}{\theta_{0}}|\nabla \theta_{R}|^2\,{\rm d}x\nonumber\\
    &\leq \big[C(\mathcal{I}_{0})+\frac{1}{16}\big]\Big[\mu\int_{\Omega}\rho_{0}|\nabla u_{R}|^2\,{\rm d}x+\kappa\int_{\Omega}\frac{\rho_{0}}{\theta_{0}}|\nabla \theta_{R}|^2\,{\rm d}x\Big]+C(\mathcal{I}_{0})\|\rho_{R}\|_{L^2}^2+C\|(\mathfrak{R}_{2},\mathfrak{R}_{3})\|_{L^2}^2\nonumber\\
    &\quad +C\|\rho_{R}\|_{L^2}\|\mathfrak{R}_{1}\|_{L^2}-\frac{1}{\var}\int_{\Omega}\rho_{0}u_{R}\cdot\nabla(\rho_{0}\theta_{R}+\rho_{R}\theta_{0})\,{\rm d}x-\frac{1}{\v}\int_{\Omega}\operatorname{div}(\rho_{0}u_{R})(\rho_{0}\theta_{R}+\theta_{0}\rho_{R})\,{\rm d}x\nonumber\\
    &\quad -\int_{\Omega}\rho_{0}u_{R}\cdot\nabla (\rho_{R}\theta_{1})\,{\rm d}x -\int_{\Omega}\rho_{0}u_{R}\cdot\tilde{F}^{\v}(\rho_{R},u_{R},\theta_{R})\,{\rm d}x-\frac{1}{\theta_{0}}\int_{\Omega}\rho_{0}\theta_{R}\tilde{G}^{\v}(\rho_{R},u_{R},\theta_{R})\,{\rm d}x\nonumber\\
    &\quad -\int_{\Omega}\theta_{0}\rho_{R}\operatorname{div}\big[\rho_{R}(u_{1}+\v(u_{2}+\tilde{u}_{R}))\big]\,{\rm d}x=:\sum\limits_{i=1}^{10}I_{i}.
\end{align}
For $I_{5}+I_{6}$, by integrating by parts, one has
\begin{align}\label{3.15-2}
	I_{5}+I_{6}=\frac{1}{\v}\int_{\Omega}\operatorname{div}(\rho_{0}u_{R})(\rho_{0}\theta_{R}+\theta_{0}\rho_{R})\,{\rm d}x-\frac{1}{\v}\int_{\Omega}\operatorname{div}(\rho_{0}u_{R})(\rho_{0}\theta_{R}+\theta_{0}\rho_{R})\,{\rm d}x=0.
\end{align}
Similarly, for $I_{7}$, a direct calculation from the H\"{o}lder and Poincar\'{e} inequalities shows that
\begin{align}\label{3.15-3}
I_{7}=\int_{\Omega}\rho_{R}\theta_{1}\operatorname{div}(\rho_{0}u_{R})\,{\rm d}x\leq C(\mathcal{I}_{0})\|\rho_{R}\|_{L^2}^2+C(\mathcal{I}_{0})\mu\int_{\Omega}\rho_{0}|\nabla u_{R}|^2\,{\rm d}x.
\end{align}
For $I_{8}$ and $I_{9}$, noting \eqref{4.5}--\eqref{4.8}, one has from the H\"{o}lder and Poincar\'{e} inequalities that
\begin{align}\label{3.15-4}
	I_{8}+I_{9}&=-\int_{\Omega}\rho_{0}u_{R}\cdot f_{1}^{\v}(u_{R},\theta_{R})\,{\rm d}x-\int_{\Omega}\frac{\rho_{0}}{\theta_{0}}\theta_{R}g_{1}^{\v}(u_{R},\theta_{R})\,{\rm d}x\nonumber\\
	&\quad-\int_{\Omega}\rho_{0}u_{R}\cdot f_{2}^{\v}(\rho_{R})\,{\rm d}x-\int_{\Omega}\frac{\rho_{0}}{\theta_{0}}\theta_{R}g_{2}^{\v}(\rho_{R})\,{\rm d}x\nonumber\\
	&\leq C(\mathcal{I}_{1})\Big[\mu\int_{\Omega}\rho_{0}|\nabla u_{R}|^2\,{\rm d}x+\kappa\int_{\Omega}\frac{\rho_{0}}{\theta_{0}}|\nabla \theta_{R}|^2\,{\rm d}x\Big]+C(\mathcal{I}_{1})\|\rho_{R}\|_{L^2}^2,
\end{align}
where we have used the fact
\begin{align*}
	\int_{\Omega}\rho_{0}u_{R}\cdot \nabla(\rho_{R}\tilde{\theta}_{R})\,{\rm d}x=-\int_{\Omega}\operatorname{div}{\rho_{0}u_{R}}(\rho_{R}\tilde{\theta}_{R})\,{\rm d}x\leq C(\mathcal{I}_{1})\|\rho_{R}\|_{L^2}^2+C(\mathcal{I}_{1})\mu\int_{\Omega}\rho_{0}|\nabla u_{R}|^2\,{\rm d}x.
\end{align*}
For $I_{10}$, by integrating by parts, one has
\begin{align}\label{3.15-5}
I_{10}&\leq -\frac{1}{2}\int_{\Omega}\theta_{0}[u_{1}+\v(u_{2}+\tilde{u}_{R})]\cdot \nabla (\rho_{R}^2)\,{\rm d}x+C(\mathcal{I}_{1})\|\rho_{R}\|_{L^2}^2\nonumber\\
&= \frac{1}{2}\int_{\Omega}\rho_{R}^2\operatorname{div}[\theta_{0}(u_{1}+\v(u_{2}+\tilde{u}_{R}))]\,{\rm d}x+C(\mathcal{I}_{1})\|\rho_{R}\|_{L^2}^2\nonumber\\
&\leq C(\mathcal{I}_{1})\|\rho_{R}\|_{L^2}^2.
\end{align}
Substituting \eqref{3.15-2}--\eqref{3.15-5} into \eqref{3.15-1} and taking $\delta_{0}$ and $\|\tilde{u}_{R}\|_{\mathcal{K}}+\|\tilde{\theta}_{R}\|_{H^3}$ small enough, we obtain
\begin{align}\label{3.15-6}
\|\nabla u_{R}\|_{L^2}^2+\|\nabla \theta_{R}\|_{L^2}^2\leq C(\mathcal{I}_{1})\|\rho_{R}\|_{L^2}^2+C\|\rho_{R}\|_{L^2}\|\mathfrak{R}_{1}\|_{L^2}+C\|(\mathfrak{R}_{2},\mathfrak{R}_{3})\|_{L^2}^2.	
\end{align}

It remains to bound $\|\rho_{R}\|_{L^2}$. For any function $\xi\in L^2(\Omega)$ with $\int_{\Omega}\xi\,{\rm d}x=0$, there exists a vector field $q\in H_{0}^{1}(\Omega)$ such that
\begin{align}\label{q}
	\operatorname{div}q=\xi\qquad \|q\|_{H^1}\leq C(\Omega)\|\xi\|_{L^2}.
\end{align}
Multiplying $\eqref{3.10}_{2}$ by $q$ and integrating the resulting equation over $\Omega$, we obtain
\begin{align}\label{3.26}
	\frac{1}{\varepsilon}\int_{\Omega}\theta_{0}\rho_{R}\operatorname{div}q\,{\rm d}x&=-\frac{1}{\varepsilon}\int_{\Omega}\rho_{0}\theta_{R}\operatorname{div}q\,{\rm d}x+\mu\int_{\Omega}\nabla u_{R}\cdot \nabla q\,{\rm d}x+\zeta\int_{\Omega}\operatorname{div}u_{R}\cdot \operatorname{div}q\,{\rm d}x\nonumber\\
	&\quad +\int_{\Omega}\rho_{0}(u_{1}\cdot \nabla u_{R}+u_{R}\cdot \nabla u_{1})\cdot q\,{\rm d}x-\int_{\Omega}\rho_{R}\theta_{1}\operatorname{div}q\,{\rm d}x-\int_{\Omega}\rho_{1}\theta_{R}\operatorname{div}q\,{\rm d}x\nonumber\\
	&\quad +\int_{\Omega}\tilde{F}^{\varepsilon}(\rho_{R},u_{R},\theta_{R})\cdot q\,{\rm d}x+\int_{\Omega}\mathfrak{R}_{2}q\,{\rm d}x.
\end{align}
By taking $\xi=\rho_{R}$, we obtain from \eqref{q}--\eqref{3.26}, \eqref{4.5}--\eqref{4.6} and the H\"{o}lder inequality that
\begin{align}\label{3.27}
	\int_{\Omega}\theta_{0}|\rho_{R}|^2\,{\rm d}x&\leq C(\mathcal{I}_{1})\int_{\Omega}\theta_{0}|\rho_{R}|^2\,{\rm d}x+C\big[\|\theta_{R}\|_{L^2}^2+\v^2\|\nabla u_{R}\|_{L^2}^2+\v^2\|\mathfrak{R}_{2}\|_{L^2}^2\big],
\end{align}
where we have used the fact that
$$
-\varepsilon^2\int_{\Omega}\nabla (\rho_{R}\tilde{\theta}_{R})\cdot q\,{\rm d}x=\varepsilon^2\int_{\Omega}\rho_{R}\tilde{\theta}_{R}\operatorname{div}q\,{\rm d}x\leq C\varepsilon^2\|\tilde{\theta}_{R}\|_{H^2}\int_{\Omega}\theta_{0}|\rho_{R}|^2\,{\rm d}x.
$$
Noting \eqref{I1} and taking $\delta_{0}$ and $\|\tilde{u}_{R}\|_{\mathcal{K}}+\|\tilde{\theta}_{R}\|_{H^3}$ small enough, we have from \eqref{3.27} that
\begin{align}\label{3.27-1}
	\|\rho_{R}\|_{L^2}^2\leq C\big[\|\theta_{R}\|_{L^2}^2+\v^2\|\nabla u_{R}\|_{L^2}^2+\v^2\|\mathfrak{R}_{2}\|_{L^2}^2\big].	
\end{align}
Combining \eqref{3.15-6} and \eqref{3.27-1}, using the H\"{o}lder and Poincar\'{e} inequalities, and taking $\delta_{0}$ and $\|\tilde{u}_{R}\|_{\mathcal{K}}+\|\tilde{\theta}_{R}\|_{H^3}$ small enough, we get \eqref{4.1}. Therefore the proof of Lemma \ref{lem2.2} is complete. $\hfill\square$

As a byproduct of \eqref{4.1}, one can obtain the uniqueness of the strong solution $(\rho_{R},u_{R},\theta_{R})$ of \eqref{3.10}. Indeed, assuming $(\rho_{R}^{(k)},u_{R}^{(k)},\theta_{R}^{(k)})\in H^2\times \mathcal{K}\times (H_{0}^1\cap H^3)\,\,(k=1,2)$ are two strong solutions of \eqref{3.10}, one can apply similar arguments as in the proof of Lemma \ref{lem2.2} to
$$
(\bar{\rho}_{R},\bar{u}_{R},\bar{\theta}_{R})=:(\rho_{R}^{(1)}-\rho_{R}^{(2)},u_{R}^{(1)}-u_{R}^{(2)},\theta_{R}^{(1)}-\theta_{R}^{(2)})
$$
to get
$$
\|\bar{\rho}_{R}\|_{L^2}^2+\|(\bar{u}_{R},\bar{\theta}_{R})\|_{H^1}^2\leq 0,
$$
which yields the uniqueness. That is, we have following corollary.

\medskip

\begin{corollary}\label{cor3.1}
	The strong solution $(\rho_{R},u_{R},\theta_{R})\in H^2\times \mathcal{K}\times (H_{0}^1\cap H^3)$ of \eqref{3.10} is unique.
\end{corollary}

\subsection{Higher order estimates} We rewrite the momentum equation $\eqref{3.10}_{2}$ as non-homogeneous Stokes problem:
\begin{equation}\label{4.21}
\left\{
\begin{aligned}
&-\mu\Delta u_{R}+\frac{\nabla (\rho_{0}\theta_{R}+\theta_{0}\rho_{R})}{\varepsilon}=-\rho_{0}(u_{1}\cdot \nabla u_{R}+u_{R}\cdot \nabla u_{1})-\nabla (\rho_{R}\theta_{1}+\rho_{1}\theta_{R})+\zeta\nabla \operatorname{div}u_{R}\\
&\qquad \qquad \qquad \qquad \qquad\qquad \quad\,\,\, -\tilde{F}^{\varepsilon}(\rho_{R},u_{R},\theta_{R})-\mathfrak{R}_{2},\\
&\operatorname{div}u_{R}=\operatorname{div}u_{R},\\
&u_{R}\vert_{\partial\Omega}=0.
\end{aligned}
\right.
\end{equation}
By the classical Stokes estimates \cite[Theorem IV.5.8]{Boyer-Fabrie}, we have
\begin{align}\label{4.22}
\|u_{R}\|_{H^2}+\|\frac{\nabla (\rho_{0}\theta_{R}+\theta_{0}\rho_{R})}{\varepsilon}\|_{L^2}&\leq C\big[\|u_{1}\|_{H^2}\|u_{R}\|_{H^1}+\|\rho_{R}\|_{H^{1}}\|\theta_{1}\|_{H^2}+\|\rho_{1}\|_{H^2}\|\theta_{R}\|_{H^1}\big] \nonumber\\
&\quad +C\big[\|\tilde{F}\|_{L^2}+\|\mathfrak{R}_{2}\|_{L^2}+\|\nabla\operatorname{div}u_{R}\|_{L^2}+\|u_{R}\|_{H^1}\big],
\end{align}
and
\begin{align}\label{4.23}
	\|u_{R}\|_{H^3}+\|\frac{\nabla (\rho_{0}\theta_{R}+\theta_{0}\rho_{R})}{\varepsilon}\|_{H^1}&\leq C\big[\|u_{1}\|_{H^2}\|u_{R}\|_{H^2}+\|\rho_{R}\|_{H^2}\|\theta_{1}\|_{H^2}+\|\rho_{1}\|_{H^2}\|\theta_{R}\|_{H^2}\big]\nonumber\\
	&\quad  +C\big[\|\tilde{F}\|_{H^1}+\|\mathfrak{R}_{2}\|_{H^1}+\|\nabla^2\operatorname{div}u_{R}\|_{L^2}+\|u_{R}\|_{H^1}\big],
\end{align}
where we have used the following interpolation inequality:
$$
\|u_{R}\|_{H^2}\leq \delta\|u_{R}\|_{H^3}+C_{\delta}\|u\|_{H^1}\qquad \text{for any }\delta>0.
$$
Taking $\mathcal{I}_{1}\ll 1$, we deduce from \eqref{F8-3}, \eqref{4.22}--\eqref{4.23} that
\begin{align}\label{4.24-0}
&\|u_{R}\|_{H^2}+\|\frac{\nabla(\rho_{0}\theta_{R}+\theta_{0}\rho_{R})}{\varepsilon}\|_{L^2}\nonumber\\
&\leq C(\mathcal{I}_{1})\|(\rho_{R},\theta_{R})\|_{H^1}+C\big[\|\mathfrak{R}_{2}\|_{L^2}+\|\nabla\operatorname{div}u_{R}\|_{L^2}+\|u_{R}\|_{H^1}\big],
\end{align}
and
\begin{align}\label{4.24}
&\varepsilon\|u_{R}\|_{H^3}+\|\rho_{R}\|_{H^2}\leq C\big[\|\theta_{R}\|_{H^2}+\|\rho_{R}\|_{H^1}+\v\|u_{R}\|_{H^1}+\v\|\nabla^2\operatorname{div}u_{R}\|_{L^2}+\v\|\mathfrak{R}_{2}\|_{H^1}\big],
\end{align}
where we have used the fact
$$
\|\nabla (\theta_{0}\rho_{R}+\rho_{0}\theta_{R})\|_{H^1}\geq \frac{1}{C}\|\rho_{R}\|_{H^2}-\|\theta_{R}\|_{H^2}-\|\rho_{R}\|_{H^1}.
$$

To close the estimate \eqref{4.24-0}--\eqref{4.24}, we still need to estimate $\|\theta_{R}\|_{H^2}$ and $\|\nabla \operatorname{div} u_{R}\|_{L^2}+\varepsilon\|\nabla^2\operatorname{div}u_{R}\|_{L^2}$. As indicated in \cite{Dou-Jiang-Jiang-Yang}, the most difficult part is to estimate $\|\nabla \operatorname{div}u_{R}\|_{L^2}+\varepsilon\|\nabla^2\operatorname{div}u_{R}\|_{L^2}$, which should be divided into the interior part and the part near the boundary.

For later use, we introduce some notations. Since $\partial\Omega$ is compact, as in \cite{Choe-Jin,Li-Liao-2019}, there are $N$ cubes $Q_{1},\cdots, Q_{N}$ such that
\begin{itemize}
	\item $\partial\Omega\subset \cup_{k=1}^N \frac{1}{4}Q_{k}$, where $\frac{1}{4}Q_{k}$ means the cube with $\frac{1}{4}$ times the size of the radius of $Q_{k}$ and with the same center as $Q_{k}$;
	\item For each $Q_{k}$, $Q_{k}$ intersects at most eight other $Q_{j}'s$;
	\item If $w_{k}\in C^4(\partial\Omega\cap Q_{k})$ is the boundary function, that is,
	\begin{align*}
		&\Omega\cap Q_{k}=\{x=(x_{1},x_{2},x_{3})\,|\,x_{3}>w_{k}(x_{1},x_{2})\}\cap Q_{k},\\
		&\partial\Omega\cap Q_{k}=\{x_{3}=w_{k}(x_{1},x_{2})\}\cap Q_{k},
	\end{align*}
then
\begin{equation}\label{5.1}
	\|\nabla w_{k}\|_{C^3(\partial\Omega\cap Q_{k})}\leq \sigma\quad \text{for some small  constant }\sigma\ll 1.
\end{equation}
\end{itemize}
Let $\Omega_{0}:=\Omega- \cup_{k=1}^{N}\frac{1}{2}Q_{k}$, then
\begin{align}\label{covering}
\Omega\subset \Omega_{0}\cup(\cup_{k=1}^{N}\frac{1}{2}Q_{k}).
\end{align}
%are all balls with radius $q_{0}$ and $q_{k}$ respectively,  $\bar{Q}_{0}\subset\Omega$ and $\partial \Omega\subset \cup_{k=1}^N\frac{1}{2}Q_{k}$ (here $\frac{1}{2}Q_{k}$ is the ball with radius $\frac{q_{k}}{2}$ and the same center). Moreover, in each $Q_{k}$, there exists a function $w_{k}\in C^3(\R^2)$ such that
%Without loss of generality, we assume $Q_{k}$ is small enough such that
\subsubsection{Interior estimate} We focus on the interior estimate of $\|\nabla\operatorname{div}u_{R}\|_{L^2}+\varepsilon\|\nabla^2\operatorname{div}u_{R}\|_{L^2}$. Let $\chi_{0}$ be a $C_{0}^{\infty}$ function with $\Omega_{0}\subset \operatorname{supp}\chi_{0}\subset \Omega_{0}'$ satisfying
$$
\chi_{0}(x)\equiv 1\quad \text{for $x\in \Omega_{0}$}\quad \text{and}\quad \chi_{0}(x)\equiv 0\quad \text{for $x\in \overline{\Omega_{0}'}^{c}$}
$$
where $\Omega_{0}'=\Omega-\cup_{k=1}^N\frac{1}{4}Q_{k}$.
\begin{lemma}[Interior estimate of $\|\nabla^2 u_{R}\|_{L^2}$]\label{lem4.2}
%Let $(\rho_{R},u_{R},\theta_{R})$ be the solution of \eqref{3.10}. Then
For any $\tau>0$, there exists a positive constants $C_{\tau}$, which depends on $\tau$ but is independent of $\v$, such that
\begin{align}\label{5.1-1}
	&\mu\|\chi_{0}\sqrt{\rho_{0}}\nabla^2u_{R}\|_{L^2}^2+\zeta\|\chi_{0}\sqrt{\rho_{0}}\nabla\operatorname{div}u_{R}\|_{L^2}^2+\kappa\|\chi_{0}\sqrt{\frac{\rho_{0}}{\theta_{0}}}\nabla^2\theta_{R}\|_{L^2}^2\nonumber\\
	&\leq \tau\|\frac{\nabla(\rho_{0}\theta_{R}+\rho_{R}\theta_{0})}{\varepsilon}\|_{L^2}^2+C_{\tau}\|u_{R}\|_{H^1}^2+C\|\theta_{R}\|_{H^1}^2\nonumber\\
	&\quad +C(\mathcal{I}_{1})\|\rho_{R}\|_{H^2}^2+C\|\rho_{R}\|_{H^1}\|\mathfrak{R}_{1}\|_{H^1}+C\|(\mathfrak{R}_{2},\mathfrak{R}_{3})\|_{L^2}^2,
\end{align}
where $C>0$ is a constant independent of $\tau$ and $\v$.
\end{lemma}

\noindent\textbf{Proof.} We apply $\nabla$ to \eqref{3.10} to get that
\begin{equation}\label{4.25}
\left\{\begin{aligned}
	&\nabla\operatorname{div}[\rho_{R}(u_{1}+\v(u_{2}+\tilde{u}_{R}))]=-\frac{\nabla \operatorname{div}(\rho_{0}u_{R})}{\v}+\nabla\mathfrak{R}_{1},\\
	%&\partial_{ji}^2(\rho_{R}u_{1}^{i})=-\frac{1}{\varepsilon}(\partial_{ji}^2(\rho_{0}u_{R}^{i}))-\partial_{ji}^2(\rho_{1}\tilde{u}_{R}^{i})-\varepsilon^2\partial_{ji}^2(\rho_{3}\tilde{u}_{R}^{i})\\
	%&\qquad\qquad\quad\,\,\, -\varepsilon\partial_{ji}^2(\rho_{R}u_{2}^{i}+\rho_{2}\tilde{u}_{R}^{i}+\rho_{R}\tilde{u}_{R}^{i})+\partial_{j}r_{1},\\
%&-\mu\partial_{jii}^3u_{R}^{k}-\zeta\partial_{jki}^3u_{R}^{i}=-\frac{\partial_{jk}^2(\rho_{0}\theta_{R}+\theta_{0}\rho_{R})}{\varepsilon}-\partial_{j}[\rho_{0}(u_{1}^{i}\partial_{i}u_{R}^{k}+u_{R}^{i}\partial_{i}u_{1}^{i})]\\
&-\mu \Delta \nabla u_{R}^{k}-\zeta \nabla (\partial_{k}\operatorname{div}u_{R})=-\frac{\nabla \partial_{k}(\rho_{0}\theta_{R}+\theta_{0}\rho_{R})}{\v}-\nabla[\rho_{0}(u_{1}\cdot \nabla )u_{R}^{k}+(u_{R}\cdot \nabla)u_{1}^{i}]\\
&\qquad\qquad\qquad\qquad\qquad\qquad\quad\,-\nabla \partial_{k}(\rho_{R}\theta_{1}+\rho_{1}\theta_{R})-\nabla \tilde{F}^{\varepsilon,k}-\nabla \mathfrak{R}_{2}^{k}\quad (k=1,2,3),\\
%&\qquad\qquad\qquad\qquad\qquad\,\, -\partial_{jk}^2(\rho_{R}\theta_{1})-\partial_{jk}^2(\rho_{1}\theta_{R})-\partial_{j}\tilde{F}^{\varepsilon,k}(\rho_{R},u_{R},\theta_{R})+\partial_{j}r_{2,k},\\
&-\frac{\kappa}{\theta_{0}}\Delta \nabla \theta_{R}=-\frac{\nabla \operatorname{div}(\rho_{0}u_{R})}{\v}-\frac{\kappa\nabla \theta_{0}}{\theta_{0}^2}\Delta \theta_{R}-P_{1}\nabla \operatorname{div}(\theta_{0}^{-1}\operatorname{div}u_{R})\\
%&-\frac{\kappa}{\theta_{0}}\partial_{jii}^3\theta_{R}=-\frac{\partial_{ji}^2(\rho_{0}u_{R}^{i})}{\varepsilon}-\frac{\kappa\partial_{j}\theta_{0}}{\theta_{0}^2}\partial_{ii}^2\theta_{R}-P_{1}\partial_{j}(\theta_{0}^{-1}\partial_{i}u_{R}^{i})-\partial_{j}[\theta_{0}^{-1}(\rho_{0}\theta_{R}+\rho_{R}\theta_{0})\partial_{i}u_{1}^{i}]\\
&\qquad\qquad\qquad\,\,-\nabla \{\theta_{0}^{-1}[\rho_{0}(u_{R}\cdot \nabla \theta_{1}+u_{1}\nabla \theta_{R})+\rho_{1}u_{R}\cdot \nabla \theta_{0}+\rho_{R}u_{1}\nabla \theta_{0}]\}\\
%&\qquad\qquad\qquad\,\, -\partial_{j}(\theta_{0}^{-1}\tilde{G}^{\varepsilon}(\rho_{R},u_{R},\theta_{R}))-\partial_{j}(\theta_{0}^{-1}r_{3})-\partial_{j}(\theta_{0}^{-1}\Psi(\nabla (u_{1}+\varepsilon(u_{2}+\tilde{u}_{R}))))\\
&\qquad\qquad\qquad\,\,
-\nabla [\theta_{0}^{-1}(\rho_{0}\theta_{R}+\rho_{R}\theta_{0})\operatorname{div}u_{1}]-\nabla (\theta_{0}^{-1}\tilde{G}^{\v})-\nabla \mathfrak{R}_{3}.
%&\qquad\qquad\qquad\,\, -\partial_{j}\{\theta_{0}^{-1}(\rho_{0}(u_{R}^{i}\partial_{i}\theta_{1}+u_{1}^{i}\partial_{i}\theta_{R}))+\rho_{1}u_{R}^{i}\partial_{i}\theta_{0}+\rho_{R}u_{1}^{i}\partial_{i}\theta_{0}\},
\end{aligned}
\right.
\end{equation}
%where the index $i$ is summed from one to three while $j,k\in \{1,2,3\}$ are fixed.
Multiplying $\eqref{4.25}_{1}$, $\eqref{4.25}_{2}$ and $\eqref{4.25}_{3}$ by $\chi_{0}^2\nabla (\theta_{0}\rho_{R})$, $\chi_{0}^2\nabla (\rho_{0}u_{R}^{k})$ and $\chi_{0}^2\nabla (\rho_{0}\theta_{R})$ respectively, then summing up the resulting equations, we obtain that
\begin{align}\label{4.26}
&\mu\|\chi_{0}\sqrt{\rho_{0}}\nabla^2u_{R}\|_{L^2}^2+\zeta\|\chi_{0}\sqrt{\rho_{0}}\nabla\operatorname{div}u_{R}\|_{L^2}^2+\kappa\|\chi_{0}\sqrt{\frac{\rho_{0}}{\theta_{0}}}\nabla^2\theta_{R}\|_{L^2}^2\nonumber\\
&\leq \frac{\mu}{16}\|\chi_{0}\sqrt{\rho_{0}}\nabla^2u_{R}\|_{L^2}^2+\frac{\zeta}{16}\|\chi_{0}\sqrt{\rho_{0}}\nabla \operatorname{div}u_{R}\|_{L^2}^2 +\frac{\kappa}{16}\|\chi_{0}\sqrt{\frac{\rho_{0}}{\theta_{0}}}\nabla ^2\theta_{R}\|_{L^2}^2+C\|(u_{R},\theta_{R})\|_{H^1}^2\nonumber\\
&\quad -\frac{1}{\v}\int_{\Omega}\big[\chi_{0}^2\nabla \operatorname{div}(\rho_{0}u_{R})\cdot \nabla (\rho_{R}\theta_{0}+\rho_{0}\theta_{R})+\chi_{0}^2\nabla (\rho_{0}u_{R}):\nabla^2(\rho_{R}\theta_{0}+\rho_{0}\theta_{R})\big]\,{\rm d}x\nonumber\\
&\quad -\int_{\Omega}\chi_{0}^2\nabla [\rho_{0}(u_{1}\cdot \nabla u_{R}+u_{R}\cdot u_{1})]:\nabla (\rho_{0}u_{R})\,{\rm d}x
-\int_{\Omega}\chi_{0}^2\nabla^2(\rho_{R}\theta_{1}+\rho_{1}\theta_{R}):\nabla (\rho_{0}u_{R})\,{\rm d}x\nonumber\\
&\quad -\int_{\Omega}\chi_{0}^2\nabla \tilde{F}(\rho_{R},u_{R},\theta_{R}):\nabla (\rho_{0}u_{R})\,{\rm d}x+\int_{\Omega}\chi_{0}^2\nabla \mathfrak{R}_{2}:\nabla (\rho_{0}u_{R})\,{\rm d}x\nonumber\\
&\quad -P_{1}\int_{\Omega}\chi_{0}^2\nabla (\theta_{0}^{-1}\operatorname{div}u_{R})\cdot \nabla(\rho_{0}\theta_{R})\,{\rm d}x-\int_{\Omega}\chi_{0}^2\nabla [\theta_{0}^{-1}(\theta_{0}\rho_{R}+\rho_{0}\theta_{R})\operatorname{div}u_{1}]\cdot \nabla (\rho_{0}\theta_{R})\,{\rm d}x\nonumber\\
&\quad -\int_{\Omega}\chi_{0}^2\nabla \{\theta_{0}^{-1}[\rho_{0}(u_{R}\cdot \nabla \theta_{1}+u_{1}\cdot \nabla \theta_{R})+\rho_{1}u_{R}\cdot \nabla \theta_{0}+\rho_{R}u_{1}\cdot \nabla \theta_{0}]\}\cdot \nabla (\rho_{0}\theta_{R})\,{\rm d}x\nonumber\\
&\quad -\int_{\Omega}\chi_{0}^2\nabla (\theta_{0}^{-1}\tilde{G}^{\varepsilon}(\rho_{R},u_{R},\theta_{R}))\cdot \nabla (\rho_{0}\theta_{R})\,{\rm d}x-\int_{\Omega}\chi_{0}^2\nabla \mathfrak{R}_{3}\cdot \nabla (\rho_{0}\theta_{R})\,{\rm d}x\nonumber\\
&\quad -\int_{\Omega}\chi_{0}^2\nabla \operatorname{div}[\rho_{R}(u_{1}+\v(u_{2}+\tilde{u}_{R}))]\cdot \nabla (\theta_{0}\rho_{R})\,{\rm d}x+\int_{\Omega}\chi_{0}^2\nabla \mathfrak{R}_{1}\cdot \nabla (\theta_{0}\rho_{R})\,{\rm d}x\nonumber\\
&=:\sum\limits_{m=1}^{16}I_{m}.
\end{align}
Now we control $I_{m}\,\,(m=5,\cdots 16)$ term by term. For $I_{5}$, integrating by parts to see that
\begin{align}\label{4.26-1}
I_{5}&=-\frac{1}{\varepsilon}\sum\limits_{k,j=1}^3\Big[\int_{\Omega}\chi_{0}^2\partial_{jk}^2(\rho_{0}u_{R}^{k})\partial_{j}(\rho_{R}\theta_{0}+\rho_{0}\theta_{R})\,{\rm d}x-\int_{\Omega}\chi_{0}^2\partial_{jk}^2(\rho_{0}u_{R}^{k})\partial_{j}(\rho_{R}\theta_{0}+\theta_{0}\rho_{R})\,{\rm d}x\nonumber\\
&\qquad\qquad\quad\,\, -\int_{\Omega}2\chi_{0}\partial_{k}\chi_{0}\partial_{j}(\rho_{0}u_{R}^{k})\partial_{j}(\rho_{0}\theta_{R}+\rho_{R}\theta_{0})\,{\rm d}x\Big]\nonumber\\
&=\frac{1}{\varepsilon}\sum\limits_{k,j=1}^3\int_{\Omega}2\chi_{0}\partial_{k}\chi_{0}\partial_{j}(\rho_{0}u_{R}^{k})\partial_{j}(\rho_{0}\theta_{R}+\rho_{R}\theta_{0})\,{\rm d}x
\nonumber\\
&\leq \tau\|\frac{\nabla(\rho_{0}\theta_{R}+\rho_{R}\theta_{0})}{\varepsilon}\|_{L^2}^2+C_{\tau}\|u_{R}\|_{H^1}^2
\end{align}
for any $\tau>0$.

For $I_{m}$ ($m=6,7,8,9$), using H\"{o}lder inequality and integrating by parts, one has
\begin{align}\label{4.26-3}
\sum\limits_{m=6}^{9}|I_{m}|&\leq \frac{\mu}{16}\|\chi_{0}\sqrt{\rho_{0}}\nabla^2u_{R}^{k}\|_{L^2}^2 +C(\mathcal{I}_{0})\|(\rho_{R},u_{R},\theta_{R})\|_{H^1}^2+C\big(\|\tilde{F}^{\varepsilon}\|_{L^2}^2+\|\mathfrak{R}_{2}\|_{L^2}^2\big)\nonumber\\
&\leq \frac{\mu}{16}\|\chi_{0}\sqrt{\rho_{0}}\nabla^2 u_{R}\|_{L^2}^2+C(\mathcal{I}_{1})\|(\rho_{R},u_{R},\theta_{R})\|_{H^1}^2+C(\|\mathfrak{R}_{2}\|_{L^2}^2)
\end{align}
where we have used \eqref{F8-3} in the last inequality.

Similarly, for $I_{m}$ ($m=10,\cdots,14$), we also integrate by parts to obtain that
\begin{align}\label{4.26-4}
\sum\limits_{m=10}^{14}|I_{m}|&\leq \frac{\kappa}{16}\|\chi_{0}\sqrt{\frac{\rho_{0}}{\theta_{0}}}\nabla^2\theta_{R}\|_{L^2}^2+C(\mathcal{I}_{0})\|(\rho_{R},u_{R},\theta_{R})\|_{H^1}^2+C(\|\tilde{G}\|_{L^2}^2+\|\mathfrak{R}_{3}\|_{L^2}^2)\nonumber\\
&\leq \frac{\kappa}{16}\|\chi_{0}\sqrt{\frac{\rho_{0}}{\theta_{0}}}\nabla^2\theta_{R}\|_{L^2}^2+C(\mathcal{I}_{1})\|(\rho_{R},u_{R},\theta_{R})\|_{H^1}^2+C\|\mathfrak{R}_{3}\|_{L^2}^2
\end{align}
where we have used \eqref{F6-1}.

For $I_{15}$, by using the H\"{o}lder inequality and Sobolev embedding, one has
\begin{align}\label{4.26-5}
	|I_{15}|\leq C(\mathcal{I}_{1})\|\rho_{R}\|_{H^2}^2. %C\|u_{1}\|_{H^2}\|\rho_{R}\|_{H^2}^2
\end{align}

Finally, for $I_{16}$, by using H\"{o}lder inequality, we obtain
\begin{align}\label{4.26-8}
	|I_{16}|\leq C\|\rho_{R}\|_{H^1}\|\mathfrak{R}_{1}\|_{H^1}.
\end{align}
Substituting \eqref{4.26-1}--\eqref{4.26-8} into \eqref{4.26}, we obtain \eqref{5.1-1}. Therefore the proof of Lemma \ref{lem4.2} is complete. $\hfill\square$

%\begin{align*}
	%&-\int_{\Omega}\partial_{ji}(\rho_{R}u_{1}^{i})\chi_{0}^2\partial_{j}(\theta_{0}\rho_{R})\,{\rm d}x=-\int_{\Omega}u_{1}^{i}\partial_{ji}\rho_{R}\chi_{0}^2\theta_{0}\partial_{j}\rho_{R}\,{\rm d}x-\int_{\Omega}u_{1}^{i}\partial_{ji}\rho_{R}\chi_{0}^2\partial_{j}\theta_{0}\rho_{R}\,{\rm d}x
	%-\int_{\Omega}\chi_{0}^2(\partial_{ij}(\rho_{R}u_{1}^{i})-u_{1}^{i}\partial_{ji}\rho_{R})\partial_{j}(\theta_{0}\rho_{R})\,{\rm d}x\\
	%&\leq \frac{1}{2}\int_{\Omega}|\partial_{j}\rho_{R}|^2\partial_{i}(\chi_{0}^2\theta_{0}u_{1}^{i})\,{\rm d}x+\int_{\Omega}\partial_{j}\rho_{R}\partial_{i}(\rho_{R}\chi_{0}^2\partial_{j}\theta_{0}u_{1}^{i})\,{\rm d}x+C\|u_{1}\|_{H^3}\|\rho_{R}\|_{H}
	%-\int_{\Omega}|\partial_{j}\rho_{R}|^2\chi_{0}^2\theta_{0}\partial_{i}u_{1}^{i}\,{\rm d}x-\int_{\Omega}u_{1}^{i}\partial_{ji}\rho_{R}\chi_{0}^2\theta_{0}\partial_{j}\rho_{R}\,{\rm d}x-\int_{\Omega}\
	%\int_{\Omega}\chi_{0}^2\partial_{j}(\rho_{R}u_{1}^{i})\partial_{ji}\theta_{0}\rho_{R}\,{\rm d}x+\int_{\Omega}\chi_{0}^2\partial_{j}(\rho_{0}u_{})
%\end{align*}

\begin{lemma}[Interior estimate of $\varepsilon\|\nabla^3 u_{R}\|_{L^2}$]\label{lem4.3}
For any $\tau>0$, there exists a positive constant $C_{\tau}$, which depends on $\tau$ but is independent of $\v$, such that
%There exists a positive constant $C$ independent of $\varepsilon$, such that
\begin{align}\label{5.1-2}
		&\varepsilon^2\mu\|\chi_{0}\sqrt{\rho_{0}}\nabla^3u_{R}\|_{L^2}^2+\varepsilon^2\zeta\|\chi_{0}\sqrt{\rho_{0}}\nabla^2\operatorname{div}u_{R}\|_{L^2}^2+\varepsilon^2\kappa\|\chi_{0}\sqrt{\frac{\rho_{0}}{\theta_{0}}}\nabla^3\theta_{R}\|_{L^2}^2\nonumber\\
		&\leq (C(\mathcal{I}_{1})\v^2+\tau)\|\rho_{R}\|_{H^2}^2+C_{\tau}\v^2\|u_{R}\|_{H^2}^2+C\|\theta_{R}\|_{H^2}^2\nonumber\\
		&\quad +C\v^2\big[\|\rho_{R}\|_{H^2}\|\mathfrak{R}_{1}\|_{H^2}+\|(\mathfrak{R}_{2},\mathfrak{R}_{3})\|_{H^1}^2\big],
	\end{align}
where $C>0$ is independent of $\tau$ and $\v$.
\end{lemma}

\textbf{Proof. } We apply $\partial^2$ to \eqref{3.10} to get that
\begin{equation}\label{4.27}
	\left\{\begin{aligned}
	&\partial^2\operatorname{div}[\rho_{R}(u_{1}+\v(u_{2}+\tilde{u}_{R}))]=-\frac{1}{\v}\partial^2\operatorname{div}(\rho_{0}u_{R})+\partial^2\mathfrak{R}_{1},\\
	&-\mu\Delta \partial^2 u_{R}^{k}-\zeta\partial^2(\partial_{k}\operatorname{div}u_{R})=-\frac{\partial^2\partial_{k}(\rho_{0}\theta_{R}+\theta_{0}\rho_{R})}{\v}-\partial^2\big[\rho_{0}(u_{1}\cdot \partial u_{R}^{k}+u_{R}\cdot \partial u_{1}^{k})\big]\\
	&\qquad\qquad\qquad\qquad\qquad\qquad\,\,\,\,\,\,\,\,\, -\partial^2\partial_{k}(\rho_{R}\theta_{1}+\rho_{1}\theta_{R})-\partial^2\tilde{F}^{\v,k}+\partial^2\mathfrak{R}_{2}^{k}\,\,(k=1,2,3),\\
	&-\frac{\kappa}{\theta_{0}}\Delta \partial^2\theta_{R}=-\frac{\partial^3\operatorname{div}(\rho_{0}u_{R})}{\v}-P_{1}\partial^2(\theta_{0}^{-1}\operatorname{div}u_{R})-\partial^2[\theta_{0}^{-1}(\rho_{0}\theta_{R}+\rho_{R}\theta_{0})\operatorname{div}u_{1}]\\
	&\qquad\qquad\qquad\,\,\,\,
		-\partial^2\big\{\theta_{0}^{-1}[\rho_{0}(u_{R}\cdot \partial \theta_{1}+u_{1}\cdot \theta_{R})+\rho_{1}u_{R}\cdot \partial \theta_{0}+\rho_{R}u_{1}\cdot \partial\theta_{0}]\big\}\\
		&\qquad\qquad\qquad\,\,\,\, -\kappa[\partial^2(\theta_{0}^{-1}\Delta \theta_{R})-\theta_{0}^{-1}\Delta \partial^2\theta_{R}]-\partial^2(\theta_{0}^{-1}\tilde{G}^{\varepsilon})-\partial^2\mathfrak{R}_{3}.
	\end{aligned}
	\right.
\end{equation}
We multiply $\eqref{4.27}_{1}$, $\eqref{4.27}_{2}$ and $\eqref{4.27}_{3}$ by $\varepsilon^2\chi_{0}^2\partial^2(\theta_{0}\rho_{R})$, $\varepsilon^2\chi_{0}^2\partial^2(\rho_{0}u_{R}^{k})$ respectively and $\varepsilon^2\chi_{0}^2\partial^2(\rho_{0}\theta_{R})$, and add the resultant equations together to get that
\begin{align}\label{4.28}
&\varepsilon^2\mu\int_{\Omega}\chi_{0}^2\rho_{0}|\nabla^3 u_{R}|^2\,{\rm d}x+\varepsilon^2\zeta\int_{\Omega}\chi_{0}^2\rho_{0}|\nabla^2\operatorname{div}u_{R}|^2\,{\rm d}x
+\varepsilon^2\kappa\int_{\Omega}\frac{\rho_{0}}{\theta_{0}}\chi_{0}^2|\nabla^3\theta_{R}|^2\,{\rm d}x\nonumber\\
&\leq \frac{\v^2\mu}{16}\int_{\Omega}\chi_{0}^2\rho_{0}|\nabla^3 u_{R}|^2\,{\rm d}x+\frac{\v^2\zeta}{16}\int_{\Omega}\chi_{0}^2\rho_{0}|\nabla^2 \operatorname{div}u_{R}|^2\,{\rm d}x+\frac{\v^2\kappa}{16}\int_{\Omega}\frac{\rho_{0}}{\theta_{0}}\chi_{0}^2|\nabla^3\theta_{R}|^2\,{\rm d}x+C\v^2\|(u_{R},\theta_{R})\|_{H^2}^2\nonumber\\
&\quad -\v\int_{\Omega}\chi_{0}^2\big[\nabla^2\operatorname{div}(\rho_{0}u_{R}):\nabla^2(\rho_{R}\theta_{0}+\rho_{0}\theta_{R})+\sum\limits_{k=1}^3\nabla^2(\rho_{0}u_{R}^{k}):\nabla^2\partial_{k}(\rho_{0}\theta_{R}+\rho_{R}\theta_{0})\big]\,{\rm d}x\nonumber\\
&\quad -\v^2\int_{\Omega}\chi_{0}^2\sum\limits_{k=1}^3\nabla^2[\rho_{0}(u_{1}\cdot \nabla u_{R}^{k}+u_{R}\cdot \nabla u_{1}^{k})]:\nabla^2(\rho_{0}u_{R}^{k})\,{\rm d}x\nonumber\\
&\quad -\v^2\int_{\Omega}\chi_{0}^2\sum\limits_{k=1}^3\nabla^2\partial_{k}(\rho_{1}\theta_{R}+\rho_{R}\theta_{1}):\nabla^2(\rho_{0}u_{R}^{k})\,{\rm d}x-\v^2\int_{\Omega}\chi_{0}^2\sum\limits_{k=1}^{3}\nabla^2\tilde{F}^{\v,k}:\nabla^2(\rho_{0}u_{R}^{k})\,{\rm d}x\nonumber\\
&\quad -\v^2\int_{\Omega}\chi_{0}^2\sum\limits_{k=1}^{3}\nabla^2\mathfrak{R}_{2}^{k}:\nabla^2(\rho_{0}u_{R}^{k})\,{\rm d}x-\v^2\int_{\Omega}\chi_{0}^2\nabla^2\mathfrak{R}_{3}:\nabla^2(\rho_{0}\theta_{R})\,{\rm d}x\nonumber\\
&\quad -\v^2\int_{\Omega}P_{1}\chi_{0}^2\nabla^2(\theta_{0}^{-1}\operatorname{div}u_{R}):\nabla^2(\rho_{0}\theta_{R})\,{\rm d}x-\v^2\int_{\Omega}\chi_{0}^2\nabla^2[\theta_{0}^{-1}(\theta_{0}\rho_{R}+\rho_{0}\theta_{R})\operatorname{div}u_{1}]:\nabla^2(\rho_{0}\theta_{R})\,{\rm d}x\nonumber\\
&\quad -\v^2\int_{\Omega}\nabla^2\big\{
\theta_{0}^{-1}[\rho_{0}(u_{R}\cdot \nabla \theta_{1}+u_{1}\cdot \nabla \theta_{R})+\rho_{1}u_{R}\cdot \theta_{0}+\rho_{R}u_{1}\cdot \nabla \theta_{0}]\big\}:\nabla^2(\rho_{0}\theta_{R})\,{\rm d}x\nonumber\\
&\quad-\v^2\int_{\Omega}\chi_{0}^2\nabla^2(\theta_{0}^{-1}\tilde{G}^{\v}):\nabla^2(\rho_{0}\theta_{R})\,{\rm d}x-\v^2\int_{\Omega}\chi_{0}^2\nabla^2\operatorname{div}[\rho_{R}(u_{1}+\v(u_{2}+\tilde{u}_{R}))]:\nabla^2(\theta_{0}\rho_{R})\,{\rm d}x\nonumber\\
&\quad +\v^2\int_{\Omega}\chi_{0}^2\nabla^2\mathfrak{R}_{1}:\nabla^2(\theta_{0}\rho_{R})\,{\rm d}x=:\sum\limits_{m=1}^{16}J_{m}.
\end{align}

Almost all $J_{m}$ can be similarly estimated as in $I_{m}$ in \eqref{4.26}. For simplicity of presentation, we only focus on terms involving $\nabla^3\rho_{R}$, that is, $J_{5}$, $J_{7}$, $J_{8}$ and $J_{15}$ . For $J_{5}$, we use the H\"{o}lder inequality and integrate by parts to see that
\begin{align}\label{4.27-1}
J_{5}&=-\varepsilon\sum\limits_{j,l,k=1}^{3}\int_{\Omega}\chi_{0}^2[\partial_{jlk}^3(\rho_{0}u_{R}^{k})\partial_{jl}^2(\rho_{0}\theta_{R}+\theta_{0}\rho_{R})-\partial_{jlk}^3(\rho_{0}\theta_{R}+\theta_{0}\rho_{R})\partial_{jl}^2(\rho_{0}u_{R}^{k})]\,{\rm d}x\nonumber\\
&=-\varepsilon\sum\limits_{j,l,k=1}^{3}\int_{\Omega}2\chi_{0}\partial_{k}\chi_{0}\partial_{jl}^2(\rho_{0}\theta_{R}+\theta_{0}\rho_{R})\partial_{jl}^2(\rho_{0}u_{R}^{k})\,{\rm d}x\nonumber\\
&\leq \tau \|\rho_{R}\|_{H^2}^2+C_{\tau}\v^2\|u_{R}\|_{H^2}^2+C\|\theta_{R}\|_{H^2}^2\quad \text{for any $\tau>0$}.
\end{align}
Similarly for $J_{7}$ and $J_{8}$, by integrating by parts, one has
\begin{align}
&|J_{7}|\leq C\varepsilon^2\sum\limits_{j,l,k=1}^3\int_{\Omega}|\partial_{jl}^2(\rho_{R}\theta_{1}+\rho_{1}\theta_{R})\partial_{k}(\chi_{0}^2\partial_{jl}^2(\rho_{0}u_{R}^{k}))|\,{\rm d}x\nonumber\\
&\quad\,\,\,\leq C(\mathcal{I}_{0})\varepsilon^2\|(\rho_{R},u_{R},\theta_{R})\|_{H^2}^2+\frac{\varepsilon^2\mu}{16}\int_{\Omega}\chi_{0}^2\rho_{0}|\nabla^3u_{R}|^2\,{\rm d}x,\label{4.27-2}
\end{align}
and
\begin{align}
&|J_{8}|\leq C(\mathcal{I}_{1})\v^2\|(\rho_{R},u_{R},\theta_{R})\|_{H^2}^2+\frac{\varepsilon^2\mu}{16}\int_{\Omega}\chi_{0}^2\rho_{0}|\nabla^3u_{R}|^2\,{\rm d}x.\label{4.27-5}
%|J_{12}|\leq C\varepsilon^2\sum\limits_{j,l,k=1}^3\int_{\Omega}|\partial_{jl}^2()\partial_{k}(\chi_{0}^2\partial_{jl}^2(\rho_{0}u_{R}^{k}))|\,{\rm d}x
%\Big\vert\varepsilon^2\int_{\Omega}\partial_{jil}^3\rho_{R}u_{1}^{i}\chi_{0}^2\partial_{jl}^2(\theta_{0}\rho_{R})\,{\rm d}x+\int_{\Omega}[\partial_{jil}^3(\rho_{R}u_{1}^{i})-\partial_{jil}^3\rho_{R}u_{1}^{i}]\chi_{0}^2\partial_{jl}^2(\rho_{0}\rho_{R})\,{\rm d}x\Big\vert\nonumber\\
%&\leq \Big\vert\frac{\varepsilon^2}{2}\int_{\Omega}|\partial_{jl}^2\rho_{R}|^2\partial_{i}(u_{1}^{i}\chi_{0}^2\theta_{0})\,{\rm d}x\Big\vert+\Big\vert\varepsilon^2\int_{\Omega}\partial_{jl}^2\rho_{R}\partial_{i}[u_{1}^{i}\chi_{0}^2(\partial_{jl}^2(\theta_{0}\rho_{R})-\theta_{0}\partial_{jl}^2\rho_{R})]\,{\rm d}x\Big\vert\nonumber\\
%&\quad +C\varepsilon^2\|u_{1}\|_{H^3}\|\rho_{R}\|_{H^{2}}^2\leq C\varepsilon^2\|u_{1}\|_{H^3}\|\rho_{R}\|_{H^{2}}^2,\label{4.27-3}.
\end{align}
For $J_{15}$, it follows from integrating by parts that
\begin{align}\label{4.27-4}
	|J_{15}|&\leq \frac{1}{2}\sum\limits_{j,l=1}^3\varepsilon^2\int_{\Omega}\theta_{0}\chi_{0}^2[u_{1}+\v(u_{2}+\tilde{u}_{R})]\cdot \nabla (\partial_{jl}^2(\rho_{R}))^2\,{\rm d}x+C(\mathcal{I}_{1})\v^2\|\rho_{R}\|_{H^2}^2\nonumber\\
	&\leq C(\mathcal{I}_{1})\v^2\|\rho_{R}\|_{H^2}^2.
\end{align}
Using \eqref{4.27-1}--\eqref{4.27-4} and similar arguments as in Lemma \ref{lem4.2}, we obtain \eqref{5.1-2}. Therefore the proof of Lemma \ref{lem4.3} is complete. $\hfill\square$

\subsubsection{Estimates near the boundary}
Motivated by \cite{Choe-Jin,Li-Liao-2019}, to bound $\nabla^2\operatorname{div}u_{R}$ in the vicinity of the boundary, we straighten the boundary by a local coordinate. Recall the covering \eqref{covering}, since $\partial\Omega$ is given locally by $x_{3}=w_{k}(x_{1},x_{2})$ (hereafter we omit the subscript $k$ for notational convenience), it is convenient to use the coordinate
$$
z_{1}=x_{1},\quad z_{2}=x_{2},\quad z_{3}=x_{3}-w(x_{1},x_{2}).
$$
then it holds that $z_{3}=0$ on $\partial\Omega\cap Q$ and
$$
\frac{\partial}{\partial x_{i}}=\frac{\partial}{\partial z_{i}}-\frac{\partial w}{\partial z_{i}}\frac{\partial}{\partial z_{3}}\,\,(i=1,2),\quad \frac{\partial}{\partial x_{3}}=\frac{\partial}{\partial z_{3}},\qquad \operatorname{det}(\frac{\partial x}{\partial z})=1.
$$
Let $\hat{Q}$ and $\hat{\Omega}$ be the images of $Q$ and $\Omega$ by the above coordinate transformation. We denote
\begin{align*}
	&\hat{f}(z)=:f(z_{1},z_{2},z_{3}+w(z_{1},z_{2}))\quad \text{for any scalar function }f\text{ on }\Omega,\\
	&\hat{D}_{m}=\frac{\partial}{\partial z_{m}},\quad \Delta_{z}=\hat{D}_{11}+\hat{D}_{22}+\hat{D}_{33},\,\,\nabla_{z}=(\hat{D}_{1},\hat{D}_{2},\hat{D}_{3}),\quad \operatorname{div}_{z}\hat{v}=\hat{D}_{1}\hat{v}_{1}+\hat{D}_{2}\hat{v}_{2}+\hat{D}_{3}\hat{v}_{3}.
\end{align*}

%For any vector field $g=(g_{1},g_{2},g_{3})$ in $x-$coordinate, we denote
%$$
%\hat{g}^{i}:=\hat{g}_{i}-\hat{g}_{3}\partial_{z_{i}}w\,\,(i=1,2)\quad \hat{g}^{3}=-\hat{g}_{1}\hat{D}_{1}w-\hat{g}_{2}\hat{D}_{2}w+\hat{g}_{3}(1+(\hat{D}_{1}w)^2+(\partial_{z_2}w)^2),
%$$
%and then
%\begin{align*}
	%&\hat{g}_{1}=\hat{g}^{1}(1-(\hat{D}_{1}w)^2)-\hat{g}^2\hat{D}_{1}w\hat{D}_{2}w+\hat{g}^{3}\hat{D}_{1}w,\\
	%&\hat{g}_{2}=-\hat{g}^{1}\hat{D}_{1}w\hat{D}_{2}w+\hat{g}^2(1-(\hat{D}_{2}w)^2)+\hat{g}^{3}\hat{D}_{2}w,\\
	%&\hat{g}_{3}=-\hat{g}^{1}\hat{D}_{1}w-\hat{g}^2\hat{D}_{2}w+\hat{g}^{3}.
%\end{align*}
%%For any $x_{0}\in \partial\Omega$, without loss of generality, we assume $x_{0}$ is the origin and $x_{3}=w(x_{1},x_{2})$ on $\partial_{\Omega}\cap Q$, where $w\in C^{4}(\R^2)$ and $Q$ is a small cube centered at the origin, and $Q$ is small enough so that
%Further, we define the local coordinate $z=(z_{1},z_{2},z_{3})$ by
Then in $\hat{\Omega}\cap \hat{Q}$, \eqref{3.10} can be reformed as
\begin{equation}\label{5.2}
\left\{
\begin{aligned}
	%&\operatorname{div}_{z}(\hat{\rho}_{R}\hat{u}_{1})=-\frac{1}{\varepsilon}\operatorname{div}_{z}(\hat{\rho}_{0}\hat{u}_{R})-\operatorname{div}_{z}(\hat{\rho}_{1}\hat{\tilde{u}}_{R})-\varepsilon\operatorname{div}_{z}(\hat{\rho}_{R}\hat{u}_{2}+\hat{\rho}_{2}\hat{\tilde{u}}_{R}+\hat{\rho}_{R}\hat{\tilde{u}}_{R})\\
	&\operatorname{div}_{z}[\hat{\rho}_{R}(\hat{u}_{1}+\v(\hat{u}_{2}+\hat{\tilde{u}}_{R}))]=-\frac{1}{\v}\operatorname{div}_{z}(\hat{\rho}_{0}\hat{u}_{R})+\frac{1}{\v}\sum\limits_{m=1}^2\hat{D}_{m}w\hat{D}_{3}(\hat{\rho}_{0}\hat{u}_{R}^{m}) +\hat{W}_{1}+\hat{\mathfrak{R}}_{1},\\
	%&\qquad\qquad \qquad +\frac{1}{\varepsilon}\sum\limits_{m=1}^2\hat{D}_{m}w\hat{D}_{3}(\hat{\rho}_{0}\hat{u}_{R}^{m})-\varepsilon^2\operatorname{div}_{z}(\hat{\rho}_{3}\hat{\tilde{u}}_{R})+\hat{W}_{1}+\hat{r}_{1},\\
	&\mu\Delta_{z} \hat{u}_{R}+\zeta\nabla_{z} \operatorname{div}_{z}\hat{u}_{R}=\frac{1}{\varepsilon}\nabla_{z}(\hat{\rho}_{0}\hat{\theta}_{R}+\hat{\rho}_{R}\hat{\theta}_{0})+\hat{\rho}_{0}(\hat{u}_{1}\cdot \nabla_{z} \hat{u}_{R}+\hat{u}_{R}\cdot \nabla_{z} \hat{u}_{1})\\
	&\qquad\qquad\qquad\qquad\qquad\,\,\,
	-\frac{1}{\varepsilon}\big(\hat{D}_{1}w\hat{D}_{3}(\hat{\rho}\hat{\theta}_{R}+\hat{\theta}_{0}\hat{\rho}_{R}),\hat{D}_{2}w\hat{D}_{3}(\hat{\rho}\hat{\theta}_{R}+\hat{\theta}_{0}\hat{\rho}_{R}),0\big)\\
	 &\qquad\qquad\qquad\qquad\qquad\,\,\, +\nabla_{z}(\hat{\rho}_{1}\hat{\theta}_{R}+\hat{\rho}_{R}\hat{\theta}_{1})+\hat{W}_{2}+\hat{\tilde{F}}^{\varepsilon}(\hat{\rho}_{R},\hat{u}_{R},\hat{\theta}_{R})+\hat{\mathfrak{R}}_{2},\\
	&\frac{\kappa}{\hat{\theta}_{0}}\Delta_{z} \hat{\theta}_{R}=\frac{1}{\varepsilon}\operatorname{div}_{z}(\hat{\rho}_{0}\hat{u}_{R})+\frac{P_{1}}{\hat{\theta}_{0}}\operatorname{div}_{z}\hat{u}_{R}+\frac{1}{\hat{\theta}_{0}}(\hat{\rho}_{0}\hat{\theta}_{R}+\hat{\rho}_{R}\hat{\theta}_{0})\operatorname{div}_{z}\hat{u}_{1}+\frac{1}{\hat{\theta}_{0}}\hat{\tilde{G}}^{\varepsilon}(\hat{\rho}_{R},\hat{u}_{R},\hat{\theta}_{R})\\
	&\qquad\qquad\,\,\,  +\frac{1}{\hat{\theta}_{0}}\big\{\hat{\rho}_{0}(\hat{u}_{R}\cdot \nabla_{z} \hat{\theta}_{1}+\hat{u}_{1}\cdot \nabla_{z} \hat{\theta}_{R})+\hat{\rho}_{1}\hat{u}_{R}\cdot \nabla_{z}\hat{\theta}_{0}+\hat{\rho}_{R}\hat{u}_{1}\cdot \nabla_{z} \hat{\theta}_{0}\big\},\\
	&\qquad\qquad\,\,\, -\frac{1}{\varepsilon}\sum\limits_{m=1}^2\hat{D}_{m}w\hat{D}_{3}(\hat{\rho}_{0}\hat{u}_{R}^{m})+\hat{W}_{3}+\hat{\mathfrak{R}}_{3},\\
	%+\frac{1}{\hat{\theta}_{0}}\hat{r}_{3}++\frac{1}{\hat{\theta}_{0}}\hat{\Psi}(\nabla (u_{1}+\varepsilon(u_{2}+\tilde{u}_{R})))\\
	&\hat{u}_{R}\vert_{\partial\hat{\Omega}\cap \Omega}=\hat{\theta}_{R}\vert_{\partial\hat{\Omega}\cap \hat{Q}}=0,
\end{aligned}
\right.
\end{equation}
where $\hat{W}_{1}$, $\hat{W}_{2}=(\hat{W}_{2}^{1},\hat{W}_{2}^{2},\hat{W}_{2}^{3})$ and $\hat{W}_{3}$ are defined as
\begin{align}\label{5.3}
%&\hat{W}_{1}=\sum\limits_{m=1}^{2}\hat{D}_{m}w\hat{D}_{3}\big[\hat{\rho}_{R}(\hat{u}_{1}^{m}+\varepsilon\hat{u}_{2}^{m}+\varepsilon\hat{u}_{3}^{m})+\hat{\tilde{u}}_{R}^{m}(\hat{\rho}_{1}+\varepsilon\hat{\rho}_{2}+\varepsilon^2\rho_{3})\big],\nonumber\\
&\hat{W}_{1}=\sum\limits_{m=1}^2\hat{D}_{m}w\hat{D}_{3}\big[\hat{\rho}_{R}(\hat{u}_{1}^{m}+\varepsilon\hat{u}_{2}^{m}+\varepsilon\hat{u}_{3}^{m})\big],\quad \hat{W}_{2}^3=-\zeta\sum\limits_{m=1}^2\hat{D}_{m}w\hat{D}_{33}^2\hat{u}_{R}^{m},\nonumber\\
&\hat{W}_{2}^{j}=\mu\sum\limits_{m=1}^2\big[\hat{D}_{mm}^2w\hat{D}_{3}\hat{u}_{R}^{j}+2\hat{D}_{m}w\hat{D}_{m3}^2\hat{u}_{R}^{j}-(\hat{D}_{m}w)^2\hat{D}_{33}^2\hat{u}_{R}^{j}\big]\nonumber\\
&\qquad\,\, +\zeta\big[\hat{D}_{j}w\hat{D}_{3}(\operatorname{div}_{z}\hat{u}_{R})-\hat{D}_{j}(\sum\limits_{m=1}^2\hat{D}_{m}w\hat{D}_{3}\hat{u}_{R}^{m})-\sum\limits_{m=1}^2\hat{D}_{m}w\hat{D}_{j}w\hat{D}_{33}^2\hat{u}_{R}^{m}\big],\nonumber\\
&\qquad\,\, -\hat{D}_{j}w[\hat{D}_{3}(\hat{\rho}_{R}\hat{\theta}_{1})+\hat{D}_{3}(\hat{\rho}_{1}\hat{\theta}_{R})] -\hat{\rho}_{0}\big[\sum\limits_{m=1}^2\hat{u}_{1}^{m}\hat{D}_{m}w\hat{D}_{3}\hat{u}_{R}^{j}+\sum\limits_{m=1}^2\hat{u}_{R}^{m}\hat{D}_{m}w\hat{D}_{3}\hat{u}_{1}^{j}\big]\qquad j=1,2,\nonumber\\
&\hat{W}_{3}=\frac{\kappa}{\hat{\theta}_{0}}\sum\limits_{m=1}^2\big[\hat{D}_{mm}^2w\hat{D}_{3}\hat{\theta}_{R}+2\hat{D}_{m}w\hat{D}_{m3}^2\hat{\theta}_{R}-(\hat{D}_{m}w)^2\hat{D}_{33}^2\hat{\theta}_{R}\big]\nonumber\\
&\qquad\,\, -\frac{P_{1}}{\hat{\theta}_{0}}\sum\limits_{m=1}^2\hat{D}_{m}w\hat{D}_{3}\hat{u}_{R}^{m}-\frac{1}{\hat{\theta}_{0}}(\hat{\rho}_{0}\hat{\theta}_{R}+\hat{\rho}_{R}\hat{\theta}_{0})\sum\limits_{m=1}^2\hat{D}_{m}w\hat{D}_{3}\hat{u}_{1}^{m},\nonumber\\
&\qquad\,\, -\frac{1}{\hat{\theta}_{0}}\sum\limits_{m=1}^2[\hat{\rho}_{0}(\hat{u}_{R}^{m}\hat{D}_{m}w\hat{D}_{3}\hat{\theta}_{1}+\hat{u}_{R}^{m}\hat{D}_{m}w\hat{D}_{3}\hat{\theta}_{0})+\hat{\rho}_{1}\hat{u}_{R}^{m}\hat{D}_{m}w\hat{D}_{3}\hat{\theta}_{0}+\hat{\rho}_{R}\hat{u}_{1}^{m}\hat{D}_{m}w\hat{D}_{3}\hat{\theta}_{0}].
\end{align}
Here in \eqref{5.2}, we slightly abuse the notations $\hat{\tilde{F}}, \hat{\tilde{G}}$ and $\hat{\mathfrak{R}}_{1}, \hat{\mathfrak{R}}_{2}, \hat{\mathfrak{R}}_{3}$ to represent the corresponding forms of $\tilde{F},\tilde{G}$ and $\mathfrak{R}_{1}, \mathfrak{R}_{2}, \mathfrak{R}_{3}$ in $z-$coordinate respectively. By \eqref{5.1}, it is easy to check that
\begin{align}
&\|\hat{W}_{2}\|_{L^2}+\|\hat{W}_{3}\|_{L^2}
\leq C(\mathcal{I}_{1})\sigma\|\rho_{R}\|_{H^1}+C\sigma \|(u_{R},\theta_{R})\|_{H^2},\quad\,\, \|\hat{W}_{1}\|_{H^1}\leq C(\mathcal{I}_{1})\sigma\|\rho_{R}\|_{H^2},\label{5.4-1}\\
&\|\hat{W}_{2}\|_{H^1}+\|\hat{W}_{3}\|_{H^1}
\leq C(\mathcal{I}_{1})\sigma\|\rho_{R}\|_{H^2}+C\sigma \|(u_{R},\theta_{R})\|_{H^3},\quad\|\hat{\mathfrak{R}}_{1}\|_{H^2}\leq C\|\mathfrak{R}_{1}\|_{H^2}, \label{5.4}\\
&\|(\hat{\tilde{F}},\hat{\tilde{G}})\|_{H^1}\leq C\|(\tilde{F},\tilde{G})\|_{H^1}\leq C(\mathcal{I}_{1})\|(\rho_{R},u_{R},\theta_{R})\|_{H^2},\quad\|(\hat{\mathfrak{R}}_{2},\hat{\mathfrak{R}}_{3})\|_{H^1}\leq C\|(\mathfrak{R}_{2},\mathfrak{R}_{3})\|_{H^1},\label{5.5}
%\\&,\quad \quad ,\label{5.5-1}
\end{align}
where we have used the fact that $\|\hat{f}(z)\|_{H_{z}^{k}}\lesssim \|f(x)\|_{H_{x}^{k}}$ for $k=1,2,3$ in the above inequalities. Hereafter, all the norms $\|f\|_{H^k}$ means a local norm $\|f\|_{H^k(Q)}$ for simplicity of notations.

\smallskip

\noindent\textbf{Step 1. Estimates of tangential derivatives.}
Let $\chi_{1}\in C_{0}^{\infty}(\hat{Q})$ with $\frac{1}{2}\hat{Q}\subset \operatorname{supp}\chi_{1}\subset \hat{Q}$ with
%satisfying
$$
\chi_{1}(z)\equiv 1\quad \text{for $z\in \frac{1}{2}\hat{Q}$}\quad \text{ and }\quad \chi_{1}(z)\equiv 0\quad \text{for $z\in \overline{\hat{Q}}^{c}$}.
$$
It is worth mentioning that $\nabla_{z} \chi_{1}\lesssim \frac{1}{\sigma}$. Similar to \eqref{4.25}, we apply $\hat{D}_{j} (j=1,2)$ to \eqref{5.2} to get
\begin{equation}\label{5.2-0}
	\left\{\begin{aligned}
		%&\hat{D}_{ji}^2(\hat{\rho}_{R}\hat{u}_{1}^{i})+\frac{1}{\varepsilon}(\hat{D}_{ji}^2(\hat{\rho}_{0}\hat{u}_{R}^{i}))=-\hat{D}_{ji}^2(\hat{\rho}_{1}\hat{\tilde{u}}_{R}^{i})-\varepsilon^2\hat{D}_{ji}^2(\hat{\rho}_{3}\hat{\tilde{u}}_{R}^{i})\\
		&\hat{D}_{j}\operatorname{div}_{z}[\hat{\rho}_{R}(\hat{u}_{1}+\v(\hat{u}_{2}+\hat{u}_{R}))]=-\frac{\hat{D}_{j}\operatorname{div}(\hat{\rho}_{0}\hat{u}_{R})}{\v}+\frac{1}{\varepsilon}\sum\limits_{m=1}^2\hat{D}_{j}(\hat{D}_{m}w\hat{D}_{3}(\hat{\rho}_{0}\hat{u}_{R}^{m}))+\hat{D}_{j}\hat{W}_{1}+\hat{D}_{j}\mathfrak{R}_{1},\\
		%&\qquad -\varepsilon\hat{D}_{ji}^2(\hat{\rho}_{R}\hat{u}_{2}^{i}+\hat{\rho}_{2}\hat{\tilde{u}}_{R}^{i}+\hat{\rho}_{R}\hat{\tilde{u}}_{R}^{i})+\hat{D}_{j}\hat{r}_{1}+\hat{D}_{j}\hat{W}_{1} +\frac{1}{\varepsilon}\sum\limits_{m=1}^2\hat{D}_{j}(\hat{D}_{m}w\hat{D}_{3}(\hat{\rho}_{0}\hat{u}_{R}^{m}))\\
		&-\mu\Delta_{z} \hat{D}_{j}\hat{u}_{R}^{k}-\zeta \hat{D}_{jk}^2\operatorname{div}_{z}u_{R}=-\frac{\hat{D}_{jk}^2(\hat{\rho}_{0}\hat{\theta}_{R}+\hat{\theta}_{0}\hat{\rho}_{R})}{\v}-\hat{D}_{j}[\hat{\rho}_{0}(\hat{u}_{1}\cdot \nabla_{z}  \hat{u}_{R}^{k}+\hat{u}_{R}\cdot \nabla_{z} u_{1}^{k})]\\
		%&-\mu\hat{D}_{jii}^3\hat{u}_{R}^{k}-\zeta\hat{D}_{jki}^3\hat{u}_{R}^{i}+\frac{\hat{D}_{jk}^2(\hat{\rho}_{0}\hat{\theta}_{R}+\hat{\theta}_{0}\hat{\rho}_{R})}{\varepsilon}=-\hat{D}_{j}[\hat{\rho}_{0}(\hat{u}_{1}^{i}\hat{D}_{i}\hat{u}_{R}^{k}+\hat{u}_{R}^{i}\hat{D}_{i}\hat{u}_{1}^{i})]\\
		&\qquad\qquad\qquad \qquad\quad	+\frac{1}{\varepsilon}\hat{D}_{j}\big(\hat{D}_{1}w\hat{D}_{3}(\hat{\rho}\hat{\theta}_{R}+\hat{\theta}_{0}\hat{\rho}_{R}),\hat{D}_{2}w\hat{D}_{3}(\hat{\rho}\hat{\theta}_{R}+\hat{\theta}_{0}\hat{\rho}_{R}),0\big)^{k}\\
		&\qquad\qquad\qquad \qquad\quad -\hat{D}_{jk}^2(\hat{\rho}_{R}\hat{\theta}_{1}+\hat{\rho}_{1}\hat{\theta}_{R})-\hat{D}_{j}\hat{\tilde{F}}^{\v,k}-\hat{D}_{j}\hat{\mathfrak{R}}_{2}^{k}-\hat{D}_{j}\hat{W}_{2}^{k}\quad k=1,2,3,\\
		%&\qquad -\hat{D}_{jk}^2(\hat{\rho}_{R}\hat{\theta}_{1})-\hat{D}_{jk}^2(\rho_{1}\theta_{R})-\hat{D}_{j}\hat{\tilde{F}}^{\varepsilon,k}(\hat{\rho}_{R},\hat{u}_{R},\hat{\theta}_{R})-\hat{D}_{j}\hat{r}_{2}^{k}-\hat{D}_{j}\hat{W}_{2}^{k}\\
		&-\frac{\kappa}{\hat{\theta}_{0}}\Delta_{z} \hat{D}_{j}\hat{\theta}_{R}=-\frac{\hat{D}_{j}\operatorname{div}_{z}(\hat{\rho}_{0}\hat{u}_{R})}{\v}-\frac{\kappa}{\hat{\theta}_{0}^2}\hat{D}_{j}\hat{\theta}_{0}\Delta_{z}\hat{\theta}_{R}-P_{1}\hat{D}_{j}(\hat{\theta}_{0}^{-1}\operatorname{div}_{z}\hat{u}_{R})\\
		%&-\frac{\kappa}{\hat{\theta}_{0}}\hat{D}_{jii}^3\hat{\theta}_{R})+\frac{\kappa\hat{D}_{j}\hat{\theta}_{0}}{\hat{\theta}_{0}^2}\hat{D}_{ii}^2\hat{\theta}_{R}\frac{\hat{D}_{ji}^2(\hat{\rho}_{0}\hat{u}_{R}^{i})}{\varepsilon}=-P_{1}\hat{D}_{j}(\hat{\theta}_{0}^{-1}\hat{D}_{i}\hat{u}_{R}^{i})-\hat{D}_{j}[\hat{\theta}_{0}^{-1}(\hat{\rho}_{0}\hat{\theta}_{R}+\hat{\rho}_{R}\hat{\theta}_{0})\hat{D}_{i}\hat{u}_{1}^{i}]\\
		&\qquad\qquad\qquad\quad-\hat{D}_{j}[\hat{\theta}_{0}^{-1}(\hat{\rho}_{0}\hat{\theta}_{R}+\hat{\rho}_{R}\hat{\theta}_{0})\operatorname{div}_{z}\hat{u}_{1}]+\frac{1}{\varepsilon}\sum\limits_{m=1}^2\hat{D}_{j}(\hat{D}_{m}w\hat{D}_{3}(\hat{\rho}_{0}\hat{u}_{R}^{m}))\\
		&\qquad\qquad\qquad\quad  -\hat{D}_{j}\big\{\hat{\theta}_{0}^{-1}[\hat{\rho}_{0}(\hat{u}_{R}\cdot \nabla_{z} \hat{\theta}_{1}+\hat{u}_{1}\cdot \nabla_{z}\hat{\theta}_{R})+\hat{\rho}_{1}\hat{u}_{R}\cdot \nabla_{z} \hat{\theta}_{0}+\hat{\rho}_{R}\hat{u}_{1}\cdot \nabla_{z}\hat{\theta}_{0}]\big\}\\
		%&\qquad -\hat{D}_{j}(\hat{\theta}_{0}^{-1}\hat{\tilde{G}}^{\varepsilon}(\hat{\rho}_{R},\hat{u}_{R},\hat{\theta}_{R}))-\hat{D}_{j}\hat{r}_{3}-\hat{D}_{j}(\hat{\theta}_{0}^{-1}\Psi(\nabla (\hat{u}_{1}+\varepsilon(\hat{u}_{2}+\hat{\tilde{u}}_{R}))))\\
		&\qquad\qquad\qquad\quad  -\hat{D}_{j}\hat{W}_{3}-\hat{D}_{j}(\hat{\theta}_{0}^{-1}\hat{\tilde{G}})-\hat{D}_{j}\hat{\mathfrak{R}}_{3}.
		%\\&\hat{D}_{jl}\hat{u}_{R}\vert_{\hat{D}\hat{\Omega}\cap \hat{Q}}=\hat{D}_{jl}\hat{\theta}_{R}\vert_{\hat{D}\hat{\Omega}\cap \hat{Q}}=0,
	\end{aligned}
	\right.
\end{equation}

We have following estimates on the tangential derivatives.
\begin{lemma}[Estimates of tangential derivatives]\label{lem5.0}
For any $\tau>0$, there exist two positive constants $C_{\sigma}$, which depends on $\sigma$, and $C_{\tau,\sigma}$, which depends only on $\tau$ and $\sigma$, both independent of $\varepsilon$, such that for any $j,l=1,2$,
\begin{align}\label{5.2-2}
&\mu\|\chi_{1}\sqrt{\hat{\rho}_{0}}\nabla_{z}\hat{D}_{j}\hat{u}_{R}\|_{L^2}^2+\zeta\|\chi_{1}\sqrt{\hat{\rho}_{0}}\hat{D}_{j}\operatorname{div}_{z}\hat{u}_{R}\|_{L^2}^2+\kappa\|\chi_{1}\sqrt{\frac{\hat{\rho}_{0}}{\hat{\theta}_{0}}}\nabla_{z}\hat{D}_{j}\hat{\theta}_{R}\|_{L^2}^2\nonumber\\
&\leq \tau\|\frac{\nabla(\rho_{0}\theta_{R}+\rho_{R}\theta_{0})}{\varepsilon}\|_{L^2}^2+(C(\mathcal{I}_{1})+C\sigma^2)\|(u_{R},\theta_{R})\|_{H^2}^2+C_{\tau,\sigma}\|u_{R}\|_{H^1}^2+C_{\sigma}\|\theta_{R}\|_{H^1}^2\nonumber\\
&\quad +C(\mathcal{I}_{1})\|\rho_{R}\|_{H^2}^2 +C\big[\|(\mathfrak{R}_{2},\mathfrak{R}_{3})\|_{L^2}^2+\|\rho_{R}\|_{H^1}\|\mathfrak{R}_{1}\|_{H^1}\big],
\end{align}
and
\begin{align}\label{5.2-3}
	&\varepsilon^2\mu\|\chi_{1}\sqrt{\hat{\rho}_{0}}\nabla_{z}\hat{D}_{jl}^2\hat{u}_{R}\|_{L^2}^2+\varepsilon^2\zeta\|\chi_{1}\sqrt{\hat{\rho}_{0}}\hat{D}_{jl}^2\operatorname{div}_{z}\hat{u}_{R}\|_{L^2}^2+\varepsilon^2\kappa\|\chi_{1}\sqrt{\frac{\hat{\rho}_{0}}{\hat{\theta}_{0}}}\nabla_{z}\hat{D}_{jl}^2\hat{\theta}_{R}\|_{L^2}^2\nonumber\\
	&\leq \big[C(\mathcal{I}_{1})C_{\sigma}\v^2+\tau\big]\|\rho_{R}\|_{H^2}^2+C_{\tau,\sigma}\v^2\|u_{R}\|_{H^2}^2+C_{\sigma}\|\theta_{R}\|_{H^2}^2\nonumber\\
	&\quad+C\v^2\big[\|(\mathfrak{R}_{2},\mathfrak{R}_{3})\|_{H^1}^2+\|\rho_{R}\|_{H^2}\|\mathfrak{R}_{1}\|_{H^2}\big]+C\sigma^2\v^2\|(u_{R},\theta_{R})\|_{H^3}^2,
\end{align}
where $C>0$ is independent of $\tau$, $\sigma$ and $\v$.
\end{lemma}

\noindent\textbf{Proof.} We first consider the estimate \eqref{5.2-2}. We shall multiply $\eqref{5.2-0}_{1}$, $\eqref{5.2-0}_{2}$ and $\eqref{5.2-0}_{3}$ by $\chi_{1}^2\hat{D}_{j}(\hat{\theta}_{0}\hat{\rho}_{R})$, $\chi_{1}^2\hat{D}_{j}(\hat{\rho}_{0}\hat{u}_{R}^{k})$ and $\chi_{1}^2\hat{D}_{j}(\hat{\rho}_{0}\hat{\theta}_{R})$ respectively, and add the resultant equations together. Almost all terms can be estimated similarly as in the proof of Lemma \ref{lem4.2}, but we would like to point out that the terms involving $\nabla_{z}\chi_{1}$ arising from the process of integrating by parts will lead to $C_{\sigma}\|(u_{R},\theta_{R})\|_{H^1}^2$. Here, we only need to focus on the estimates of new terms, that is, the terms caused by boundary function $w$.
%\begin{align}\label{5.5}
%&-\frac{1}{\varepsilon}\int_{\hat{\Omega}\cap\hat{Q}}\chi_{1}^2[\hat{D}_{i\tau\xi}^3(\hat{\rho}_{0}\hat{u}_{R}^{i})\hat{D}_{\tau\xi}^2(\hat{\rho}_{0}\hat{\theta}_{R}+\hat{\theta}_{0}\hat{\rho}_{R})-\hat{D}_{\tau\xi k}^3(\hat{\rho}_{0}\hat{\theta}_{R}+\hat{\theta}_{0}\hat{\rho}_{R})\hat{D}_{\tau\xi}^2(\hat{\rho}_{0}\hat{u}_{R}^{k})]\,{\rm d}x\nonumber\\
%&=-\frac{1}{\varepsilon}\int_{\hat{\Omega}\cap\hat{Q}}2\chi_{1}\hat{D}_{k}\chi_{1}\hat{D}_{\tau\xi}^2(\hat{\rho}_{0}\hat{\theta}_{R}+\hat{\theta}_{0}\hat{\rho}_{R})\hat{D}_{j\xi}^2(\hat{\rho}_{0}\hat{u}_{R}^{k})\,{\rm d}z\nonumber\\
%&\leq \delta\|\frac{\nabla(\hat{\rho}_{0}\hat{\theta}_{R}+\hat{\theta}_{0}\hat{\rho}_{R})}{\varepsilon}\|_{H^1}^2+C_\delta\|\hat{u}_{R}\|_{H^2}^2,
%\end{align}
First, we integrate by parts to see that
\begin{align}\label{5.6}
	&\frac{1}{\varepsilon}\int_{\hat{\Omega}\cap \hat{Q}}\chi_{1}^2\sum\limits_{m=1}^2\Big[\hat{D}_{j}(\hat{D}_{m}w\hat{D}_{3}(\hat{\rho}_{0}\hat{u}_{R}^m))\hat{D}_{j}(\hat{\rho}_{0}\hat{\theta}_{R}+\hat{\theta}_{0}\hat{\rho}_{R})+\hat{D}_{m}w\hat{D}_{3j}^2(\hat{\rho}_{0}\hat{\theta}_{R}+\hat{\theta}_{0}\hat{\rho}_{R})\hat{D}_{j}(\hat{\rho}_{0}\hat{u}_{R}^{m})\Big]\,{\rm d}z\nonumber\\
	&\quad +\frac{1}{\varepsilon}\int_{\hat{\Omega}\cap \hat{Q}}\chi_{1}^2\sum\limits_{m=1}^2\hat{D}_{jm}^2w\hat{D}_{3}(\hat{\rho}_{0}\hat{\theta}_{R}+\hat{\theta}_{0}\hat{\rho}_{R})\hat{D}_{j}(\hat{\rho}_{0}\hat{u}_{R}^{m})\,{\rm d}z\nonumber\\
	&=-\frac{1}{\varepsilon}\int_{\hat{\Omega}\cap \hat{Q}}2\chi_{1}\hat{D}_{3}\chi_{1}\sum\limits_{m=1}^2\Big[\hat{D}_{m}w\hat{D}_{j}(\hat{\rho}_{0}\hat{\theta}_{R}+\hat{\theta}_{0}\hat{\rho}_{R})\hat{D}_{j}(\hat{\rho}_{0}\hat{u}_{R}^{m})\Big]\,{\rm d}z\nonumber\\
	%&\quad -\frac{1}{\varepsilon}\int_{\hat{\Omega}\cap \hat{Q}}\chi_{1}^2\sum\limits_{m=1}^2[\hat{D}_{m3}w\hat{D}_{j}(\hat{\rho}_{0}\hat{\theta}_{R}+\hat{\theta}_{0}\hat{\rho}_{R})\hat{D}_{j}(\hat{\rho}_{0}\hat{u}_{R}^{m})]\,{\rm d}z\nonumber\\
	&\quad +\frac{1}{\varepsilon}\int_{\hat{Q}\cap \hat{\Omega}}\chi_{1}^2\sum\limits_{m=1}^2\hat{D}_{jm}^2w\Big[\hat{D}_{3}(\hat{\rho}_{0}\hat{u}_{R}^{m})\hat{D}_{j}(\hat{\rho}_{0}\hat{\theta}_{R}+\hat{\theta}_{0}\hat{\rho}_{R})+\hat{D}_{3}(\hat{\rho}_{0}\hat{\theta}_{R}+\hat{\theta}_{0}\hat{\rho}_{R})\hat{D}_{j}(\hat{\rho}_{0}\hat{u}_{R}^{m})\Big]\,{\rm d}z\nonumber\\
	&\leq \tau\|\frac{\nabla (\rho_{0}\theta_{R}+\theta_{0}\rho_{R})}{\v}\|_{L^2}+C_{\tau,\sigma}\|u_{R}\|_{H^1}^2
\end{align}
for any $\tau>0$, where we have used the vanishing of $\hat{D}_{j}(\hat{\rho}_{0}\hat{u}_{R})$ on $\partial\hat{\Omega}\cap\hat{Q}$ for $j=1,2$.
%Secondly, by using \eqref{5.4} and integrating by parts that
Similarly, by integrating by parts and using \eqref{5.4-1}, one has
\begin{align}\label{5.6-1}
&\int_{\hat{\Omega}\cap \hat{Q}}\chi_{1}^2[\sum\limits_{k=1}^2\hat{D}_{j}\hat{W}_{2}^{k}\hat{D}_{j}(\hat{\rho}_{0}\hat{u}_{R}^{k})
+\hat{D}_{j}\hat{W}_{3}\hat{D}_{j}(\hat{\rho}_{0}\hat{\theta}_{R})]\,{\rm d}z\nonumber\\
&\leq \frac{\mu}{16}\|\chi_{1}\sqrt{\hat{\rho}_{0}}\nabla_{z}\hat{D}_{j}\hat{u}_{R}\|_{L^2}^2+\frac{\kappa}{16}\|\chi_{1}\sqrt{\frac{\hat{\rho}_{0}}{\hat{\theta}_{0}}}\nabla_{z} \hat{D}_{j}\hat{\theta}_{R}\|_{L^2}^2+C_{\sigma}\|(u_{R},\theta_{R})\|_{H^1}^2 +C\|(\hat{W}_{2},\hat{W}_{3})\|_{L^2}^2\nonumber\\
&\leq \frac{\mu}{16}\|\chi_{1}\sqrt{\rho}_{0}\nabla_{z}\hat{D}_{j}\hat{u}_{R}\|_{L^2}^2+\frac{\kappa}{16}\|\chi_{1}\sqrt{\frac{\rho_{0}}{\theta_{0}}}\nabla_{z} \hat{D}_{j}\hat{\theta}_{R}\|_{L^2}^2+C_{\sigma}\|(u_{R},\theta_{R})\|_{H^1}^2\nonumber\\
&\quad +C(\mathcal{I}_{1})\sigma^2\|\rho_{R}\|_{H^1}^2+C\sigma^2\|(u_{R},\theta_{R})\|_{H^2}^2.
\end{align}
For the term involving $\hat{W}_{1}$, similar to \eqref{4.26-5}, we have from \eqref{5.4-1} that
\begin{align}\label{5.6-2}
&\int_{\hat{\Omega}\cap \hat{Q}}\chi_{1}^2\hat{D}_{j}\hat{W}_{1}\hat{D}_{j}(\hat{\theta}_{0}\hat{\rho}_{R})\,{\rm d}z\leq C\|\hat{W}_{1}\|_{H^1}\|\rho_{R}\|_{H^1}\leq C(\mathcal{I}_{1})\sigma \|\rho_{R}\|_{H^2}^2.
\end{align}
Combining \eqref{5.6}--\eqref{5.6-2} and applying similar arguments as in Lemma \ref{lem4.2}, we deduce \eqref{5.2-2}.

Next, we consider the estimate \eqref{5.2-3}. Similar to \eqref{4.27}, we apply $\hat{D}_{jl}^2$ to \eqref{5.2}, and then multiply the first, second and third equation of the resultant system by $\v^2\chi_{1}^2\hat{D}_{jl}^2(\theta_{0}\rho_{R})$, $\v^2\chi_{1}^2\hat{D}_{jl}^2(\rho_{0}u_{R}^{k})$ and $\v^2\chi_{1}^2\hat{D}_{jl}^2(\rho_{0}\theta_{R})$ respectively. Most of the calculations of \eqref{5.2-3} are very similar to \eqref{5.6}--\eqref{5.6-1} and the proof of Lemma \ref{lem4.3} except the following term:
\begin{align}\label{5.6-3}
&\v^2\int_{\hat{Q}\cap \hat{\Omega}}\chi_{1}^2\hat{D}_{jl}^2(\hat{\theta}_{0}\hat{\rho}_{R})\hat{D}_{jl}^2\Big[\operatorname{div}_{z}[\hat{\rho}_{R}(\hat{u}_{1}+\v(\hat{u}_{2}+\hat{\tilde{u}}_{R}))]-\hat{W}_{1}\Big]\,{\rm d}z\nonumber\\
&=\v^2\int_{\hat{Q}\cap \hat{\Omega}}\chi_{1}^2\hat{D}_{jl}^2(\hat{\theta}_{0}\hat{\rho}_{R})\hat{D}_{jl}^2\operatorname{div}_{z}[\hat{\rho}_{R}(\hat{u}_{1}+\v(\hat{u}_{2}+\hat{\tilde{u}}_{R}))]\nonumber\\
&\quad -\sum\limits_{m=1}^2\v^2\int_{\hat{Q}\cap \hat{\Omega}}\chi_{1}^2\hat{D}_{jl}^2(\hat{\theta}_{0}\hat{\rho}_{R})\hat{D}_{jl}^2[\hat{D}_{m}w(\hat{D}_{3}(\hat{\rho}_{R}(\hat{u}_{1}^{m}+\v\hat{u}_{2}^m+\v\hat{\tilde{u}}_{R}^m)))]\,{\rm d}z.
\end{align}
For terms involving with $\hat{u}_{2}$ and $\hat{\tilde{u}}_{R}$, applying the similar  integrating by parts arguments in \eqref{4.27-5}, they can be controlled by $C(\mathcal{I}_{1})C_{\sigma}\v^2\|\rho_{R}\|_{H^2}^2$ due to $(\hat{u}_{2},\hat{\tilde{u}}_{R})\vert_{z_{3}=0}=(0,0)$.
%Unlike the arguments used in \eqref{4.27-5},

For terms involving $\hat{u}_{1}$ in \eqref{5.6-3}, since we only have $u_{1}\cdot \vec{n}\vert_{\partial\Omega}=0$ in $x-$ coordinate, we should be very careful when integrating by parts to remove the terms involving the third order derivative of $\hat{\rho}_{R}$ in $z-$coordinate, that is, the following term
%$z-$coordinate, %since we only have $\hat{u}_{1}\cdot \vec{n}=0$ on $\partial \Omega$ in $x-$coordinate.
%The terms involved with the third derivatives of $\rho_{R}$ in \eqref{5.6-3} is
\begin{align}\label{5.6-5}
&\int_{\hat{Q}\cap\hat{\Omega}}\chi_{1}^2\hat{D}_{jl}^2(\hat{\theta}_{0}\hat{\rho}_{R})[\hat{u}_{1}\cdot \nabla_{z} (\hat{D}_{jl}^2\hat{\rho}_{R})-\sum\limits_{m=1}^2\hat{D}_{m}w\hat{u}_{1}^{m}\hat{D}_{3}(\hat{D}_{jl}^2\hat{\rho}_{R})]\,{\rm d}z.
\end{align}
Noting that $\vec{n}=\frac{(-\partial_{x_{1}}w,-\partial_{x_{2}}w, 1)}{1+(\partial_{x_{1}}w)^2+(\partial_{z_{2}}w)^2}$ and $u_{1}\cdot \vec{n}\vert_{\partial\Omega}=0$ in $x-$coordinate, one has $-\hat{D}_{1}w\hat{u}_{1}^{1}-\hat{D}_{2}w\hat{u}_{1}^{2}+\hat{u}_{1}^{3}=0$ on $z_{3}=0$.
%By using the fact that $\operatorname{supp}\chi_{1}\subset \hat{Q}$, one has
It is clear that $\int_{\hat{Q}\cap \hat{\Omega}}\hat{D}_{m}(\chi_{1}^2\cdots)\,{\rm d}z=0$ for $m=1,2$ due to $\operatorname{supp}\chi_{1}\subset \hat{Q}$. Now we integrate by parts to obtain
\begin{align}\label{5.6-6}
&\v^2\int_{\hat{Q}\cap\hat{\Omega}}\chi_{1}^2\hat{D}_{jl}^2(\hat{\theta}_{0}\hat{\rho}_{R})[\hat{u}_{1}\cdot \nabla_{z} (\hat{D}_{jl}^2\hat{\rho}_{R})-\sum\limits_{m=1}^2\hat{D}_{m}w\hat{u}_{1}^{m}\hat{D}_{3}(\hat{D}_{jl}^2\hat{\rho}_{R})]\,{\rm d}z\nonumber\\
&=\frac{1}{2}\v^2\int_{\hat{Q}\cap\hat{\Omega}}\chi_{1}^2\hat{\theta}_{0}[\hat{u}_{1}\cdot \nabla_{z} (\hat{D}_{jl}^2\hat{\rho}_{R})^2-\sum\limits_{m=1}^2\hat{D}_{m}w\hat{u}_{1}^{m}\hat{D}_{3}(\hat{D}_{jl}^2\hat{\rho}_{R})^2]\,{\rm d}z\nonumber\\
&\quad
+\v^2\int_{\hat{Q}\cap\hat{\Omega}}\chi_{1}^2[\hat{u}_{1}\cdot \nabla_{z} (\hat{D}_{jl}^2\hat{\rho}_{R})-\sum\limits_{m=1}^2\hat{D}_{m}w\hat{u}_{1}^{m}\hat{D}_{3}(\hat{D}_{jl}^2\hat{\rho}_{R})](\hat{D}_{jl}^2(\hat{\theta}_{0}\hat{\rho}_{R})-\hat{\theta}_{0}\hat{D}_{jl}^2\hat{\rho}_{R})\,{\rm d}z\nonumber\\
&=-\frac{1}{2}\v^2\int_{\hat{Q}\cap\partial \hat{\Omega}}\chi_{1}^2\hat{\theta}_{0}(\hat{D}_{jl}^2\hat{\rho}_{R})^2\cancel{\big[\hat{u}_{1}^3-\sum\limits_{m=1}^2\hat{D}_{m}w\hat{u}_{1}^{m}\big]}\,{\rm d}z\nonumber\\
&\quad -\v^2\int_{\hat{Q}\cap\partial \hat{\Omega}}\chi_{1}^2\hat{D}_{jl}^2(\hat{\rho}_{R})\cancel{\big[\hat{u}_{1}^3-\sum\limits_{m=1}^2\hat{D}_{m}w\hat{u}_{1}^{m}\big]}(\hat{D}_{jl}^2(\hat{\theta}_{0}\hat{\rho}_{R})-\hat{\theta}_{0}\hat{D}_{jl}^2\hat{\rho}_{R})\,{\rm d}z\nonumber\\
&\quad -\frac{1}{2}\v^2\int_{\hat{Q}\cap\hat{\Omega}}(\hat{D}_{jl}^2\hat{\rho}_{R})^2\operatorname{div}_{z}\big[\chi_{1}^2\hat{\theta}_{0}\hat{u}_{1}\big]\,{\rm d}z +\frac{1}{2}\sum\limits_{m=1}^2\int_{\hat{Q}\cap\hat{\Omega}}(\hat{D}_{jl}^2\hat{\rho}_{R})^2\hat{D}_{3}\big[\chi_{1}^2\hat{\theta}_{0}\hat{D}_{m}w\hat{u}_{1}^{m}\big]\,{\rm d}z\nonumber\\
&\quad
-\v^2\int_{\hat{Q}\cap\hat{\Omega}}\hat{D}_{jl}^2\hat{\rho}_{R}\operatorname{div}_{z}\big[\chi_{1}^2\hat{u}_{1}(\hat{D}_{jl}^2(\hat{\theta}_{0}\hat{\rho}_{R})-\hat{\theta}_{0}\hat{D}_{jl}^2\hat{\rho}_{R})\big]\,{\rm d}z\nonumber\\
&\quad
-\sum\limits_{m=1}^2\v^2\int_{\hat{Q}\cap\hat{\Omega}}\hat{D}_{jl}^2\hat{\rho}_{R}\hat{D}_{3}\big[\chi_{1}^2\hat{D}_{m}w\hat{u}_{1}^{m}(\hat{D}_{jl}^2(\hat{\rho}_{R})-\hat{\theta}_{0}\hat{D}_{jl}^2\hat{\rho}_{R})\big]\,{\rm d}z\nonumber\\
&\leq C_{\sigma}\v^2\|\rho_{R}\|_{H^2}^2\|u_{1}\|_{H^3}\leq C(\mathcal{I}_{1})C_{\sigma}\v^2\|\rho_{R}\|_{H^2}^2.
\end{align}
Therefore the proof of Lemma \ref{lem5.0} is complete. $\hfill\square$
\smallskip

\noindent\textbf{Step 2. Estimates of $\varepsilon\|\hat{D}_{3}\operatorname{div}_{z}\hat{u}_{R}\|_{L^2}^2$}. Since  $\partial_{z_{3}}\hat{u}_{R}\vert_{z_{3}=0}$ and $\partial_{z_{3}}\hat{\theta}_{R}\vert_{z_{3}=0}$ are unknown,
%the normal derivatives of $u_{R}$ and $\theta_{R}$,
we cannot apply the standard energy method directly to obtain higher order normal derivatives.
%get the estimates of more than second order normal derivatives.
Motivated by \cite{Choe-Jin, Dou-Jiang-Jiang-Yang, Li-Liao-2019}, in order to get $\hat{D}_{3}\operatorname{div}u_{R}$, we
have from
%denote $\hat{e}_{3}=(0,0,1)$ and note that
$\eqref{5.2}_{2}\cdot \hat{e}_{3}$ that
\begin{align}\label{5.7}
	&-(\mu+\zeta)\hat{D}_{3}(\operatorname{div}_{z}\hat{u}_{R})+\frac{1}{\varepsilon}\hat{D}_{3}(\hat{\rho}_{0}\hat{\theta}_{R}+\hat{\rho}_{R}\hat{\theta}_{0})\nonumber\\
	&=-\hat{\rho}_{0}(\hat{u}_{1}\cdot \nabla_{z} \hat{u}_{R}+\hat{u}_{R}\cdot \nabla_{z} \hat{u}_{1})\cdot \hat{e}_{3}-\hat{D}_{3}(\hat{\rho}_{R}\hat{\theta}_{1}+\hat{\rho}_{1}\hat{\theta}_{R})+\mu(\Delta_{z}\hat{u}_{R}-\nabla_{z}\operatorname{div}_{z}\hat{u}_{R})\cdot \hat{e}_{3}\nonumber\\
	&\quad -\hat{W}_{2}^3-\hat{\tilde{F}}^{\varepsilon,3}(\hat{\rho}_{R},\hat{u}_{R},\hat{\theta}_{R})-\hat{\mathfrak{R}}_{2}^{3},
\end{align}
where $\hat{e}_{3}=(0,0,1)^{t}$. It is worth pointing out that
\begin{align}\label{5.8}
	\mu(\Delta_{z}\hat{u}_{R}-\nabla_{z}\operatorname{div}_{z}\hat{u}_{R})\cdot \hat{e}_{3}=\mu[\hat{D}_{1}^2\hat{u}_{R}^{3}+\hat{D}_{2}^2\hat{u}_{R}^{3}-\hat{D}_{13}^2\hat{u}_{R}^{1}-\hat{D}_{23}^2\hat{u}_{R}^2],
\end{align}
where $\hat{D}_{33}^2\hat{u}_{R}^{3}$ is canceled.

\begin{lemma}\label{lem5.0-1}
For any $\tau>0$, there exist two positive constants $C_{\sigma}$, which depends on $\sigma$, and $C_{\tau,\sigma}$, which depends on $\tau$ and $\sigma$, both independent of $\varepsilon$, such that
\begin{align}\label{5.7-17}
	&(\mu+\zeta)\int_{\hat{Q}\cap \hat{\Omega}}\chi_{1}^2\hat{\rho}_{0}|\hat{D}_{3}\operatorname{div}_{z}\hat{u}_{R}|^2\,{\rm d}z+\kappa\int_{\hat{Q}\cap \hat{\Omega}}\chi_{1}^2\frac{\hat{\rho}_{0}}{\hat{\theta}_{0}}|\hat{D}_{33}^2\hat{\theta}_{R}|^2\,{\rm d}z\nonumber\\
	&\leq (\tau+\sigma)\|\frac{\nabla(\rho_{0}\theta_{R}+\theta_{0}\rho_{R})}{\varepsilon}\|_{L^2}^2+(C(\mathcal{I}_{1})+C\sigma)\|(u_{R},\theta_{R})\|_{H^2}^2+C\sum\limits_{i=1}^2\|\nabla_{z} \hat{D}_{i}(\hat{u}_{R},\hat{\theta}_{R})\|_{L^2}^2\nonumber\\
	&\quad +C(\mathcal{I}_{1})\|\rho_{R}\|_{H^2}^2+C\|\mathfrak{R}_{1}\|_{H^1}\|(\rho_{R},\theta_{R})\|_{H^2}
	+C\|(\mathfrak{R}_{2},\mathfrak{R}_{3})\|_{L^2}^2+C_{\tau,\sigma}\|u_{R}\|_{H^1}^2,
\end{align}
where $C>0$ is independent of $\tau$, $\sigma$ and $\v$.
\end{lemma}

\noindent\textbf{Proof}. We multiply \eqref{5.7} by $-\chi_{1}^2\hat{D}_{3}\operatorname{div}_{z}(\hat{\rho}_{0}\hat{u}_{R})$ to get
\begin{align}\label{5.7-1}
&(\mu+\zeta)\int_{\hat{Q}\cap \hat{\Omega}}\hat{\rho}_{0}\chi_{1}^2|\hat{D}_{3}\operatorname{div}_{z}\hat{u}_{R}|^2\,{\rm d}z-\frac{1}{\varepsilon}\int_{\hat{Q}\cap \hat{Q}}\chi_{1}^2\hat{D}_{3}(\hat{\rho}_{0}\hat{\theta}_{R}+\hat{\rho}_{R}\hat{\theta}_{0})\hat{D}_{3}\operatorname{div}_{z}(\hat{\rho}_{0}\hat{u}_{R})\,{\rm d}z\nonumber\\
&=-(\mu+\zeta)\int_{\hat{\Omega}\cap \hat{Q}}\chi_{1}^2\hat{D}_{3}(\hat{u}_{R}\cdot \nabla \hat{\rho}_{0})\hat{D}_{3}\operatorname{div}_{z}\hat{u}_{R}\,{\rm d}z-(\mu+\zeta)\int_{\hat{\Omega}\cap \hat{Q}}\chi_{1}^2\hat{D}_{3}\operatorname{div}_{z}\hat{u}_{R}(\hat{D}_{3}\hat{\rho}_{0})\operatorname{div}_{z}\hat{u}_{R}\,{\rm d}z\nonumber\\
&\quad +\int_{\hat{Q}\cap \hat{\Omega}}\chi_{1}^2\hat{D}_{3}\operatorname{div}_{z}(\hat{\rho}_{0}\hat{u}_{R})[\hat{\rho}_{0}(\hat{u}_{R}\cdot \hat{u}_{1}\cdot \nabla_{z}\hat{u}_{R}+\hat{u}_{R}\cdot \nabla \hat{u}_{1})\cdot \hat{e}_{3}]\,{\rm d}z\nonumber\\
&\quad +\int_{\hat{Q}\cap \hat{\Omega}}\chi_{1}^2\hat{D}_{3}\operatorname{div}(\hat{\rho}_{0}\hat{u}_{R})\hat{D}_{3}[\hat{\rho}_{R}\hat{\theta}_{1}+\hat{\rho}_{1}\hat{\theta}_{R}]\,{\rm d}z+\int_{\hat{Q}\cap \hat{\Omega}}\chi_{1}^2\hat{D}_{3}\operatorname{div}(\hat{\rho}_{0}\hat{u}_{R})[\hat{W}_{2}^3+\hat{F}^{\varepsilon,3}+\hat{\mathfrak{R}}_{2}^{3}]\,{\rm d}z\nonumber\\
&\quad -\mu\int_{\hat{Q}\cap \hat{\Omega}}\chi_{1}^2\hat{D}_{3}\operatorname{div}(\hat{\rho}_{0}\hat{u}_{R})[(\Delta_{z}\hat{u}_{R}-\nabla_{z}\operatorname{div}_{z}\hat{u}_{R})\cdot \hat{e}_{3}]\,{\rm d}z=:\sum\limits_{m=1}^6I_{m}.
\end{align}
For $I_{m}(m=1,2,\cdots,5)$, using the H\"{o}lder inequality, one has that
\begin{align}\label{5.7-2}
	\sum\limits_{m=1}^{5}|I_{m}|&\leq \frac{\mu+\zeta}{4}\int_{\hat{Q}\cap \hat{\Omega}}\hat{\rho}_{0}\chi_{1}^2|\hat{D}_{3}\operatorname{div}_{z}\hat{u}_{R}|^2\,{\rm d}z+C(\mathcal{I}_{0})\|(\rho_{R},u_{R},\theta_{R})\|_{H^1}^2\nonumber\\
	&\quad+C(\|\hat{W}_{2}^3\|_{L^2}^2+\|\hat{F}^{\var,3}\|_{L^2}^2+C\|\mathfrak{R}_{2}^3\|_{L^2}^2)\nonumber\\
	&\leq \frac{\mu+\zeta}{4}\int_{\hat{Q}\cap \hat{\Omega}}\hat{\rho}_{0}\chi_{1}^2|\hat{D}_{3}\operatorname{div}_{z}\hat{u}_{R}|^2\,{\rm d}z +C(\mathcal{I}_{1})\|(\rho_{R},u_{R},\theta_{R})\|_{H^2}^2\nonumber\\
	&\quad +C\sigma^2\|(u_{R},\theta_{R})\|_{H^2}^2+C\|\mathfrak{R}_{2}\|_{L^2}^2,
\end{align}
where we have used \eqref{5.4-1}--\eqref{5.5} and \eqref{5.5} in the last inequality.
For $I_{6}$, noting \eqref{5.8}, we have
\begin{align}\label{5.7-3}
|I_{6}|&\leq \frac{\mu+\zeta}{4}\int_{\hat{Q}\cap \hat{\Omega}}\hat{\rho}_{0}\chi_{1}^2|\hat{D}_{3}\operatorname{div}_{z}\hat{u}_{R}|^2\,{\rm d}z+C(\mathcal{I}_{0})\|u_{R}\|_{H^1}^2+C\sum\limits_{i=1}^2\|\nabla_{z} \hat{D}_{i}\hat{u}_{R}\|_{L^2}^2.
\end{align}
Substituting \eqref{5.7-2}--\eqref{5.7-3} into \eqref{5.7-1} to get
\begin{align}\label{5.7-4}
&(\mu+\zeta)\int_{\hat{Q}\cap \hat{\Omega}}\hat{\rho}_{0}\chi_{1}^2|\hat{D}_{3}\operatorname{div}_{z}\hat{u}_{R}|^2\,{\rm d}z-\frac{1}{\varepsilon}\int_{\hat{Q}\cap \hat{Q}}\chi_{1}^2\hat{D}_{3}(\hat{\rho}_{0}\hat{\theta}_{R}+\hat{\rho}_{R}\hat{\theta}_{0})\hat{D}_{3}\operatorname{div}_{z}(\hat{\rho}_{0}\hat{u}_{R})\,{\rm d}z\nonumber\\
&\leq C\sum\limits_{i=1}^2\|\nabla_{z} \hat{D}_{i}\hat{u}_{R}\|_{L^2}^2+C(\mathcal{I}_{1})\|(\rho_{R},u_{R},\theta_{R})\|_{H^2}^2+C\sigma^2\|(u_{R},\theta_{R})\|_{H^2}^2+C\|\mathfrak{R}_{2}\|_{L^2}^2.
%&\leq C(\|u_{R}\|_{H^1}^2+\|\theta_{R}\|_{H^1}^2+\|r_{2}\|_{L^2}^2)+C\|\rho_{R}\|_{H^1}^2\|\theta_{1}\|_{H^2}^2+C\sum\limits_{i=1}^2\|\nabla \hat{D}_{i}\hat{u}_{R}\|_{L^2}^2\nonumber\\
%&\quad +C\varepsilon^2\|u_{R}\|_{H^2}^2(\|u_{1}\|_{H^2}^2+\|u_{2}\|_{H^2}^2+\|\tilde{u}_{R}\|_{H^2}^2)+C\varepsilon^2\|\theta_{R}\|_{H^2}^2(\|\rho_{2}\|_{H^2}^2+\|\rho_{3}\|_{H^2}^2)\nonumber\\
%&\quad +C\varepsilon^2\|\rho_{R}\|_{H^2}^2((\|u_{1}\|_{H^2}^4+\|\varepsilon u_2\|_{H^2}^4+\|\varepsilon\tilde{u}_{R}\|_{H^2}^4)+\|\tilde{\theta}_{R}\|_{H^2}^2)\nonumber\\
%&\quad +C\sigma(\|\rho_{R}\|_{H^1}^2+\|u_{R}\|_{H^2}^2+\|\theta_{R}\|_{H^2}^2).
\end{align}
To cope with the singular term on the left hand side (LHS) of \eqref{5.7-4}, we
%first apply
%$\hat{D}_{3}$ to $\eqref{5.2}_{1}$ and multiply the resultant equation by $\chi_{1}^2\hat{D}_{3}(\hat{\theta}_{0}\hat{\rho}_{R})$. In the meanwhile, we multiply $\eqref{5.2}_{3}$ by $\chi_{1}^2\hat{D}_{33}^2(\hat{\rho}_{0}\hat{\theta}_{R})$, and then adding the resultant equations together to see
multiply $\eqref{5.2}_{1}$  and $\eqref{5.2}_{3}$ by  $-\chi_{1}^2\hat{D}_{33}^2(\hat{\theta}_{0}\hat{\rho}_{R})$ and $\chi_{1}^2\hat{D}_{33}^2(\hat{\rho}_{0}\hat{\theta}_{R})$ respectively to get
\begin{align}\label{5.7-5}
&\kappa\int_{\hat{Q}\cap \hat{\Omega}}\chi_{1}^2\frac{\hat{\rho}_{0}}{\hat{\theta}_{0}}|\hat{D}_{33}^2\hat{\theta}_{R}|^2\,{\rm d}z
-\frac{1}{\v}\int_{\hat{Q}\cap \hat{\Omega}}\chi_{1}^2\operatorname{div}_{z}(\hat{\rho}_{0}\hat{u}_{R})\hat{D}_{33}^2(\hat{\rho}_{0}\hat{\theta}_{R}+\hat{\theta}_{0}\hat{\rho}_{R})\,{\rm d}z\nonumber\\
&\leq \frac{\kappa}{16}\int_{\hat{Q}\cap \hat{\Omega}}\chi_{1}^2\frac{\hat{\rho}_{0}}{\hat{\theta}_{0}}|\hat{D}_{33}^2\hat{\theta}_{R}|^2\,{\rm d}z+C\|\theta_{R}\|_{H^1}^2
-\kappa\sum\limits_{i=1}^2\int_{\hat{Q}\cap \hat{\Omega}}\chi_{1}^2\frac{\hat{D}_{ii}^2\hat{\theta}_{R}}{\hat{\theta}_{0}}\hat{D}_{33}^2(\hat{\rho}_{0}\hat{\theta}_{R})\,{\rm d}x\nonumber\\
 &\quad +\frac{P_{1}}{\hat{\theta}_{0}}\int_{\hat{Q}\cap \hat{\Omega}}\chi_{1}^2\operatorname{div}_{z}\hat{u}_{R}\hat{D}_{33}^2(\hat{\rho}_{0}\hat{\theta}_{R})\,{\rm d}z-\frac{1}{\varepsilon}\sum\limits_{m=1}^2\int_{\hat{Q}\cap \hat{\Omega}}\chi_{1}^2\hat{D}_{m}w\hat{D}_{3}(\hat{\rho}_{0}\hat{u}_{R}^{m})\hat{D}_{33}^2(\hat{\rho}_{0}\hat{\theta}_{R}+\hat{\theta}_{0}\hat{\rho}_{R})\,{\rm d}z\nonumber\\
 &\quad +\int_{\hat{Q}\cap \hat{\Omega}}\chi_{1}^2\frac{\hat{\rho}_{0}\hat{\theta}_{R}+\hat{\rho}_{R}\hat{\theta}_{0}}{\hat{\theta}_{0}}\operatorname{div}_{z}\hat{u}_{1}\hat{D}_{33}^2(\hat{\rho}_{0}\hat{\theta}_{R})\,{\rm d}z+\int_{\hat{Q}\cap \hat{\Omega}}\chi_{1}^2\hat{D}_{33}^2(\hat{\rho}_{0}\hat{\theta}_{R})[\hat{\mathfrak{R}}_{3}+\hat{W}_{3}+\frac{1}{\hat{\theta}_{0}}\hat{G}^{\var}]\,{\rm d}z\nonumber\\
%+\int_{\hat{Q}\cap \hat{\Omega}}\chi_{1}^2\hat{D}_{33}^2(\hat{\rho}_{0}\hat{\theta}_{R})\frac{1}{\hat{\theta}_{0}}\hat{\Psi}(\nabla (u_{1}+\varepsilon(u_{2}+\tilde{u}_{R})))\nonumber\\
&\quad +\int_{\hat{Q}\cap \hat{\Omega}}\chi_{1}^2\hat{D}_{33}^2(\hat{\rho}_{0}\hat{\theta}_{R})\frac{1}{\hat{\theta}_{0}}\big[\hat{\rho}_{0}(\hat{u}_{R}\cdot \nabla_{z}\hat{\theta}_{1}+\hat{u}_{1}\cdot \nabla_{z}\hat{\theta}_{R})+\hat{\rho}_{1}\hat{u}_{R}\cdot \nabla_{z}\hat{\theta}_{0}+\hat{\rho}_{R}\hat{u}_{1}\cdot \nabla_{z}\hat{\theta}_{0}\big]\,{\rm d}z\nonumber\\
&\quad +\int_{\hat{Q}\cap \hat{\Omega}}\chi_{1}^2\hat{D}_{33}^2(\hat{\theta}_{0}\hat{\rho}_{R})\operatorname{div}_{z}[\hat{\rho}_{R}(\hat{u}_{1}+\varepsilon(\hat{u}_{2}+\hat{\tilde{u}}_{R}))]\,{\rm d}z\nonumber\\
&\quad +\int_{\hat{Q}\cap \hat{\Omega}}\chi_{1}^2\hat{D}_{33}^2(\hat{\theta}_{0}\hat{\rho}_{R})(\hat{W}_{1}+\hat{\mathfrak{R}}_{1})\,{\rm d}z=:\sum\limits_{m=1}^{10}J_{m}.
\end{align}
For the second term on the LHS of \eqref{5.7-5}, we integrate by parts to see
\begin{align}\label{5.7-6}
&-\frac{1}{\varepsilon}\int_{\hat{Q}\cap \hat{\Omega}}\chi_{1}^2\operatorname{div}_{z}(\hat{\rho}_{0}\hat{u}_{R})\hat{D}_{33}^2(\hat{\rho}_{0}\hat{\theta}_{R}+\hat{\theta}_{0}\hat{\rho}_{R})\,{\rm d}z-\frac{1}{\varepsilon}\int_{\hat{Q}\cap \hat{\Omega}}\chi_{1}^2\hat{D}_{3}\operatorname{div}_{z}(\hat{\rho}_{0}\hat{u}_{R})\hat{D}_{3}(\hat{\rho}_{0}\hat{\theta}_{R}+\hat{\theta}_{0}\hat{\rho}_{R})\,{\rm d}z\nonumber\\
&=\frac{1}{\varepsilon}\int_{\hat{Q}\cap \partial \hat{\Omega}}\chi_{1}^2\operatorname{div}_{z}(\hat{\rho}_{0}\hat{u}_{R})\hat{D}_{3}(\hat{\rho}_{0}\hat{\theta}_{R}+\hat{\theta}_{0}\hat{\rho}_{R})\,{\rm d}z_{1}{\rm d}z_{2}\nonumber\\ &\quad +\frac{1}{\varepsilon}\int_{\hat{Q}\cap \hat{\Omega}}\chi_{1}\hat{D}_{3}\chi_{1}\operatorname{div}_{z}(\hat{\rho}_{0}\hat{u}_{R})\hat{D}_{3}(\hat{\rho}_{0}\hat{\theta}_{R}+\hat{\theta}_{0}\hat{\rho}_{R})\,{\rm d}z.
\end{align}
For the first term on the RHS of \eqref{5.7-6}, we use $\eqref{5.2}_{1}$ to see that
\begin{align}\label{5.7-7}
&\frac{1}{\varepsilon}\int_{\hat{Q}\cap \partial\hat{\Omega}}\chi_{1}^2\operatorname{div}_{z}(\hat{\rho}_{0}\hat{u}_{R})\hat{D}_{3}(\hat{\rho}_{0}\hat{\theta}_{R}+\hat{\theta}_{0}\hat{\rho}_{R})\,{\rm d}z_{1}{\rm d}z_{2}\nonumber\\
&=\frac{1}{\varepsilon}\sum\limits_{m=1}^2\int_{\hat{Q}\cap \partial\hat{\Omega}}\chi_{1}^2\hat{D}_{m}w\hat{D}_{3}(\hat{\rho}_{0}\hat{u}_{R}^{m})\hat{D}_{3}(\hat{\rho}_{0}\hat{\theta}_{R}+\hat{\theta}_{0}\hat{\rho}_{R})\,{\rm d}z_{1}{\rm d}z_{2}\nonumber\\
&\quad -\int_{\hat{Q}\cap \partial\hat{\Omega}}\chi_{1}^2\operatorname{div}_{z}[\hat{\rho}_{R}(\hat{u}_{1}+\varepsilon(\hat{u}_{2}+\hat{\tilde{u}}_{R}))]\hat{D}_{3}(\hat{\rho}_{0}\hat{\theta}_{R}+\hat{\theta}_{0}\hat{\rho}_{R})\,{\rm d}z_{1}{\rm d}z_{2}\nonumber\\
%&\quad -\int_{\hat{Q}\cap \partial \hat{\Omega}}\operatorname{div}_{z}[\hat{\tilde{u}}_{R}(\hat{\rho}_{1}+\varepsilon\hat{\rho}_{2}+\varepsilon^2\hat{\rho}_{3})]\hat{D}_{3}(\hat{\rho}_{0}\hat{\theta}_{R}+\hat{\theta}_{0}\hat{\rho}_{R})\,{\rm d}z_{1}{\rm d}z_{2}\nonumber\\
&\quad -\int_{\hat{Q}\cap \partial\hat{\Omega}}\chi_{1}^2(\hat{W}_{1}+\hat{\mathfrak{R}}_{1})\hat{D}_{3}(\hat{\rho}_{0}\hat{\theta}_{R}+\hat{\theta}_{0}\hat{\rho}_{R})\,{\rm d}z_{1}{\rm d}z_{2}.
\end{align}
Using trace inequality, we can bound the last two terms on the RHS of \eqref{5.7-7} are bounded by:
\begin{align}\label{5.7-10}
&C(\|\hat{\rho}_{R}\|_{H^2}+\|\hat{\theta}_{R}\|_{H^2})\Big[\|\hat{\rho}_{R}(u_{1}+\varepsilon(u_{2}+\hat{\tilde{u}}_{R}))\|_{H^2}
%+\|\hat{\tilde{u}}_{R}(\rho_{1}+\varepsilon\rho_{2}+\varepsilon^2\rho_{3})\|_{H^2}
+\|\hat{W}_{1}\|_{H^1}+\|\hat{\mathfrak{R}}_{1}\|_{H^1}\Big]\nonumber\\
&\leq C(\mathcal{I}_{1})(\|\rho_{R}\|_{H^2}^2+\|\theta_{R}\|_{H^2}^2)
%(\|u_{1}\|_{H^2}+\varepsilon\|u_{2}\|_{H^2}+\varepsilon\|\tilde{u}_{R}\|)\nonumber\\
+C\|(\rho_{R},\theta_{R})\|_{H^2}\|\mathfrak{R}_{1}\|_{H^1}.
%\big[(\|\rho_{1}\|_{H^2}+\varepsilon\|\rho_{2}\|_{H^2}+\varepsilon^2\|\rho_{3}\|_{H^2})\|\tilde{u}_{R}\|_{H^2}+\|\mathfrak{R}_{1}\|_{H^1}\big].
\end{align}
For the second term on the RHS of \eqref{5.7-6}, by the H\"{o}lder inequality, for any $\tau>0$, one has
\begin{align}\label{5.7-7-1}
\frac{1}{\varepsilon}\int_{\hat{Q}\cap \hat{\Omega}}\chi_{1}\hat{D}_{3}\chi_{1}\operatorname{div}_{z}(\hat{\rho}_{0}\hat{u}_{R})\hat{D}_{3}(\hat{\rho}_{0}\hat{\theta}_{R}+\hat{\theta}_{0}\hat{\rho}_{R})\,{\rm d}z\leq \tau\|\frac{\nabla (\rho_{0}\theta_{R}+\theta_{0}\rho_{R})}{\varepsilon}\|_{L^2}^2+C_{\tau,\sigma}\|u_{R}\|_{H^1}^2.
\end{align}
Combining \eqref{5.7-6}--\eqref{5.7-7-1}, one has
\begin{align}\label{5.7-9}
&-\frac{1}{\varepsilon}\int_{\hat{Q}\cap \hat{\Omega}}\chi_{1}^2\operatorname{div}_{z}(\hat{\rho}_{0}\hat{u}_{R})\hat{D}_{33}^2(\hat{\rho}_{0}\hat{\theta}_{R}+\hat{\theta}_{0}\hat{\rho}_{R})\,{\rm d}z-\frac{1}{\varepsilon}\int_{\hat{Q}\cap\hat{\Omega}}\chi_{1}^2\hat{D}_{3}\operatorname{div}_{z}(\hat{\rho}_{0}\hat{u}_{R})\hat{D}_{3}(\hat{\rho}_{0}\hat{\theta}_{R}+\hat{\theta}_{0}\hat{\rho}_{R})\,{\rm d}z\nonumber\\
&\leq \tau\|\frac{\nabla (\rho_{0}\theta_{R}+\theta_{0}\rho_{R})}{\varepsilon}\|_{L^2}^2+C_{\tau,\sigma}\|u_{R}\|_{H^1}^2 +C(\mathcal{I}_{1})\|(\rho_{R},\theta_{R})\|_{H^2}^2+
%(\|u_{1}\|_{H^2}+\varepsilon\|u_{2}\|_{H^2}+\varepsilon\|\tilde{u}_{R}\|)\nonumber\\
+C\|(\rho_{R},\theta_{R})\|_{H^2}\|\mathfrak{R}_{1}\|_{H^1}\nonumber\\
&\quad +\frac{1}{\varepsilon}\sum\limits_{m=1}^2\int_{\hat{Q}\cap \partial\hat{\Omega}}\chi_{1}^2\hat{D}_{m}w\hat{D}_{3}(\hat{\rho}_{0}\hat{u}_{R}^{m})\hat{D}_{3}(\hat{\rho}_{0}\hat{\theta}_{R}+\hat{\theta}_{0}\hat{\rho}_{R})\,{\rm d}z_{1}{\rm d}z_{2}.
%&\quad +C(\|\rho_{R}\|_{H^2}^2+\|\theta_{R}\|_{H^2}^2)(\|u_{1}\|_{H^2}+\varepsilon\|u_{2}\|_{H^2}+\varepsilon\|\tilde{u}_{R}\|)
%&\quad+C(\|\rho_{R}\|_{H^2}+\|\theta_{R}\|_{H^2})\big[(\|\rho_{1}\|_{H^2}+\varepsilon\|\rho_{2}\|_{H^2}+\varepsilon^2\|\rho_{3}\|_{H^2})\|\tilde{u}_{R}\|_{H^2}+\|r_{1}\|_{H^1}\big].
\end{align}
For $J_{3}$ in \eqref{5.7-5}, a direct calculation shows that
\begin{align}\label{5.7-11}
|J_{3}|&\leq \frac{\kappa}{16}\int_{\hat{Q}\cap \hat{\Omega}}\chi_{1}^2\frac{\hat{\rho}}{\hat{\theta_{0}}}|\hat{D}_{33}^2\hat{\theta}_{R}|^2\,{\rm d}z+C\sum\limits_{i=1}^2\|\nabla_{z} \hat{D}_{i}\hat{\theta}_{R}|_{L^2}^2.
\end{align}
For $J_{4}$ and $J_{6}$, by using the H\"{o}lder inequality, one has
\begin{align}\label{5.7-12}
|J_{4}|+|J_{6}|&\leq \frac{\kappa}{16}\int_{\hat{Q}\cap \hat{\Omega}}\chi_{1}^2\frac{\hat{\rho}_{0}}{\hat{\theta}_{0}}|\hat{D}_{33}^2\hat{\theta}_{R}|^2\,{\rm d}z+C\|(u_{R},\theta_{R})\|_{H^1}^2+C\|\rho_{R}\|_{L^2}^2\|u_{1}\|_{H^3}^2.
\end{align}
For $J_{5}$, we integrate by parts to get
\begin{align}\label{5.7-8}
&J_{5}=\frac{1}{\varepsilon}\sum\limits_{m=1}^2\int_{\hat{Q}\cap \partial\hat{\Omega}}\chi_{1}^2\hat{D}_{m}w\hat{D}_{3}(\hat{\rho}_{0}\hat{u}_{R}^{m})\hat{D}_{3}(\hat{\rho}_{0}\hat{\theta}_{R}+\hat{\theta}_{0}\hat{\rho}_{R})\,{\rm d}z_{1}{\rm d}z_{2}\nonumber\\
&\qquad +\frac{1}{\varepsilon}\sum\limits_{m=1}^2\int_{\hat{Q}\cap \hat{\Omega}}2\chi_{1}\hat{D}_{3}\chi_{1}\hat{D}_{m}w\hat{D}_{3}(\hat{\rho}_{0}\hat{u}_{R}^{m})\hat{D}_{3}(\hat{\rho}_{0}\hat{\theta}_{R}+\hat{\theta}_{0}\hat{\rho}_{R})\,{\rm d}z\nonumber\\
&\qquad +\frac{1}{\varepsilon}\sum\limits_{m=1}^2\int_{\hat{Q}\cap \hat{\Omega}}\chi_{1}^2\hat{D}_{m}w\hat{D}_{33}^2(\hat{\rho}_{0}\hat{u}_{R}^{m})\hat{D}_{3}(\hat{\rho}_{0}\hat{\theta}_{R}+\hat{\theta}_{0}\hat{\rho}_{R})\,{\rm d}z,
\end{align}
which, together with the H\"{o}lder inequality and \eqref{5.7-9}, implies that
\begin{align}\label{5.7-13}
&-J_{5}-\frac{1}{\varepsilon}\int_{\hat{Q}\cap \hat{\Omega}}\chi_{1}^2\operatorname{div}_{z}(\hat{\rho}_{0}\hat{u}_{R})\hat{D}_{33}^2(\hat{\rho}_{0}\hat{\theta}_{R}+\hat{\theta}_{0}\hat{\rho}_{R})\,{\rm d}z-\frac{1}{\varepsilon}\int_{\hat{Q}\cap\hat{\Omega}}\chi_{1}^2\hat{D}_{3}\operatorname{div}_{z}(\hat{\rho}_{0}\hat{u}_{R})\hat{D}_{3}(\hat{\rho}_{0}\hat{\theta}_{R}+\hat{\theta}_{0}\hat{\rho}_{R})\,{\rm d}z\nonumber\\
%&\leq %\frac{1}{\varepsilon}\sum\limits_{m=1}^2\int_{\hat{Q}\cap \partial\hat{\Omega}}\chi_{1}^2\hat{D}_{m}w\hat{D}_{3}(\hat{\rho}_{0}\hat{u}_{R}^{m})\hat{D}_{3}(\hat{\rho}_{0}\hat{\theta}_{R}+\hat{\theta}_{0}\hat{\rho}_{R})\,{\rm d}z_{1}{\rm d}z_{2}\nonumber\\
&\leq (\tau+\sigma)\|\frac{\nabla (\rho_{0}\theta_{R}+\theta_{0}\rho_{R})}{\varepsilon}\|_{L^2}^2+C\sigma \sum\limits_{m=1}^2\int_{\hat{Q}\cap \hat{\Omega}}\chi_{1}^2\hat{\rho}_{0}|\hat{D}_{33}^2\hat{u}_{R}^{m}|^2\,{\rm d}z+C_{\tau,\sigma}\|u_{R}\|_{H^1}^2\nonumber\\
&\quad +C(\mathcal{I}_{1})\|(\rho_{R},\theta_{R})\|_{H^2}^2+C\|(\rho_{R},\theta_{R})\|_{H^2}\|\mathfrak{R}_{1}\|_{H^1}\nonumber\\
%\nonumber\\
%&\leq %\frac{1}{\varepsilon}\sum\limits_{m=1}^2\int_{\hat{Q}\cap \partial\hat{\Omega}}\chi_{1}^2\hat{D}_{m}w\hat{D}_{3}(\hat{\rho}_{0}\hat{u}_{R}^{m})\hat{D}_{3}(\hat{\rho}_{0}\hat{\theta}_{R}+\hat{\theta}_{0}\hat{\rho}_{R})\,{\rm d}z_{1}{\rm d}z_{2}\nonumber\\
&\leq (\sigma +\tau)\|\frac{\nabla (\rho_{0}\theta_{R}+\theta_{0}\rho_{R})}{\varepsilon}\|_{L^2}^2+C\sigma\|u_{R}\|_{H^2}^2 +C_{\tau,\sigma}\|u_{R}\|_{H^1}^2+C(\mathcal{I}_{1})\|(\rho_{R},\theta_{R})\|_{H^2}^2\nonumber\\
&\quad+C\|(\rho_{R},\theta_{R})\|_{H^2}\|\mathfrak{R}_{1}\|_{H^1}.
\end{align}
%where we have used \eqref{4.24-0} to obtain
%\begin{align*}
%	C\sigma \sum\limits_{m=1}^2\int_{\hat{Q}\cap \hat{\Omega}}\chi_{1}^2\hat{\rho}_{0}|\hat{D}_{33}^2\hat{u}_{R}^{m}|^2\,{\rm d}z
%	&\leq C(\mathcal{I}_{1})\sigma\|(\rho_{R},\theta_{R})\|_{H^2}^2+C\sigma[\|\mathfrak{R}_{2}\|_{L^2}^2+\|u_{R}\|_{H^1}^2+\|\nabla \operatorname{div}u_{R}\|_{L^2}^2].
%\end{align*}
For $J_{7}$, we get directly from the H\"{o}lder inequality that
\begin{align}\label{5.7-14}
|J_{7}|&\leq \frac{\kappa}{16}\int_{\hat{Q}\cap \hat{\Omega}}\chi_{1}^2\frac{\hat{\rho}_{0}}{\hat{\theta}_{0}}|\hat{D}_{33}^2\hat{\theta}_{R}|^2\,{\rm d}z
%\nonumber\\
%&\quad +C(\|r_{3}\|_{L^2}^2+\|\hat{W}_{3}\|_{L^2}^2+\|\Psi(\nabla (u_{1}+\varepsilon(u_{2}+\tilde{u}_{R})))\|_{L^2}^2+\|G\|_{L^2}^2)+C\|\theta_{R}\|_{H^1}^2\nonumber\\
+C\|\theta_{R}\|_{H^1}^2+C\big[\|\hat{\mathfrak{R}}_{3}\|_{L^2}^2+\|\hat{\tilde{G}}\|_{L^2}^2+\|\hat{W}_{3}\|_{L^2}^2\big]\nonumber\\
%\nonumber\\&\quad +C\|r_{3}\|_{L^2}^2+C\|\theta_{R}\|_{H^1}^2+C\sigma(\|\rho_{R}\|_{H^1}^2+\|u_{R}\|_{H^2}^2+\|\theta_{R}\|_{H^2}^2)\nonumber\\
%&\quad+C(\|u_{1}\|_{H^2}^4+\varepsilon^4(\|u_{2}\|_{H^2}^4+\|u_{R}\|_{H^2}^4))+C\varepsilon^2\|u_{R}\|_{H^2}^2(\|\nabla \theta_{0}\|_{H^1}^2+\|\theta_{1}\|_{H^2}^2+\|\tilde{\theta}_{R}\|_{H^2}^2)\nonumber\\
%&\quad +C\varepsilon^2\|\theta_{R}\|_{H^2}^2(\|u_{1}\|_{H^2}^2+\|u_{2}\|_{H^2}^2+\|\tilde{u}_{R}\|_{H^2}^2)\nonumber\\
&\leq \frac{\kappa}{16}\int_{\hat{Q}\cap \hat{\Omega}}\chi_{1}^2\frac{\hat{\rho}_{0}}{\hat{\theta}_{0}}|\hat{D}_{33}^2\hat{\theta}_{R}|^2\,{\rm d}z+C(\mathcal{I}_{1})\|(\rho_{R},u_{R},\theta_{R})\|_{H^2}^2+C\sigma^2\|(u_{R},\theta_{R})\|_{H^2}^2\nonumber\\
&\quad +C(\mathcal{I}_{0})\|\theta_{R}\|_{H^1}^2+C\|\mathfrak{R}_{3}\|_{L^2}^2,
\end{align}
where we have used \eqref{5.4-1}--\eqref{5.5} in the last inequality.

For $J_{8}$, it follows from the H\"older inequality that
\begin{align}\label{5.7-15}
	|J_{8}|&\leq \frac{\kappa}{16}\int_{\hat{Q}\cap \hat{\Omega}}\frac{\hat{\rho}_{0}}{\theta_{0}}|\hat{D}_{33}^2\hat{\theta}_{R}|^2\,{\rm d}z+C(\mathcal{I}_{0})\|(u_{R},\theta_{R})\|_{H^1}^2+C\|\rho_{R}\|_{L^2}^2\|u_{1}\|_{H^3}^2.
	%+C(\|\theta_{1}\|_{H^3}^2+\|\rho_{1}\|_{H^2}^2)\|u_{R}\|_{L^2}^2\nonumber\\
	%&\quad +C\|u_{1}\|_{H^2}^2(\|u_{R}\|_{L^2}^2+\|\rho_{R}\|_{L^2}^2)+C\|\nabla \rho_{0}\|_{H^2}^2\|\theta_{R}\|_{H^1}^2.
\end{align}
For $J_{9}$ and $J_{10}$, one has from the Sobolev embedding and H\"{o}lder inequality that
\begin{align}\label{5.7-16}
|J_{9}|+|J_{10}|&\leq C(\mathcal{I}_{1})\|\rho_{R}\|_{H^2}^2+C\|\mathfrak{R}_{1}\|_{L^2}\|\rho_{R}\|_{H^2}.
%&\leq C(\|u_{1}\|_{H^2}+\varepsilon(\|u_{2}\|_{H^2}+\|\tilde{u}_{R}\|_{H^2}) \|\rho_{R}\|_{H^2}^2\nonumber\\
%&\quad +C\|\rho_{R}\|_{H^2}\|u_R\|_{H^1}(\|\rho_{1}\|_{H^3}+\varepsilon\|\rho_{2}\|_{H^3}+\varepsilon^2\|\rho_{3}\|_{H^3})+C\|\mathfrak{R}_{1}\|_{H^1}\|\rho_{R}\|_{H^2}.
\end{align}
Combining \eqref{5.7-1}--\eqref{5.7-16},
%and noting the smallness of $\sigma$
%and the singular terms in \eqref{5.7-4}, \eqref{5.7-9} and \eqref{5.7-13} cancel each other out,
we obtain \eqref{5.7-17}. Therefore the proof of Lemma \ref{lem5.0-1} is complete. $\hfill\square$
\smallskip
%Combining \eqref{4.24-0}, Lemmas \ref{lem4.1}--\ref{lem4.2} and \ref{lem5.0} and using the smallness of $\sigma$,
%\begin{corollary}\label{cor5.1}
	%Lemmas \ref{lem4.1}, \ref{lem4.2} and \ref{lem5.0}, we conclude
	%\begin{align*}
	%	\|u_{R}\|_{H^2}+\|\theta_{R}\|_{H^2}+\|\frac{\nabla(\rho_{R}\theta_{0}+\rho_{R}\theta_{0})}{\varepsilon}\|_{L^2}^2\leq
	%\end{align*}
%where $C_{11}$ is a positive constant, which may depend on $\tau$ but is independent of $\varepsilon$.
%\end{corollary}

\smallskip

\noindent\textbf{Step 3. Estimates of $\varepsilon\|\hat{D}_{j3}^2(\operatorname{div}_{z}\hat{u}_{R})\|_{L^2}\,\,(j=1,2)$.}
\begin{lemma}\label{lem5.1}
	There exists a positive constant $C_{\sigma}$, which depends on $\sigma$ but is independent of $\varepsilon$, such that
	\begin{align}\label{5.9-9}
		&(\mu+\zeta)\varepsilon^2\sum\limits_{j=1}^2\int_{\hat{\Omega}\cap\hat{Q}}\chi_{1}^2\hat{\rho}_{0}| \hat{D}_{j3}^2(\operatorname{div}_{z}\hat{u}_{R})|^2\,{\rm d}z+\kappa\varepsilon^2\sum\limits_{j=1}^2\int_{\hat{Q}\cap \hat{\Omega}}\chi_{1}^2\frac{\hat{\rho}_{0}}{\hat{\theta}_{0}}|\hat{D}_{j33}^3\hat{\theta}_{R}|^2\,{\rm d}z\nonumber\\
		&\leq C\big[\|(u_{R},\theta_{R})\|_{H^2}^2+\v^2\|(\mathfrak{R}_{2},\mathfrak{R}_{3})\|_{H^1}^2\big]+C\varepsilon^2\sum\limits_{j,m=1}^2\big(|\nabla_{z} \hat{D}_{jm}^2\hat{u}_{R}\|_{L^2}^2+\|\hat{D}_{jmm}^3\hat{\theta}_{R}\|_{L^2}^2\big)\nonumber\\
		&\quad +[C(\mathcal{I}_{1})C_{\sigma}\v^2+C\sigma]\|\rho_{R}\|_{H^2}^2+C\sigma^2 \v^2\|(u_{R},\theta_{R})\|_{H^3}^2+C\v^2\|\rho_{R}\|_{H^2}\|\mathfrak{R}_{1}\|_{H^2},
	\end{align}
where $C>0$ is independent of $\sigma$ and $\v$.
\end{lemma}

\noindent\textbf{Proof}. We apply $\hat{D}_{j}(j=1,2)$ to \eqref{5.7} and multiply the resultant equation by $-\varepsilon^2\chi_{1}^2\hat{D}_{j3}^2\operatorname{div}_{z}(\hat{\rho}_{0}\hat{u}_{R})$ to get
\begin{align}\label{5.9}
	&(\mu+\zeta)\varepsilon^2\int_{\hat{\Omega}\cap \hat{Q}}\hat{\rho}_{0}\chi_{1}^2|\hat{D}_{j 3}^2(\operatorname{div}_{z}\hat{u}_{R})|^2\,{\rm d}z\nonumber\\
	&=\varepsilon\int_{\hat{\Omega}\cap \hat{Q}}\chi_{1}^2\hat{D}_{j 3}^2(\hat{\rho}_{0}\hat{\theta}_{R})\hat{D}_{j 3}^2\operatorname{div}_{z}(\hat{\rho}_{0}\hat{u}_{R})\,{\rm d}z+\varepsilon\int_{\hat{\Omega}\cap \hat{Q}}\chi_{1}^2\hat{D}_{j3}^2(\hat{\rho}_{R}\hat{\theta}_{0})\hat{D}_{j 3}^2\operatorname{div}_{z}(\hat{\rho}_{0}\hat{u}_{R})\,{\rm d}z\nonumber\\
	&\quad -(\mu+\zeta)\varepsilon^2\int_{\hat{\Omega}\cap \hat{Q}}\chi_{1}^2\hat{D}_{j 3}^2(\operatorname{div}_{z}\hat{u}_{R})[\hat{D}_{j3}^2\operatorname{div}_{z}(\hat{\rho}_{0}\hat{u}_{R})-\hat{\rho}_{0}\hat{D}_{j3}^2\operatorname{div}_{z}\hat{u}_{R}]\,{\rm d}z\nonumber\\
	&\quad +\varepsilon^2\int_{\hat{Q}\cap \hat{\Omega}}\chi_{1}^2\hat{D}_{j3}^2(\hat{\rho}_{R}\hat{\theta}_{1}+\hat{\rho}_{1}\hat{\theta}_{R})\hat{D}_{j3}^2\operatorname{div}_{z}(\hat{\rho}_{0}\hat{u}_{R})\,{\rm d}z\nonumber\\
	%+\varepsilon^2\int_{\hat{Q}\cap \hat{\Omega}}\chi_{1}^2\hat{D}_{j3}^2(\hat{\rho}_{1}\hat{\theta}_{R})\hat{D}_{j3}^2\operatorname{div}_{z}(\hat{\rho}_{0}\hat{u}_{R})\,{\rm d}z\nonumber\\
	&\quad -\mu\varepsilon^2\int_{\hat{Q}\cap \hat{\Omega}}\chi_{1}^2\hat{D}_{j}[(\Delta_{z}\hat{u}_{R}-\nabla_{z}\operatorname{div}_{z}\hat{u}_{R})\cdot \hat{e}_{3}]\hat{D}_{j3}^2\operatorname{div}_{z}(\hat{\rho}_{0}\hat{u}_{R}){\rm d}z\nonumber\\
	&\quad +\varepsilon^2\int_{\hat{\Omega}\cap \hat{Q}}\chi_{1}^2\hat{D}_{j}[\hat{\rho}_{0}(\hat{u}_{1}\cdot \nabla_{z} \hat{u}_{R}+\hat{u}_{R}\cdot \nabla_{z} \hat{u}_{1})\cdot \hat{e}_{3}]\hat{D}_{j 3}^2(\operatorname{div}_{z}(\hat{\rho}_{0}\hat{u}_{R}))\,{\rm d}z\nonumber\\
	&\quad +\varepsilon^2\int_{\hat{Q}\cap \hat{\Omega}}\chi_{1}^2\hat{D}_{j}(\hat{W}_{2}^3+\hat{F}^{\varepsilon,3}+\hat{\mathfrak{R}}_{2}^3)\hat{D}_{j3}^2\operatorname{div}_{z}(\hat{\rho}_{0}\hat{u}_{R})\,{\rm d}z=:\sum\limits_{m=1}^7J_{m}.
\end{align}
It is clear that
%For $J_{1}$, by using the H\"{o}lder equality, one has
\begin{align}\label{5.9-1}
|J_{1}|\leq \frac{(\mu+\zeta)\varepsilon^2}{16}\int_{\hat{Q}\cap \hat{\Omega}}\chi_{1}^2\hat{\rho}_{0}|\hat{D}_{j3}^2\operatorname{div}_{z}\hat{u}_{R}|^2\,{\rm d}z+C\varepsilon^2\|u_{R}\|_{H^2}^2+C\|\theta_{R}\|_{H^2}^2,
\end{align}
and
%For $J_{3}$, a direct calculation shows that
\begin{align}\label{5.9-4}
|J_{3}|\leq \frac{(\mu+\zeta)\varepsilon^2}{16}\int_{\hat{Q}\cap \hat{\Omega}}\chi_{1}^2\hat{\rho}_{0}|\hat{D}_{j3}^2\operatorname{div}_{z}\hat{u}_{R}|^2\,{\rm d}z+C\varepsilon^2\|u_{R}\|_{H^2}^2.
\end{align}
Similarly, for $J_{m}(m=4,\cdots,7)$, we get from the H\"{o}lder inequality that
\begin{align}\label{5.9-5}
&\sum\limits_{m=4}^{7}|J_{m}|\leq \frac{(\mu+\zeta)\varepsilon^2}{16}\int_{\hat{Q}\cap \hat{\Omega}}\chi_{1}^2\hat{\rho}_{0}|\hat{D}_{j3}^2\operatorname{div}_{z}\hat{u}_{R}|^2\,{\rm d}z+C(\mathcal{I}_{0})\varepsilon^2\|(\rho_{R},u_{R},\theta_{R})\|_{H^2}^2\nonumber\\
%\|u_{R}\|_{H^2}^2+C\varepsilon^2\|\rho_{1}\|_{H^2}^2\|\theta_{R}\|_{H^2}^2+\|\rho_{R}\|_{H^2}^2\|\theta_{1}\|_{H^2}^2+\|u_{1}\|_{H^2}^2\|u_{R}\|_{H^2}^2\nonumber\\
&\qquad\qquad\quad  +C\varepsilon^2\Big[\sum\limits_{m=1}^2\|\nabla_{z} \hat{D}_{jm}^2\hat{u}_{R}\|_{L^2}^2+\|W_{2}\|_{H^1}^2+\|F^{\varepsilon}\|_{H^1}^2+\|\mathfrak{R}_{2}\|_{H^1}^2\Big]\nonumber\\
&\leq \frac{(\mu+\zeta)\varepsilon^2}{16}\int_{\hat{Q}\cap \hat{\Omega}}\chi_{1}^2\hat{\rho}_{0}|\hat{D}_{j3}^2\operatorname{div}_{z}(\hat{\rho}_{0}\hat{u}_{R})\,{\rm d}z+C\v^2\sum\limits_{m=1}^2\|\nabla_{z} \hat{D}_{jm}^2\hat{u}_{R}\|_{L^2}^2+C\v^2\|\mathfrak{R}_{2}\|_{H^1}^2\nonumber\\
&\quad +C(\mathcal{I}_{1})\v^2\|(\rho_{R},u_{R},\theta_{R})\|_{H^2}^2+C\v^2\sigma^2\{\|\rho_{R}\|_{H^2}^2+\|(u_{R},\theta_{R})\|_{H^3}^2\}.
%&\quad+C\varepsilon^2\Big[\|\rho_{R}\|_{H^2}^2\|\theta_{1}\|_{H^2}^2+\|\rho_{1}\|_{H^2}^2\|\theta_{R}\|_{H^2}^2+\sum\limits_{m=1}^2\|\nabla \hat{D}_{jm}^2\hat{u}_{R}\|_{L^2}^2+\|u_{R}\|_{H^2}^2+\|r_{2}\|_{H^1}^2\Big]\nonumber\\
%&\quad +C\varepsilon^2\sigma(\|\rho_{R}\|_{H^2}^2+\|u_{R}\|_{H^3}^2+\|\theta_{R}\|_{H^3}^2)+C\varepsilon^4\|u_{R}\|_{H^2}^2(\|u_{1}\|_{H^2}^2+\|u_{2}\|_{H^2}^2+\|\tilde{u}_{R}\|_{H^2}^2)\nonumber\\
%&\quad+C\varepsilon^4\Big[\|\theta_{R}\|_{H^2}^2(\|\rho_{2}\|_{H^2}^2+\|\rho_{3}\|_{H^2}^2)+\|\rho_{R}\|_{H^2}^2\big[(\|u_{1}\|_{H^2}^2+\|\| u_{2}\|_{H^2}^2+\|\tilde{u}_{R}\|_{H^2}^2)^2+\|\tilde{\theta}_{R}\|_{H^2}^2\big]\Big],
\end{align}
where we have used \eqref{5.4-1}--\eqref{5.5}.

To deal with $J_{2}$, we apply $\hat{D}_{j3}^2$ to $\eqref{5.2}_{1}$ and multiply the resulting equation by $\varepsilon^2\chi_{1}^2\hat{D}_{j3}^2(\hat{\theta}_{0}\hat{\rho}_{R})$ to get
\begin{align}\label{5.9-2}
J_{2}&=\varepsilon\sum\limits_{m=1}^2\int_{\hat{Q}\cap \hat{\Omega}}\chi_{1}^2\hat{D}_{j3}^2(\hat{\theta}_{0}\hat{\rho}_{R})\hat{D}_{j3}^2(\hat{D}_{m}w\hat{D}_3(\hat{\rho}_{0}\hat{u}_{R}^{m}))\,{\rm d}z\nonumber\\
&\quad-\varepsilon^2\int_{\hat{Q}\cap \hat{\Omega}}\chi_{1}^2\hat{D}_{j3}^2(\hat{\theta}_{0}\hat{\rho}_{R})\hat{D}_{j3}^2\big\{\operatorname{div}_{z}[\hat{\rho}_{R}(\hat{u}_{1}+\varepsilon(\hat{u}_{2}+\hat{\tilde{u}}_{R}))]+\hat{W}_{1}\big\}\,{\rm d}z\nonumber\\
&\quad -\v^2\int_{\hat{Q}\cap \hat{\Omega}}\chi_{1}^2\hat{D}_{j3}^2(\hat{\theta}_{0}\hat{\rho}_{R})\hat{D}_{j3}^2\hat{\mathfrak{R}}_{1}\,{\rm d}z=:\sum\limits_{j=1}^{3}J_{2,j}.
%&\quad-\varepsilon^2\int_{\hat{Q}\cap \hat{\Omega}}\chi_{1}^2\hat{D}_{j3}^2(\hat{\theta}_{0}\hat{\rho}_{R})\hat{D}_{j3}^2\operatorname{div}_{z}[\hat{\tilde{u}}_{R}(\hat{\rho}_{1}+\varepsilon\hat{\rho}_{2}+\varepsilon^2\hat{\rho}_{3})]\,{\rm d}z\nonumber\\
%&\quad -\varepsilon^2\int_{\hat{Q}\cap \hat{\Omega}}\chi_{1}^2\hat{D}_{j3}^2(\hat{\theta}_{0}\hat{\rho}_{R})\hat{D}_{j3}^2\hat{r}_{1}\,{\rm d}z=:\sum\limits_{j=1}^{4}J_{2,j}.
\end{align}
For $J_{2,1}$ and $J_{2,3}$, a direct calculations shows that
\begin{align}\label{5.9-6}
|J_{2,1}|+|J_{2,3}|&\leq \frac{(\mu+\zeta)\varepsilon^2}{16}\int_{\hat{Q}\cap \hat{\Omega}}\chi_{1}^2\hat{\rho}_{0}|\hat{D}_{j33}^2\hat{u}_{R}|^2\,{\rm d}z+C\sigma\varepsilon^2\|u_{R}\|_{H^2}^2\nonumber\\
&\quad +C\sigma\|\rho_{R}\|_{H^2}^2+C\varepsilon^2\|\mathfrak{R}_{1}\|_{H^2}\|\rho_{R}\|_{H^2}.
\end{align}
For $J_{2,2}$, by similar arguments as in \eqref{5.6-5}--\eqref{5.6-6}, we have
\begin{align}\label{5.9-7}
|J_{2,2}|&\leq C(\mathcal{I}_{1})C_{\sigma}\v^2\|\rho_{R}\|_{H^2}^2.
%&\leq C\varepsilon^2\big[\|\rho_{R}\|_{H^2}^2\big(\|u_{1}\|_{H^3}+\varepsilon(\|u_{2}\|_{H^3}+\|\tilde{u}_{R}\|_{H^3})\big)\nonumber\\
%&\qquad\quad +\|\rho_{R}\|_{H^2}\|\tilde{u}_{R}\|_{H^3}(\|\rho_{1}\|_{H^3}+\varepsilon\|\rho_{2}\|_{H^3}+\varepsilon^2\|\rho_{3}\|_{H^3})\big].
\end{align}

We apply $\hat{D}_{j}$ to $\eqref{5.2}_{3}$ and multiply the resultant equation by $\varepsilon^2\chi_{1}^2\hat{D}_{j33}^3(\hat{\rho}_{0}\hat{\theta}_{R})$ to obtain
\begin{align}\label{5.9-3}
&\kappa\varepsilon^2\int_{\hat{Q}\cap \hat{\Omega}}\chi_{1}^2\frac{\hat{\rho}_{0}}{\hat{\theta}_{0}}|\hat{D}_{j33}^3\hat{\theta}_{R}|^2\,{\rm d}z\nonumber\\
&=\varepsilon\int_{\hat{Q}\cap \hat{\Omega}}\chi_{1}^2\hat{D}_{j33}^3(\hat{\rho}_{0}\hat{\theta}_{R})\hat{D}_{j}\operatorname{div}_{z}(\hat{\rho}_{0}\hat{u}_{R})\,{\rm d}z -\varepsilon\sum\limits_{m=1}^2\int_{\hat{Q}\cap \hat{\Omega}}\chi_{1}^2\hat{D}_{j33}^3(\hat{\rho}_{0}\hat{\theta}_{R})\hat{D}_{j}(\hat{D}_{m}w\hat{D}_{3}(\hat{\rho}_{0}\hat{u}_{R}^{m}))\,{\rm d}z\nonumber\\
&\quad -\kappa\varepsilon^2\int_{\hat{Q}\cap \hat{\Omega}}\chi_{1}^2\hat{D}_{j33}^3(\hat{\rho}_{0}\hat{\theta}_{R})[\hat{D}_{j}(\frac{\Delta_{z}\hat{\theta}_{R}}{\hat{\theta}_{0}})-\frac{1}{\hat{\theta}}_{0}\hat{D}_{j33}^3\hat{\theta}_{R}]\,{\rm d}z\nonumber\\
&\quad +\varepsilon^2P_{1}\int_{\hat{Q}\cap \hat{\Omega}}\chi_{1}^2\hat{D}_{j33}^3(\hat{\rho}_{0}\hat{\theta}_{R})\hat{D}_{j}(\frac{1}{\hat{\theta}_0}\operatorname{div}_{z}\hat{u}_{R})\,{\rm d}z\nonumber\\
&\quad +\varepsilon^2\int_{\hat{Q}\cap \hat{\Omega}}\chi_{1}^2\hat{D}_{j33}^3(\hat{\rho}_{0}\hat{\theta}_{R})\hat{D}_{j}\big(\frac{1}{\hat{\theta}_{0}}(\hat{\rho}_{0}\hat{\theta}_{R}+\hat{\rho}_{R}\hat{\theta}_{0})\operatorname{div}_{z}\hat{u}_{1}\big)\,{\rm d}z\nonumber\\
&\quad +\varepsilon^2\int_{\hat{Q}\cap \hat{\Omega}}\chi_{1}^2\hat{D}_{j33}^3(\hat{\rho}_{0}\hat{\theta}_{R})\hat{D}_{j}\big(\frac{1}{\hat{\theta}_{0}}\hat{G}^{\varepsilon}+\hat{\mathfrak{R}}_{3}+\hat{W}_{3}\big)\,{\rm d}z\nonumber\\
%&\quad +\varepsilon^2\int_{\hat{Q}\cap \hat{\Omega}}\chi_{1}^2\hat{D}_{j33}^3(\hat{\rho}_{0}\hat{\theta}_{R})\hat{D}_{j}(\frac{1}{\hat{\theta}_{0}}\hat{\Psi}(\nabla (u_{1}+\varepsilon(\hat{u}_{2}+\tilde{u}_{R}))))\,{\rm d}z\nonumber\\
&\quad +\varepsilon^2\int_{\hat{Q}\cap \hat{\Omega}}\chi_{1}^2\hat{D}_{j33}^3\big(\hat{\rho}_{0}\hat{\theta}_{R})\hat{D}_{j}\Big\{\frac{1}{\hat{\theta}_{0}}\{\hat{\rho}_{0}(\hat{u}_{R}\cdot \nabla_{z}\hat{\theta}_{1}+\hat{u}_{1}\cdot \nabla_{z}\hat{\theta}_{R})\nonumber\\
&\qquad \qquad\qquad +\hat{\rho}_{1}\hat{u}_{R}\cdot \nabla_{z}\hat{\theta}_{0}+\hat{\rho}_{R}\hat{u}_{1}\cdot \nabla_{z}\hat{\theta}_{0}\}\Big\}\,{\rm d}z=:\sum\limits_{m=8}^{14}J_{m}.
\end{align}
Applying the H\"{o}lder inequality and Sobolev embedding, we have
\begin{align}\label{5.9-8}
&\sum\limits_{m=8}^{14}J_{m}\leq \frac{\kappa\varepsilon^2}{2}\int_{\hat{Q}\cap \hat{\Omega}}\chi_{1}^2\frac{\hat{\rho}_{0}}{\hat{\theta}_{0}}|\hat{D}_{j33}^3\hat{\theta}_{R}|^2\,{\rm d}z+C\varepsilon^2\sum\limits_{l=1}^2\|\hat{D}_{jll}^3\theta_{R}\|_{L^2}^2+C\|u_{R}\|_{H^2}^2\nonumber\\
&\qquad\qquad\quad +C\varepsilon^2\big[\|\theta_{R}\|_{H^2}^2+\|\rho_{R}\|_{H^1}^2\|u_{1}\|_{H^3}^2+\|\hat{G}^{\varepsilon}\|_{H^1}^2+\|\hat{W}_{3}\|_{H^1}^2+\|\hat{\mathfrak{R}}_{3}\|_{H^1}^2\big]\nonumber\\%+\|r_{3}\|_{H^1}^2+\|\Psi(\nabla (u_{1}+\varepsilon(u_{2}+\tilde{u}_{R})))\|_{H^1}^2\big]\nonumber\\
&\qquad\qquad\quad\leq \frac{\kappa\varepsilon^2}{2}\int_{\hat{Q}\cap \hat{\Omega}}\chi_{1}^2\frac{\hat{\rho}_{0}}{\hat{\theta}_{0}}|\hat{D}_{j33}^3\hat{\theta}_{R}|^2\,{\rm d}z+C\sum\limits_{l=1}^{2}\varepsilon^2\|\hat{D}_{jll}^3\theta_{R}\|_{L^2}^2+C\v^2\sigma^2\|(u_{R},\theta_{R})\|_{H^3}^2\nonumber\\
&\qquad\qquad\quad\,\,\,\,+C(\mathcal{I}_{1})\v^2\|\rho_{R}\|_{H^2}^2+C\|u_{R}\|_{H^2}^2+C\v^2\|\theta_{R}\|_{H^2}^2+C\v^2\|\mathfrak{R}_{3}\|_{H^1}^2 ,
%&\quad +C\varepsilon^2 (\|u_{1}\|_{H^3}^4+\varepsilon^4(\|u_{2}\|_{H^3}^4+\|\tilde{u}_{R}\|_{H^3}^4))\nonumber\\
%&\quad +C\varepsilon^4\|u_{R}\|_{H^2}^2(\|\nabla \theta_{0}\|_{H^1}+\|\theta_{1}\|_{H^2}+\|\tilde{\theta}_{R}\|_{H^2})^2+C\varepsilon^4\|\theta_{R}\|_{H^2}^2(\|u_{1}\|_{H^2}+\|u_{2}\|_{H^2}+\|\tilde{u}_{R}\|_{H^2})^2\nonumber\\
%&\quad +C\varepsilon^4\|\rho_{R}\|_{H^2}^2(\|u_{1}\|_{H^2}^2+\|u_{2}\|_{H^2}^2+\|\tilde{u}_{R}\|_{H^2}^2)(\|\nabla \theta_{0}\|_{H^1}^2+\|\theta_{1}\|_{H^2}^2+\|\tilde{\theta}_{R}\|_{H^2}^2).
\end{align}
where we have used \eqref{5.4-1}--\eqref{5.5}. Combining \eqref{5.9}--\eqref{5.9-8}, we obtain \eqref{5.9-9}. This completes the proof of Lemma \ref{lem5.1}. $\hfill\square$

%\noindent\textbf{Proof. } We apply $\hat{D}_{\tau} (\tau=1,2)$ to \eqref{5.7}, then multiply $-\chi_{1}^2\hat{D}_{\tau3}(\operatorname{div}(\hat{\rho}_{0}\hat{u}_{R}))$ to get

%On the other hand, we apply $\hat{D}_{\tau 3}$ to $\eqref{5.2}_{1}$ and $\eqref{5.2}_{3}$, and the multiplying the resultant equations by $\chi_{1}^2\hat{D}_{\tau 3}^2(\hat{\rho}_{0}\hat{\theta}_{R})$ and $\chi_{1}^2\hat{D}_{\tau 3}^2(\hat{\theta}_{0}\hat{\rho}_{R})$, and finally adding together, one has
%\begin{align}\label{5.10}
%\kappa\int_{\hat{\Omega}\cap \hat{Q}}\frac{\rho_{0}}{\theta_{0}}\chi_{1}^2|\nabla(\hat{D}_{\tau 3}^2(\hat{\theta}_{R}))|^2\,{\rm d}z+\frac{1}{\varepsilon}\int_{\hat{\Omega}\cap \hat{Q}}\chi_{1}^2\hat{D}_{\tau 3}^2(\operatorname{div}_{z}(\hat{\rho}_{0}\hat{u}_{R}))\hat{D}_{\tau3}^2(\hat{\rho}_{0}\hat{\theta}_{R}+\hat{\theta}_{0}\hat{\rho}_{R})\,{\rm d}z
%\end{align}
%Summing \eqref{5.9} and \eqref{5.10}, then we can get
%\begin{align}\label{5.11}
%&(\mu+\zeta)\int_{\hat{\Omega}\cap \hat{Q}}\hat{\rho}_{0}\chi_{1}^2|\hat{D}_{\tau 3}(\operatorname{div}_{z}\hat{u}_{R})|^2\,{\rm d}z+\kappa\int_{\hat{\Omega}\cap \hat{Q}}\frac{\rho_{0}}{\theta_{0}}\chi_{1}^2|\nabla(\hat{D}_{\tau 3}^2(\hat{\theta}_{R}))|^2\,{\rm d}z
%\end{align}

\medskip

\noindent\textbf{Step 4. Estimate of $\varepsilon\|\hat{D}_{j33}^3\hat{u}_{R}\|_{L^2}\,\,(j=1,2)$.}
\begin{lemma}\label{lem5.2}
	There exists a positive constant $C_{\sigma}$, which depends on $\sigma$ but is independent of $\varepsilon$, such that
	\begin{align}\label{5.16-0}
		\varepsilon\sum\limits_{j=1}^2\|\nabla_{z}^2\hat{D}_{j}\hat{u}_{R}\|_{L^2}
		&\leq
		%+\varepsilon\|\rho_{R}\|_{H^2}\|\theta_{1}\|_{H^2}+\varepsilon\|\rho_{1}\|_{H^2}\|\theta_{R}\|_{H^2}
		C\varepsilon\sum\limits_{j=1}^2\|\nabla_{z} \hat{D}_{j}\operatorname{div}_{z}u_{R}\|_{L^2}+C_{\sigma}\v\|u_{R}\|_{H^2}+C(\mathcal{I}_{1})\v\|(\rho_{R},u_{R},\theta_{R})\|_{H^2}\nonumber\\
		&\quad+C_{\sigma}\|(\rho_{R},\theta_{R})\|_{H^1}+C\sigma\v\|(u_{R},\theta_{R})\|_{H^3}+C\varepsilon\|\mathfrak{R}_{2}\|_{H^1},
		%\nonumber\\
		&%\quad +C_{12}\varepsilon^2\|u_{R}\|_{H^2}\big(\|u_{1}\|_{H^2}+\|u_{2}\|_{H^2}+\|\tilde{u}_{R}\|_{H^2}\big)+C_{12}\varepsilon^2\|\theta_{R}\|_{H^2}\big(\|\rho_{2}\|_{H^2}+\|\rho_{3}\|_{H^2}\big)\nonumber\\
		%&\quad +C_{12}\varepsilon^2\|\rho_{R}\|_{H^2}\big((\|u_{1}\|_{H^2}^2+\|\varepsilon u_2\|_{H^2}^2+\|\varepsilon\tilde{u}_{R}\|_{H^2}^2)+\|\tilde{\theta}_{R}\|_{H^2}\big).
	\end{align}
where $C>0$ is independent of $\v$ and $\sigma$.
\end{lemma}

\noindent\textbf{Proof}. %To control the term $\hat{D}_{j33}^3\hat{u}_{R}(j=1,2)$,
Recalling $\eqref{5.2}_{2}$, we turn to following auxiliary Stokes problem
\begin{equation}\label{5.16}
	\left\{\begin{aligned}
		&-\varepsilon\mu\Delta_{z}(\chi_{1}\hat{D}_{j}\hat{u}_{R})+\nabla_{z}[\chi_{1}\hat{D}_{j}(\hat{\rho}_{0}\hat{\theta}_{R}+\hat{\theta}_{0}\hat{\rho}_{R})]=\varepsilon S_{1},\\
		&\varepsilon\operatorname{div}(\chi_{1}\hat{D}_{j}\hat{u}_{R})=\varepsilon S_{2},\\
		&\chi_{1}\hat{D}_{j}\hat{u}_{R}\vert_{\partial\hat{\Omega}\cap \hat{Q}}=0,
	\end{aligned}
	\right.
\end{equation}
where
\begin{align}
	&S_{1}=-\chi_{1}\hat{D}_{j}\big[\hat{\rho}_{0}[\hat{u}_{1}\cdot \nabla_{z}\hat{u}_{R}+\hat{u}_{R}\cdot \nabla_{z}\hat{u}_{1}]\big]-\chi_{1}\nabla_{z}(\hat{D}_{j}(\hat{\rho}_{R}\hat{\theta}_1)+\hat{D}_{3}(\hat{\rho}_{1}\hat{\theta}_{R}))\nonumber\\
	&\qquad -\frac{1}{\varepsilon}\chi_{1}\hat{D}_{j}[(\hat{D}_{1}w\hat{D}_{3}(\hat{\rho}_{0}\hat{\theta}_{R}+\hat{\theta}_{0}\hat{\rho}_{R}),\hat{D}_{2}w\hat{D}_{3}(\hat{\rho}_{0}\hat{\theta}_{R}+\hat{\theta}_{0}\hat{\rho}_{R}),0)]^{t}\nonumber\\
	&\qquad -\chi_{1}\hat{D}_{j}\hat{W}_{2}-\chi_{1}\hat{D}_{j}\hat{\tilde{F}}^{\varepsilon}-\chi_{1}\hat{D}_{j}\hat{\mathfrak{R}}_{2}-\frac{1}{\varepsilon}\hat{D}_{j}(\hat{\rho}_{0}\hat{\theta}_{R}+\hat{\rho}_{R}\hat{\theta}_{0})\nabla_{z}\chi_{1}\nonumber\\
	&\qquad+\zeta\chi_{1}\nabla_{z}\operatorname{div}_{z}(\hat{D}_{j}\hat{u}_{R})-2\mu\nabla_{z}\chi_{1}\cdot \nabla_{z}(\hat{D}_{j}\hat{u}_{R})^{t}-\mu\hat{D}_{j}\hat{u}_{R}\Delta_{z}\chi_{1},\label{5.16-1}\\
	&S_{2}=\chi_{1}\operatorname{div}_{z}(\hat{D}_{j}\hat{u}_{R})+\hat{D}_{j}\hat{u}_{R}\cdot \nabla_{z}\chi_{1}.\label{5.16-2}
\end{align}

By direct calculations, we can obtain
\begin{align}\label{5.16-3}
\varepsilon\|S_{1}\|_{L^2}&\leq %\varepsilon\|u_{1}\|_{H^2}\|u_{R}\|_{H^2}+\varepsilon\|\rho_{R}\|_{H^2}\|\theta_{1}\|_{H^2}+\varepsilon\|\rho_{1}\|_{H^2}\|\theta_{R}\|_{H^2}
C(\mathcal{I}_{1})\v \|(\rho_{R},u_{R},\theta_{R})\|_{H^2}+C\sigma\|\chi_{1}\nabla_{z}\hat{D}_{j}(\hat{\rho}_{0}\hat{\theta}_{R}+\hat{\theta}_{0}\hat{\rho}_{R})\|_{L^2}+C\sigma\v\|(u_{R},\theta_{R})\|_{H^3}
%+C\varepsilon\|\nabla_{z}^2\operatorname{div}_{z}u_{R}\|_{L^2}+C\varepsilon\sigma\|\hat{D}_{j33}^3u_{R}\|_{L^2}
\nonumber\\
&\quad+C\varepsilon\|\nabla_{z}\hat{D}_{j}\operatorname{div}_{z}u_{R}\|_{L^2} +C_{\sigma}\varepsilon\|u_{R}\|_{H^2}+C\varepsilon\|\mathfrak{R}_{2}\|_{H^1}+C_{\sigma}\|(\rho_{R},\theta_{R})\|_{H^1},
\end{align}
and
\begin{align}\label{5.16-4}
\varepsilon\|S_{2}\|_{H^1}\leq C\varepsilon\|\nabla_{z} \hat{D}_{j}\operatorname{div}_{z}u_{R}\|_{L^2}+C_{\sigma}\varepsilon\|u_{R}\|_{H^2},
\end{align}
where we have used \eqref{5.4} and \eqref{5.5} in \eqref{5.16-3}.

Then applying the classical Stokes's estimates \cite[Theorem IV.5.8]{Boyer-Fabrie}, we have
\begin{align}\label{5.17}
	&\varepsilon\|\nabla_{z}^2(\hat{D}_{j}\hat{u}_{R})\|_{L^2}+\|\nabla_{z} (\chi_{1}\hat{D}_{j}(\hat{\rho}_{0}\hat{\theta}_{R}+\hat{\theta}_{0}\hat{\rho}_{R}))\|_{L^2}\nonumber\\
	&\leq C\v\|S_{1}\|_{L^2}+C\v\|S_{2}\|_{H^{1}}\nonumber\\
	&\leq C(\mathcal{I}_{1})\v\|(\rho_{R},u_{R},\theta_{R})\|_{H^2}+C\varepsilon\|\nabla_{z} \hat{D}_{j}\operatorname{div}_{z}u_{R}\|_{L^2}+C_{\sigma}\varepsilon\|u_{R}\|_{H^2}\nonumber\\
	&\quad+C\varepsilon\|\mathfrak{R}_{2}\|_{H^1} +C_{\sigma}\|(\rho_{R},\theta_{R})\|_{H^1}+C\sigma\v\|(u_{R},\theta_{R})\|_{H^3}
	%\nonumber\\
	%&\leq C\Big[\varepsilon\|u_{R}\|_{H^2}+\varepsilon\|\rho_{R}\|_{H^2}\|\theta_{1}\|_{H^2}+\varepsilon\|\nabla_{z} \hat{D}_{j}\operatorname{div}_{z}u_{R}\|_{L^2}\nonumber\\
	%&\quad+\varepsilon\sum\limits_{j,m=1}^2\| \hat{D}_{jm3}^2u_{R}\|_{L^2}+\varepsilon\|r_{2}\|_{H^1} +\|\rho_{R}\|_{H^1}+\|\theta_{R}\|_{H^1}\Big]\nonumber\\
	%&\quad +C\varepsilon^2\|u_{R}\|_{H^2}(\|u_{1}\|_{H^2}+\|u_{2}\|_{H^2}+\|\tilde{u}_{R}\|_{H^2})+C\varepsilon^2\|\theta_{R}\|_{H^2}(\|\rho_{2}\|_{H^2}+\varepsilon\|\rho_{3}\|_{H^2})\nonumber\\
	%&\quad +C\varepsilon^2\|\rho_{R}\|_{H^2}((\|u_{1}\|_{H^2}^2+\|\varepsilon u_2\|_{H^2}^2+\|\varepsilon\tilde{u}_{R}\|_{H^2}^2)+\|\tilde{\theta}_{R}\|_{H^2}).
\end{align}
where we have used $\sigma\ll 1$. Therefore the proof of Lemma \ref{lem5.2} is complete. $\hfill\square$
\smallskip
%where we point that the estimate $\nabla_{z} (\hat{D}_{j}\operatorname{div}\hat{u}_{R})$ depends only on the estimates established in Steps 1--2, which doesn't depend on the estimate of $\hat{D}_{33}\operatorname{div}_{\hat{u}_{R}}$.

\smallskip

\noindent\textbf{Step 5. Estimate of $\varepsilon\|\hat{D}_{33}^2(\operatorname{div}_{z}\hat{u}_{R})\|_{L^2}$.}
\begin{lemma}\label{lem5.3}
	There exists a positive constant $C_{\sigma}$, which depends on $\sigma$ but is independent of $\v$, such that
	\begin{align}\label{5.19}
	&(\mu+\zeta)\varepsilon^2\int_{\hat{\Omega}\cap\hat{Q}}\chi_{1}^2\hat{\rho}_{0}| \hat{D}_{33}^2(\operatorname{div}_{z}\hat{u}_{R})|^2\,{\rm d}z+\kappa\varepsilon^2\int_{\hat{Q}\cap \hat{\Omega}}\chi_{1}^2\frac{\hat{\rho}_{0}}{\hat{\theta}_{0}}|\hat{D}_{333}^3\hat{\theta}_{R}|^2\,{\rm d}z\nonumber\\
	&\leq C\big(\|(u_{R},\theta_{R})\|_{H^2}^2+\v^2\|(\mathfrak{R}_{2},\mathfrak{R}_{3})\|_{H^1}^2\big)+ C\varepsilon^2\sum\limits_{j=1}^2\big(\|\nabla_{z}^2 \hat{D}_{j}\hat{u}_{R}\|_{L^2}^2+\|\hat{D}_{jj3}^3\hat{\theta}_{R}\|_{L^2}^2\big)\nonumber\\
	&\quad +\big[C(\mathcal{I}_{1})C_{\sigma}\v^2+C\sigma\big]\|\rho_{R}\|_{H^2}^2+C\sigma^2\v^2\|(u_{R},\theta_{R})\|_{H^3}^2+C\v^2\|\rho_{R}\|_{H^2}\|\mathfrak{R}_{1}\|_{H^2},
	\end{align}
where $C>0$ is independent of $\sigma$ and $\v$.
\end{lemma}

\noindent\textbf{Proof}. We apply $\hat{D}_{3}$ to \eqref{5.7} and multiply the resultant equation by $-\varepsilon^2\hat{D}_{33}^2\operatorname{div}_{z}(\hat{\rho}_{0}\hat{u}_{R})$ to get
%\begin{align}\label{5.12}
%&-(\mu+\zeta)\hat{D}_{33}^2(\operatorname{div}\hat{u}_{R})+\frac{1}{\varepsilon}\hat{D}_{33}^2(\hat{\rho}_{0}\hat{\theta}_{R}+\hat{\rho}_{R}\hat{\theta}_{0})=-\hat{D}_{3}[\hat{\rho}_{0}(\hat{u}_{1}\cdot \nabla_{z}\hat{u}_{R}+\hat{u}_{R}\cdot \nabla_{z}\hat{u}_{1})\cdot \hat{e}_{3}]\nonumber\\
%&\quad +\hat{D}_{33}^2(\hat{\rho}_{R}\hat{\theta}_{1})+\hat{D}_{3}\hat{W}_{2}^3+\hat{D}_{3}[(\hat{\tilde{F}}^{\varepsilon}+\hat{r}_{2})\cdot \hat{e}_{3}]+\mu\hat{D}_{3}[(\Delta_{z}\hat{u}_{R}-\nabla_{z}\operatorname{div}_{z}\hat{u}_{R})\cdot \hat{e}_{3}],
%\end{align}
%Multiplying \eqref{5.12} by $-\chi_{1}^2\hat{D}_{33}^2(\operatorname{div}_{z}(\hat{\rho}_{0}\hat{u}_{R}))$, one obtains that
\begin{align}\label{5.13}
&(\mu+\zeta)\varepsilon^2\int_{\hat{\Omega}\cap \hat{Q}}\hat{\rho}_{0}\chi_{1}^2|\hat{D}_{33}^2(\operatorname{div}_{z}\hat{u}_{R})|^2\,{\rm d}z\nonumber\\
&=\varepsilon\int_{\hat{\Omega}\cap \hat{Q}}\chi_{1}^2\hat{D}_{3 3}^2(\hat{\rho}_{0}\hat{\theta}_{R})\hat{D}_{3 3}^2\operatorname{div}_{z}(\hat{\rho}_{0}\hat{u}_{R})\,{\rm d}z+\varepsilon\int_{\hat{\Omega}\cap \hat{Q}}\chi_{1}^2\hat{D}_{33}^2(\hat{\rho}_{R}\hat{\theta}_{0})\hat{D}_{3 3}^2\operatorname{div}_{z}(\hat{\rho}_{0}\hat{u}_{R})\,{\rm d}z\nonumber\\
&\quad -(\mu+\zeta)\varepsilon^2\int_{\hat{\Omega}\cap \hat{Q}}\chi_{1}^2\hat{D}_{3 3}^2(\operatorname{div}_{z}\hat{u}_{R})[\hat{D}_{33}^2\operatorname{div}_{z}(\hat{\rho}_{0}\hat{u}_{R})-\hat{\rho}_{0}\hat{D}_{33}^2\operatorname{div}_{z}\hat{u}_{R}]\,{\rm d}z\nonumber\\
&\quad +\varepsilon^2\int_{\hat{Q}\cap \hat{\Omega}}\chi_{1}^2\hat{D}_{33}^2(\hat{\rho}_{R}\hat{\theta}_{1}+\hat{\rho}_{1}\hat{\theta}_{R})\hat{D}_{33}^2\operatorname{div}_{z}(\hat{\rho}_{0}\hat{u}_{R})\,{\rm d}z\nonumber\\
%-\varepsilon^2\int_{\hat{Q}\cap \hat{\Omega}}\chi_{1}^2\hat{D}_{33}^2(\hat{\rho}_{1}\hat{\theta}_{R})\hat{D}_{33}^2\operatorname{div}_{z}(\hat{\rho}_{0}\hat{u}_{R})\,{\rm d}z\nonumber\\
&\quad -\mu\varepsilon^2\int_{\hat{Q}\cap \hat{\Omega}}\chi_{1}^2\hat{D}_{3}[(\Delta_{z}\hat{u}_{R}-\nabla_{z}\operatorname{div}_{z}\hat{u}_{R})\cdot \hat{e}_{3}]\hat{D}_{33}^2\operatorname{div}_{z}(\hat{\rho}_{0}\hat{u}_{R}){\rm d}z\nonumber\\
&\quad +\varepsilon^2\int_{\hat{\Omega}\cap \hat{Q}}\chi_{1}^2\hat{D}_{3}[\hat{\rho}_{0}(\hat{u}_{1}\cdot \nabla \hat{u}_{R}+\hat{u}_{R}\cdot \nabla \hat{u}_{1})\cdot \hat{e}_{3}]\hat{D}_{3 3}^2(\operatorname{div}_{z}(\hat{\rho}_{0}\hat{u}_{R}))\,{\rm d}z\nonumber\\
&\quad +\varepsilon^2\int_{\hat{Q}\cap \hat{\Omega}}\chi_{1}^2\hat{D}_{3}(\hat{W}_{2}^3+\hat{F}^{\varepsilon,3}+\hat{\mathfrak{R}}_{2}^3)\hat{D}_{33}^2\operatorname{div}_{z}(\hat{\rho}_{0}\hat{u}_{R})\,{\rm d}z.
\end{align}

We apply $\hat{D}_{33}^2$ to $\eqref{5.2}_{1}$ and multiply the resultant equation by $\varepsilon^2\chi_{1}^2\hat{D}_{33}^2(\hat{\theta}_{0}\hat{\rho}_{R})$ to get
\begin{align}\label{5.12}
	&\varepsilon\int_{\hat{Q}\cap \hat{\Omega}}\chi_{1}^2\hat{D}_{33}^2\operatorname{div}_{z}(\hat{\rho}_{0}\hat{u}_{R})\hat{D}_{33}^2(\hat{\theta}_{0}\hat{\rho}_{R})\,{\rm d}z\nonumber\\
	&=\varepsilon\sum\limits_{m=1}^2\int_{\hat{Q}\cap \hat{\Omega}}\chi_{1}^2\hat{D}_{33}^2(\hat{\theta}_{0}\hat{\rho}_{R})\hat{D}_{33}^2(\hat{D}_{m}w\hat{D}_3(\hat{\rho}_{0}\hat{u}_{R}^{m}))\,{\rm d}z\nonumber\\
	&\quad-\varepsilon^2\int_{\hat{Q}\cap \hat{\Omega}}\chi_{1}^2\hat{D}_{33}^2(\hat{\theta}_{0}\hat{\rho}_{R})\hat{D}_{33}^2\{\operatorname{div}_{z}[\hat{\rho}_{R}(\hat{u}_{1}+\varepsilon(\hat{u}_{2}+\hat{u}_{R}))]+\hat{W}_{1}\}\,{\rm d}z\nonumber\\
	&\quad -\varepsilon^2\int_{\hat{Q}\cap \hat{\Omega}}\chi_{1}^2\hat{D}_{33}^2(\hat{\theta}_{0}\hat{\rho}_{R})\hat{D}_{33}\hat{\mathfrak{R}}_{1}\,{\rm d}x.
	%&\quad-\varepsilon^2\int_{\hat{Q}\cap \hat{\Omega}}\chi_{1}^2\hat{D}_{33}^2(\hat{\theta}_{0}\hat{\rho}_{R})\hat{D}_{33}^2\operatorname{div}_{z}[\hat{\tilde{u}}_{R}(\hat{\rho}_{1}+\varepsilon\hat{\rho}_{2}+\varepsilon^2\hat{\rho}_{3})]\,{\rm d}z\nonumber\\
	%&\quad -\varepsilon^2\int_{\hat{Q}\cap \hat{\Omega}}\chi_{1}^2\hat{D}_{33}^2(\hat{\theta}_{0}\hat{\rho}_{R})\hat{D}_{33}^2\hat{r}_{1}\,{\rm d}z.
\end{align}

Finally, we apply $\hat{D}_{3}$ to $\eqref{5.2}_{3}$ and multiply the resultant equation by $\varepsilon^2\chi_{1}^2\hat{D}_{333}^3(\hat{\rho}_{0}\hat{\theta}_{R})$ to obtain
\begin{align}\label{5.15}
	&\kappa\varepsilon^2\int_{\hat{Q}\cap \hat{\Omega}}\chi_{1}^2\frac{\hat{\rho}_{0}}{\hat{\theta}_{0}}|\hat{D}_{333}^3\hat{\theta}_{R}|^2\,{\rm d}z\nonumber\\
	&=\varepsilon\int_{\hat{Q}\cap \hat{\Omega}}\chi_{1}^2\hat{D}_{333}^3(\hat{\rho}_{0}\hat{\theta}_{R})\hat{D}_{3}\operatorname{div}_{z}(\hat{\rho}_{0}\hat{u}_{R})\,{\rm d}z -\varepsilon\sum\limits_{m=1}^2\int_{\hat{Q}\cap \hat{\Omega}}\chi_{1}^2\hat{D}_{333}^3(\hat{\theta}_{0}\hat{\rho}_{R})\hat{D}_{3}(\hat{D}_{m}w\hat{D}_{3}(\hat{\rho}_{0}\hat{u}_{R}^{m}))\,{\rm d}z\nonumber\\
	&\quad -\kappa\varepsilon^2\int_{\hat{Q}\cap \hat{\Omega}}\chi_{1}^2\hat{D}_{333}^3(\hat{\rho}_{0}\hat{\theta}_{R})\big[\hat{D}_{3}(\frac{\Delta_{z}\hat{\theta}_{R}}{\hat{\theta}_{0}})-\frac{1}{\hat{\theta}}_{0}\hat{D}_{333}^3\hat{\theta}_{R}\big]\,{\rm d}z\nonumber\\
	&\quad +\varepsilon^2P_{1}\int_{\hat{Q}\cap \hat{\Omega}}\chi_{1}^2\hat{D}_{333}^3(\hat{\rho}_{0}\hat{\theta}_{R})\hat{D}_{3}\big(\frac{1}{\hat{\theta}_0}\operatorname{div}_{z}\hat{u}_{R}\big)\,{\rm d}z\nonumber\\
	&\quad +\varepsilon^2\int_{\hat{Q}\cap \hat{\Omega}}\chi_{1}^2\hat{D}_{333}^3(\hat{\rho}_{0}\hat{\theta}_{R})\hat{D}_{3}\big(\frac{1}{\hat{\theta}_{0}}(\hat{\rho}_{0}\hat{\theta}_{R}+\hat{\rho}_{R}\hat{\theta}_{0})\operatorname{div}_{z}\hat{u}_{1}\big)\,{\rm d}z\nonumber\\
	&\quad +\varepsilon^2\int_{\hat{Q}\cap \hat{\Omega}}\chi_{1}^2\hat{D}_{333}^3(\hat{\rho}_{0}\hat{\theta}_{R})\hat{D}_{3}\big(\frac{1}{\hat{\theta}_{0}}(\hat{G}^{\varepsilon})+\hat{\mathfrak{R}}_{3}+\hat{W}_{3}\big)\,{\rm d}z\nonumber\\
	%&\quad +\varepsilon^2\int_{\hat{Q}\cap \hat{\Omega}}\chi_{1}^2\hat{D}_{333}^3(\hat{\rho}_{0}\hat{\theta}_{R})\hat{D}_{3}(\frac{1}{\hat{\theta}_{0}}\hat{\Psi}(\nabla (u_{1}+\varepsilon(\hat{u}_{2}+\tilde{u}_{R}))))\,{\rm d}z\nonumber\\
	&\quad +\varepsilon^2\int_{\hat{Q}\cap \hat{\Omega}}\chi_{1}^2\hat{D}_{333}^3\big(\hat{\rho}_{0}\hat{\theta}_{R})\hat{D}_{3}\Big\{\frac{1}{\hat{\theta}_{0}}[\hat{\rho}_{0}(\hat{u}_{R}\cdot \nabla_{z}\hat{\theta}_{1}+\hat{u}_{1}\cdot \nabla_{z}\hat{\theta}_{R})\nonumber\\
	&\qquad \qquad
	\qquad +\hat{\rho}_{1}\hat{u}_{R}\cdot \nabla_{z}\hat{\theta}_{0}+\hat{\rho}_{R}\hat{u}_{1}\cdot \nabla_{z}\hat{\theta}_{0}]\Big\}\,{\rm d}z.
\end{align}

It follows from \eqref{5.8} and the H\"{o}lder inequality that
\begin{align}\label{5.14}
	&\varepsilon^2\int_{\hat{\Omega}\cap \hat{Q}}\chi_{1}^2\hat{D}_{33}^2(\operatorname{div}_{z}(\hat{\rho}_{0}\hat{u}_{R}))\hat{D}_{3}[(\Delta_{z}\hat{u}_{R}-\nabla_{z}\operatorname{div}_{z}\hat{u}_{R})\cdot \hat{e}_{3}]\,{\rm d}z\nonumber\\
	&\leq \frac{(\mu+\zeta)}{16}\varepsilon^2\int_{\hat{\Omega}\cap \hat{Q}}\hat{\rho}_{0}\chi_{1}^2|\hat{D}_{33}^2\operatorname{div}\hat{u}_{R}|^2\,{\rm d}z+C\varepsilon^2\sum\limits_{j=1}^2\int_{\hat{\Omega}\cap \hat{Q}}\chi_{1}^2|\nabla \hat{D}_{3j}^2\hat{u}_{R}|^2\,{\rm d}z.
\end{align}
%The other terms in \eqref{5.13} can controlled by similar arguments as in Lemma \ref{lem5.1}.
Noting \eqref{5.14}, then applying similar arguments as in Lemma \ref{lem5.1} to \eqref{5.13}--\eqref{5.14}, we conclude \eqref{5.19}. This completes the proof of Lemma \ref{lem5.3}. $\hfill\square$

\medskip

\noindent\textbf{Step 6. {\it A priori} uniform estimates of $\|\rho_{R}\|_{H^2}+\|u_{R}\|_{\mathcal{K}}+\|\theta_{R}\|_{H^3}$}.

Recall \eqref{covering}. Replacing $Q$ by $Q_{k}$ in Lemmas \ref{lem5.0}--\ref{lem5.3}, summing up $k=1,\cdots,N$, and using Lemma \ref{lem4.2}, \eqref{4.1} and \eqref{4.24-0}, we conclude
\begin{align}\label{5.18-0}
&\|(u_{R},\theta_{R})\|_{H^2}^2+\|\frac{\nabla (\rho_{0}\theta_{R}+\theta_{0}\rho_{R})}{\v}\|_{L^2}^2\nonumber\\
&\leq C(\mathcal{I}_{1})C_{\sigma}\|\rho_{R}\|_{H^2}^2+C\|\mathfrak{R}_{1}\|_{H^1}^2+C\|(\mathfrak{R}_{2},\mathfrak{R}_{3})\|_{L^2}^2+C_{\tau,\sigma}\|(u_{R},\theta_{R})\|_{H^1}^2\nonumber\\
&\leq C(\mathcal{I}_{1})C_{\sigma}\|\rho_{R}\|_{H^2}^2+C\|\mathfrak{R}_{1}\|_{H^1}^2+C\|(\mathfrak{R}_{2},\mathfrak{R}_{3})\|_{L^2}^2,
\end{align}
where we have taken $\tau$ and $\sigma$ small enough first, and then $\mathcal{I}_{1}$ small enough. Noting that
$$
\|\rho_{R}\|_{H^1}^2\leq C\Big(\|\theta_{R}\|_{H^1}^2+\|\rho_{R}\|_{L^2}^2+\|\frac{\nabla (\rho_{0}\theta_{R}+\theta_{0}\rho_{R})}{\v}\|_{L^2}^2\Big),
$$
we obtain from \eqref{5.18-0} and \eqref{4.1} that
\begin{align}\label{5.18-1}
\|\rho_{R}\|_{H^1}^2\leq C(\mathcal{I}_{1})C_{\sigma}\|\rho_{R}\|_{H^2}^2+C\|\mathfrak{R}_{1}\|_{H^1}^2+C\|(\mathfrak{R}_{2},\mathfrak{R}_{3})\|_{L^2}^2.
\end{align}
Furthermore, using Lemmas \ref{lem4.3}, \ref{lem5.1}--\ref{lem5.4}, one has
%Since $\partial\Omega$ is compact, there exists finite number of $Q_{k}$'s such that
%\begin{itemize}
	%\item [a)] each $Q_{k}$ corresponds to a $Q$ above;
	%\item [b)] $\partial\Omega\subset \frac{1}{2}Q_{k}$ (here $\frac{1}{2}Q_{k}$ is the cube with half the radius of $Q_{k}$ and with the same center);
	%\item [c)] each $Q_{k}$ intersects at most eight other $Q_j$'s;
	%\item [d)] $w_{k}\in C^{4}(\partial\Omega\cap Q_{k}$ is the corresponding local strengthening function and $\|w_{k}\|_{C^4(\partial\Omega\cap Q_{k})}\leq \delta$.
%\end{itemize}
%Let $\hat{Q}_{k}$ be the corresponding domain in local coordinate. Define $\chi_{k}\in C_{0}^{\infty}(\hat{Q}_{k})$ to be the cutoff function satisfying
%$$
%\chi_{k}=1\quad \text{in }\frac{1}{2}\hat{Q}_{k},\quad \chi_{k}=0\quad \text{in }\hat{Q}_{k}^{c},
%$$
%and
%$$
%\chi_{0}=1\quad \text{in }\Omega\backslash \bigcup\frac{1}{2}Q_{k},\quad \chi_{0}=0\quad \text{in }(\Omega\backslash \bigcup\frac{1}{2}Q_{k})^{c}.
%$$
%Combining \eqref{4.24} and all above estimates established in Steps 1--4 together, we can get
\begin{align}\label{5.18-2}
&\v^2\|\nabla^2\operatorname{div}u_{R}\|_{L^2}^2+\v^2\|\nabla^3 \theta_{R}\|_{L^2}^2\nonumber\\
&\leq [C(\mathcal{I}_{1})C_{\sigma}\v^2+C\sigma+C\tau]\|\rho_{R}\|_{H^2}^2+C_{\tau,\sigma}\|(u_{R},\theta_{R})\|_{H^2}^2+C_{\sigma}\|\rho_{R}\|_{H^1}^2\nonumber\\
&\quad +C\sigma^2\v^2\|(u_{R},\theta_{R})\|_{H^3}^2+C\v^2\big[\|(\mathfrak{R}_{2},\mathfrak{R}_{3})\|_{H^1}^2+\|\rho_{R}\|_{H^2}\|\mathfrak{R}_{1}\|_{H^2}\big],
\end{align}
which, together with using \eqref{5.18-0}--\eqref{5.18-1} and \eqref{4.24}, yields that
\begin{align}\label{5.18}
&\|(\rho_{R},u_{R},\theta_{R})\|_{H^2}^2+\varepsilon^2(\|u_{R}\|_{H^3}^2+\|\theta_{R}\|_{H^3}^2)+\|\frac{\nabla (\rho_{0}\theta_{R}+\theta_{0}\rho_{R})}{\varepsilon}\|_{L^2}^2\nonumber\\
&\leq C\big[\|(\mathfrak{R}_{2},\mathfrak{R}_{3})\|_{L^2}^2+\|\mathfrak{R}_{1}\|_{H^1}^2\big]+C\v^2\big[\|(\mathfrak{R}_{2},\mathfrak{R}_{3})\|_{H^1}^2+\v^2\|\mathfrak{R}_{1}\|_{H^2}^2\big],
\end{align}
where we have taken $\tau$ and $\sigma$ small enough first, and then taken $\mathcal{I}_{1}$ small enough. Here $C$ is a positive constant independent of $\varepsilon$.

Noting \eqref{5.20-1}, and using the classical elliptic estimate, we obtain from \eqref{5.18} that
\begin{align}\label{5.20}
\|\theta_{R}\|_{H^3}&\leq C(\mathcal{I}_{1})\|(\rho_{R},u_{R},\theta_{R})\|_{H^2}+C\|(\mathfrak{R}_{1},\mathfrak{R}_{3})\|_{H^1}\nonumber\\
&\leq C\big[\|\mathfrak{R}_{2}\|_{L^2}+\|(\mathfrak{R}_{1},\mathfrak{R}_{3})\|_{H^1}\big]+C\v\big[\|\mathfrak{R}_{2}\|_{H^1}+\v\|\mathfrak{R}_{1}\|_{H^2}\big].
%&\leq C\big[\|(\mathfrak{R}_{2},\mathfrak{R}_{3})\|_{L^2}^2+\|(\mathfrak{R}_{1},\mathfrak{R}_{3})\|_{H^1}^2\big]+C\v^2\big[\|\mathfrak{R}_{2}\|_{H^1}^2+\v^2\|\mathfrak{R}_{1}\|_{H^2}^2\big].
\end{align}
%by H\"{o}lder inequality, we deduce from \eqref{5.18} that
%\begin{align}\label{5.20}
%&\|\rho_{R}\|_{H^2}^2+\|u_{R}\|_{H^2}^2+\|\theta_{R}\|_{H^2}^2+\varepsilon^2(\|u_{R}\|_{H^3}^2+\|\theta_{R}\|_{H^3}^2)+\|\frac{\nabla (\rho_{0}\theta_{R}+\theta_{0}\rho_{R})}{\varepsilon}\|_{L^2}^2\nonumber\\
%&\leq C_{15}\big(\|r_{2}\|_{H^1}^2+\|r_{3}\|_{H^1}^2\big)+C_{15}\big(\|u_{1}\|_{H^2}^4+\varepsilon^{4}(\|u_{2}\|_{H^2}^4+\|\tilde{u}_{R}\|_{H^2}^4)\big)\nonumber\\
%&\quad+ C_{15}\big[(\|\tilde{u}_{R}\|_{H^2}^2+\varepsilon^4\|\|\tilde{u}_{R}\|_{H^3}^2)(\|\rho_{1}\|_{H^3}^2+\varepsilon^2\|\rho_{2}\|_{H^3}^2+\varepsilon^4\|\rho_{3}\|_{H^3}^2)+\|r_{1}\|_{H^2}^2]\nonumber\\
%&\quad +C_{15}\varepsilon^2\big(\|u_{1}\|_{H^3}^4+\varepsilon^{4}(\|u_{2}\|_{H^3}^4+\|\tilde{u}_{R}\|_{H^3}^4)\big),
%\end{align}
Hence, combining Lemma \ref{lem2.2} and the higher order uniform estimates \eqref{5.18} and \eqref{5.20}, we have following proposition.
\begin{proposition}\label{lem5.4}
Let $\tilde{u}_{R}\in \mathcal{K}$ and $\tilde{\theta}_{R}\in (H_{0}^1\cap H^3)$. For any $0<\v<1$, let $\tilde{\delta}_{0}$ be given in Lemma \ref{lem2.2}. There exists a constant $\tilde{\delta}_{1}>0$ with $\tilde{\delta}_{1}\leq \tilde{\delta}_{0}$ such that if $\delta_{0}\leq \tilde{\delta}_{1}$ and $\|\tilde{u}_{R}\|_{\mathcal{K}}+\|\tilde{\theta}_{R}\|_{H^3}\leq \tilde{\delta}_{1}$, Any strong solutions $(\rho_{R},u_{R},\theta_{R})\in H^2\times \mathcal{K}\times (H_{0}^1\cap H^3)$ of \eqref{3.10} satisfies following {\it a priori} uniform estimate:
\begin{align}\label{5.21}
	&\|\rho_{R}\|_{H^2}+\|u_{R}\|_{\mathcal{K}}+\|\theta_{R}\|_{H^3}+\|\frac{\nabla (\rho_{0}\theta_{R}+\theta_{0}\rho_{R})}{\varepsilon}\|_{L^2}\nonumber\\
	&\leq C\big[\|\mathfrak{R}_{2}\|_{L^2}+\|(\mathfrak{R}_{1},\mathfrak{R}_{3})\|_{H^1}\big]+C\v\big[\|\mathfrak{R}_{2}\|_{H^1}+\v\|\mathfrak{R}_{1}\|_{H^2}\big],
\end{align}
where $C$ is a positive constant independent of $\varepsilon$ and depends continuously on $(\mu,\zeta,\kappa)$. Here $\mathcal{K}$ is the space defined in \eqref{K}.
\end{proposition}

\section{Existence of Strong Solution of \eqref{3.10}}
With the help of {\it a priori} uniform estimate \eqref{5.21}, as in \cite{Valli-1987}, we shall apply a continuity method to show the existence of strong solution $(\rho_{R},u_{R},\theta_{R})\in H^2\times \mathcal{K}\times (H_{0}^1\cap H^3)$ of \eqref{3.10} for any given viscous coefficients pair $(\mu,\zeta)$.
%As indicated in \cite{Valli-1987}, the approximate scheme \eqref{3.11} fail to get the desired estimate due to the term $u_{R}\cdot \nabla \rho_{R}$ in $\eqref{3.11}_{1}$, which breaks down the ellipticity of \eqref{3.10}. Hence we will follow the idea of \cite{Valli-1987} to establish the existence of strong solution of \eqref{3.10}.

%Motivated by \cite{Valli-1987}, we shall apply Leary-Schaulder fixed point theorem to get the above result for all viscous constants $\mu$ and $\zeta$.

\begin{lemma}\label{lem7.2}
	For any fixed $\mu>0, \zeta>0$ and $0<\var<1$. Let $\tilde{\delta}_{1}$ be given in Lemma \ref{lem5.4}. Then there exists $\tilde{\delta}_{2}$ with $\tilde{\delta}_{2}\leq \tilde{\delta}_{1}\leq \tilde{\delta}_{0}$ such that if $\delta_{0}\leq \tilde{\delta}_{2}$ and $\|\tilde{u}_{R}\|_{\mathcal{K}}+\|\tilde{\theta}_{R}\|_{H^3}\leq \tilde{\delta}_{2}$, there exists a unique strong solution $(\rho_{R},u_{R},\theta_{R})\in H^2\times \mathcal{K}\times (H_{0}^1\cap H^3)$. Moreover, it holds that
	\begin{align}\label{g2}
		\|\rho_{R}\|_{H^2}+\|u_{R}\|_{\mathcal{K}}+\|\theta_{R}\|_{H^3}\leq C\big[\|\mathfrak{R}_{2}\|_{L^2}+\|(\mathfrak{R}_{1},\mathfrak{R}_{3})\|_{H^1}\big]+C\v\big[\|\mathfrak{R}_{2}\|_{H^1}+\v\|\mathfrak{R}_{1}\|_{H^2}\big],
	\end{align}
where $C$ is a positive constant independent of $\varepsilon$.
\end{lemma}
\noindent\textbf{Proof.} Let $\mu_{0}$ and $\zeta_{0}$ be given in Lemma \ref{lem7.1}. We take $\delta_{0}$ and $\|\tilde{u}_{R}\|_{\mathcal{K}}+\|\tilde{\theta}_{R}\|_{H^3}$ small enough such that Lemma \ref{lem7.1} holds. For any $t\in [0,1]$, we define
\begin{align*}
	&\mu_{t}=(1-t)\mu_{0}+t\mu,\quad \zeta_{t}=(1-t)\zeta_{0}+t\zeta,\quad \mathfrak{R}=(\mathfrak{R}_{1},\mathfrak{R}_{2},\mathfrak{R}_{3}),\nonumber\\
	&\mathscr{L}_{t}(\rho_{R},u_{R},\theta_{R}):=(\mathscr{L}_{t}^{(1)},\mathscr{L}_{t}^{(2)},\mathscr{L}_{t}^{(3)})(\rho_{R},u_{R},\theta_{R}),
\end{align*}
with
\begin{align}\label{H2}
	&\mathscr{L}_{t}^{(1)}(\rho_{R},u_{R},\theta_{R})=\operatorname{div}[\rho_{R}(u_{1}+\v(u_{2}+\tilde{u}_{R}))]+\frac{1}{\v}\operatorname{div}(\rho_{0}u_{R}),\nonumber\\
	&\mathscr{L}_{t}^{(2)}(\rho_{R},u_{R},\theta_{R})=\mu_{t}\Delta u_{R}+\zeta_{t}\nabla \operatorname{div}u_{R}-\frac{1}{\v}\nabla (\rho_{0}\theta_{R}+\rho_{R}\theta_{0})-\rho_{0}(u_{1}\cdot \nabla u_{R}+u_{R}\cdot \nabla u_{1})\nonumber\\
	&\qquad \qquad \qquad\qquad-\nabla (\rho_{R}\theta_{1}+\rho_{1}\theta_{R})+\tilde{F}^{\v}(\rho_{R},u_{R},\theta_{R}),\nonumber\\
	&\mathscr{L}_{t}^{(3)}(\rho_{R},u_{R},\theta_{R})=\frac{\kappa}{\theta_{0}}\Delta\theta_{R}-\frac{1}{\v}\operatorname{div}(\rho_{0}u_{R})-\frac{P_{1}}{\theta_{0}}\operatorname{div}u_{R}\nonumber\\
	&\qquad \qquad \qquad\qquad-\frac{1}{\theta_{0}}(\rho_{0}\theta_{R}+\rho_{R}\theta_{0})\operatorname{div}u_{1}-\frac{1}{\theta_{0}}\tilde{G}(\rho_{R},u_{R},\theta_{R})\nonumber\\
	&\qquad \qquad \qquad\qquad-\frac{1}{\theta_{0}}\big[\rho_{0}(u_{R}\cdot \nabla \theta_{1}+u_{1}\cdot \nabla \theta_{R})+\rho_{1}u_{R}\cdot \nabla \theta_{0}+\rho_{R}u_{1}\cdot \nabla \theta_{0}\big],
	%&\mathfrak{R}=(\mathfrak{R}_{1},\mathfrak{R}_{2},\mathfrak{R}_{3}),\quad \mathfrak{R}_{1}=:-\operatorname{div}[\tilde{u}_{R}(\rho_{1}+\v\rho_{2}+\v^2\rho_{3})]+r_{1},\quad \mathfrak{R}_{2}=r_{2},\nonumber\\
	%&\mathfrak{R}_{3}=\frac{1}{\theta_{0}}r_{3}-\frac{1}{\theta_{0}}\Psi(\nabla(u_{1}+\v(u_{2}+\tilde{u}_{R}))),
\end{align}
and
\begin{align}\label{H3}
	\begin{aligned}
	&\mathscr{X}=\Big\{(\rho_{R},u_{R},\theta_{R})\in H^2\times \mathcal{K}\times (H_{0}^1\cap H^3)\,\vert\,\int_{\Omega}\rho_{R}\,{\rm d}x=0\Big\},\\
	&\mathscr{Y}=\Big\{\mathfrak{R}=(\mathfrak{R}_{1},\mathfrak{R}_{2},\mathfrak{R}_{3})\in H^2\times H^1\times H^1\,|\,\int_{\Omega}\mathfrak{R}_{1}\,{\rm d}x=0\Big\}.
	\end{aligned}
\end{align}
We shall prove that the set
\begin{align}\label{H4}
	\mathscr{T}=\Big\{t\in [0,1]\,|\,&\text{for each $\mathfrak{R}\in \mathscr{Y}$ there exists a unique solution $(\rho_{R},u_{R},\theta_{R})\in \mathscr{X}$}\nonumber\\
	&\text{of $\mathscr{L}_{t}(\rho_{R},u_{R},\theta_{R})=\mathfrak{R}$}\Big\}
\end{align}
is not empty, open and closed, i.e., $\mathscr{T}=[0,1]$.

It follows from Lemma \ref{lem7.1} that $0\in \mathscr{T}$, which implies $\mathscr{T}\neq \emptyset$. Let now $t_{0}\in \mathscr{T}$, then one has from {\it a priori} uniform estimate \eqref{5.21} that
%from the {a priori} estimate \eqref{5.21} that
$$
\|\mathscr{L}_{t_{0}}^{-1}\|_{\mathscr{Y}\mapsto \mathscr{X}}\leq M_{1}
$$
for some constants $M_{1}=M_{1}(\Omega,\kappa,\mu_{t_{0}},\zeta_{t_{0}})$. For any $\sigma>0$, we can rewrite
$$
\mathscr{L}_{t_{0}+\sigma}(\rho_{R},u_{R},\theta_{R})=\mathfrak{R}
$$
in the following form:
$$
[\mathbf{I}-\sigma\mathscr{L}_{t_{0}}^{-1}(\mathscr{L}_{0}-\mathscr{L}_{1})](\rho_{R},u_{R},\theta_{R})=\mathscr{L}_{t_{0}}^{-1}(\mathfrak{R}).
$$
Hence, it can be solved if
$$
\sigma \|\mathscr{L}_{t_{0}}^{-1}(\mathscr{L}_{0}-\mathscr{L}_{1})\|_{\mathscr{X}\mapsto \mathscr{Y}}\leq M_{1}\sigma\|\mathscr{L}_{0}-\mathscr{L}_{1}\|_{\mathscr{X}\mapsto \mathscr{Y}}\ll 1
$$
i.e., $|\sigma|\ll M_{1}^{-1}\|\mathscr{L}_{0}-\mathscr{L}_{1}\|_{\mathscr{X}\mapsto \mathscr{Y}}$, where we have used the fact that $\mathscr{L}_{0}-\mathscr{L}_{1}$
%$$
%\mathscr{L}_{0}-\mathscr{L}_{1}(\rho_{R},u_{R},\theta_{R})=(0,(\mu_{0}-\mu)\Delta u_{R}+(\zeta_{0}-\zeta)\nabla \operatorname{div}u_{R},0)
%$$
is a bounded linear transformation from $\mathscr{X}$ to $\mathscr{Y}$. Hence $\mathscr{T}$ is open.

It remains to show that $\mathscr{T}$ is closed. Let $t_{n}\in \mathscr{T}$ and $t_{n}\to t_{0}$ as $n\to \infty$. Noting from \eqref{5.21} that for $\delta_{0}$ and $\|\tilde{u}_{R}\|_{\mathcal{K}}+\|\tilde{\theta}_{R}\|_{H^3}$ small enough,
\begin{align}\label{H5}
	\sup_{t_{n}\in [0,1]}\|\mathscr{L}_{t_{n}}^{-1}\|_{\mathscr{Y}\mapsto \mathscr{X}}\leq \max_{t_{n}\in [0,1]}M_{1}(\Omega,\kappa,\mu_{t_{n}},\zeta_{t_{n}})<C(\Omega, \kappa, \mu_{0}, \zeta_{0}, \mu, \zeta)<\infty
\end{align}
Set $(\rho_{R}^{(n)},u_{R}^{(n)},\theta_{R}^{(n)}) =\mathscr{L}_{t_{n}}^{-1}(\mathfrak{R})$, \eqref{H5} implies that there exists a subsequence $(\rho_{R}^{(n_{k})},u_{R}^{(n_{k})},\theta_{R}^{(n_{k})})$ such that $(\rho_{R}^{(n_{k})},u_{R}^{(n_{k})},\theta_{R}^{(n_{k})})\to (\rho_{R},u_{R},\theta_{R})$ weakly in $\mathscr{X}$. %Thus, $\mathscr{L}_{t_{n_{k}}}(\rho_{R}^{(n_{k})},u_{R}^{(n_{k})},\theta_{R}^{(n_{k})})\to \mathscr{L}_{t_{0}}(\rho_{R},u_{R},\theta_{R})$ weakly in $\mathscr{Y}$.
Then it is clear to check that $(\rho_{R},u_{R},\theta_{R})$ satisfy $\mathfrak{R}=\mathscr{L}_{t_{0}}(\rho_{R},u_{R},\theta_{R})$.
%Therefore, $\mathfrak{R}=\mathscr{L}_{t_{0}}(\rho_{R},u_{R},\theta_{R})$, which yields that $t_{0}\in \mathscr{T}$. Hence $\mathscr{T}$ is closed.
Therefore $\mathscr{T}\equiv [0,1]$, which completes the proof of Lemma \ref{lem7.2}. $\hfill\square$
%It is easy to check that $(\rho_{R}^{(k+1)},\theta_{R}^{(k+1)},v_{R}^{(k+1)},q_{R}^{(k+1)})\in H^2\times H_{0}^3\times H_{0}^{3}\times H^{4}$

\medskip

\section{Proof of Theorem \ref{thm1.1}}
We are in a position to establish the existence of nonlinear coupled systems \eqref{lge2}--\eqref{r2}.
Motivated by \cite{Dou-Jiang-Jiang-Yang}, we use the following Tikhonov's fixed point theorem to achieve this.
%\begin{theorem}[The Schauder Fixed Point Theorem {\cite{GT}}]\label{thm6.1}
%Let $X$ be a normal vector space, and let $K\subset X$ be a non-empty, bounded and convex set. Then any given continuous mapping $T: K\mapsto K$, with $T(K)$ being pre-compact, there exists a fixed point $x\in K$ such that $T(x)=x$.
%\end{theorem}
\begin{theorem}[Tikhonov's fixed point theorem {\cite[Page 72, 1.2.11.6]{Novotny-Straskraba}}]\label{thm6.2}
Let $K$ be a nonempty closed convex subset of a separable reflective Banach space $X$, and let $F:K\to K$ be a weakly continuous mapping ${\rm (}${\it i.e.}, if $x_{n}\in K$ and $x_{n}\rightharpoonup x$ weakly in $X$, then $F(x_{n})\rightharpoonup F(x)$ weakly in $X$ as well ${\rm )}$. Then $F$ has at least one fixed point in $K$.
\end{theorem}
For later use, we define a function space $X$ by
$$
X=H_0^1\times H_{0}^1.
$$
Then, for any a small constant $A$, we define a subset $\mathbf{K}_{A}$ of $X$ by
$$
\mathbf{K}_{A}=\{(v,\vartheta)\in \mathcal{K}\times (H_{0}^1\cap H^3) \,|\,\|v\|_{\mathcal{K}}+\|\vartheta\|_{H^3}\leq A\},
$$
where $\mathcal{K}$ is the space defined \eqref{3.1}. By the lower semi-continuity of norms, it is easy to check that $\mathbf{K}_{A}$ is convex.

For any $(\tilde{u}_{R},\tilde{\theta}_{R})\in \mathbf{K}_{A}$ with $A\leq \tilde{\delta}_{2}$ and $\delta_{0}\leq \tilde{\delta}_{2}$, where $\tilde{\delta}_{2}$ is the one given in Lemma \ref{lem7.2}, letting $(\rho_{1},u_{2},\theta_{1},P_{1})$ be the solution of \eqref{3.1} given in Lemma \ref{lem2.1}, $(\rho_{2},\rho_3)$ be given in \eqref{3.7-4}, and $(\rho_{R},u_{R},\theta_{R})$ be the solution of \eqref{3.10} established in Lemma \ref{lem7.2}, we can define a nonlinear operator $\mathcal{Q}$ from $\mathbf{K}_{A}$ to $X$ by
$$
\mathcal{Q}(\tilde{u}_{R},\tilde{\theta}_{R})=(u_{R},\theta_{R}).
$$
To apply Theorem \ref{thm6.2}, we shall show $\mathcal{Q}(\mathbf{K}_{A})\subset \mathbf{K}_{A}$ and $\mathcal{Q}$ is weakly continuous for some small constant $A>0$.

\begin{lemma}\label{lem6.1}
For any fixed $0<\var<1$. Let $\delta_{0}\leq \tilde{\delta}_{2}$ with $\tilde{\delta}_{2}$ given in Lemma \ref{lem7.2}. Then there exists a small positive constant $A_{0}$ satisfying $A_{0}\leq \min\{\tilde{\delta}_{2},1\}$ such that
for any given $A\leq A_{0}$, there exists $\tilde{\delta}_{3}>0$ with $\tilde{\delta}_{3}\leq \tilde{\delta}_{2}$ such that if $\delta_{0}\leq \tilde{\delta}_{3}$, then $\mathcal{Q}(\mathbf{K}_{A})\subset \mathbf{K}_{A}$ and $\mathcal{Q}$ is weakly continuous.
\end{lemma}

\noindent\textbf{Proof}. Since the proof is very long, we divide into two steps.

{\it Step 1.} Recalling Lemma \ref{lem1}, \eqref{r1-1}--\eqref{r3}, \eqref{3.1-1}, \eqref{3.7-2}--\eqref{3.7-3}, \eqref{3.10-1} and \eqref{I0}--\eqref{I1}, one can obtain
%For any $\|\tilde{u}_{R}\|_{\mathcal{K}}+\|\tilde{\theta}_{R}\|_{\mathcal{K}}\leq A$,
\begin{align}
	&\|\mathfrak{R}_{1}\|_{H^1}+\v\|\mathfrak{R}_{1}\|_{H^2}\leq
	C\|\tilde{u}_{R}\|_{\mathcal{K}}(\|\rho_{1}\|_{H^3}+\v\|\rho_{2}\|_{H^3}+\v^2\|\rho_{3}\|_{H^3})+C(\mathcal{I}_{0})\nonumber\\
	&\qquad \qquad \qquad\qquad\,\,\,\leq C\delta_{0}(1+\|\tilde{u}_{R}\|_{\mathcal{K}}),\label{6.1}\\
	%C(\|\rho_{1}\|_{H^3}\|u_{2}\|_{H^3}+\|\rho_{2}\|_{H^3}\|u_{1}\|_{H^3})+ C\varepsilon(\|u_{2}\|_{H^3}\|\rho_{2}\|_{H^3}+\|\rho_{3}\|_{H^3}\|u_{1}\|_{H^3})\nonumber\\
	%&\qquad\qquad +C\varepsilon^2\|\rho_{3}\|_{H^3}\|u_{2}\|_{H^3} \leq C\delta_{0}(\|\tilde{u}_{R}\|_{H^2}+\delta_{0}),\label{6.1}\\
	%&\|r_{2}\|_{H^1}+\|r_{3}\|_{H^1}\leq C\delta_{0}^2\|\tilde{u}_{R}\|_{H^2},\label{6.2}
	&\|\mathfrak{R}_{2}\|_{L^2}+\v\|\mathfrak{R}_{2}\|_{H^1}+\|\mathfrak{R}_{3}\|_{H^1}\leq C(\mathcal{I}_{0})+C[\|u_{1}\|_{H^3}^2+\v^2(\|u_{2}\|_{H^3}^2+\|\tilde{u}_{R}\|_{H^3}^2)]\nonumber\\
	&\qquad\qquad \qquad \qquad \qquad \qquad\quad \leq C\delta_{0}(1+\|\tilde{u}_{R}\|_{H^2})+C\|\tilde{u}_{R}\|_{\mathcal{K}}^2.\label{6.2}
\end{align}
Substituting \eqref{6.1}--\eqref{6.2} into \eqref{5.21}, we obtain from $(\tilde{u}_{R},\tilde{\theta}_{R})\in \mathbf{K}_{A}$ that
\begin{align}\label{6.3}
\|u_{R}\|_{\mathcal{K}}+\|\theta_{R}\|_{H^3}\leq C\delta_{0}(1+\|\tilde{u}_{R}\|_{\mathcal{K}})+C\|\tilde{u}_{R}\|_{\mathcal{K}}^2\leq C\delta_{0}(1+A)+CA^2. %C(\|u_{1}\|_{H^2}^2+\varepsilon\|u_{1}\|_{H^3}^2)+C\delta_{0}(\|\tilde{u}_{R}\|_{\mathcal{K}}+\delta_{0})\leq C(\delta_{0}^2+\delta_{0}A).
\end{align}
Taking $A\ll 1$ and then $\delta_{0}\ll 1$, we deuce from \eqref{6.3} that $(u_{R},\theta_{R})\in \mathbf{K}_{A}$, which implies that $\mathcal{Q}(\mathbf{K}_{A})\subset \mathbf{K}_{A}$.

{\it Step 2.} Recalling the definition of weakly continuous, to show $\mathcal{Q}$ is weakly continuous, it suffices to show $\mathcal{Q}$ is continuous on $\mathbf{K}_{A}$ with respect to the norm of $X$.

Let $(u_{R}^{(i)},\theta_{R}^{(i)})=\mathcal{Q}(\tilde{u}_{R}^{(i)},\tilde{\theta}_{R}^{(i)})$, $i=1,2$. In particular, Let $(\rho_{1}^{(i)},u_{2}^{(i)},\theta_{1}^{(i)},P_{3}^{(i)})\in H^{4}\times H^4\times H^4\times H^3$ be the solution of \eqref{3.1}, and $(\rho_{R}^{(i)},u_{R}^{(i)},\theta_{R}^{(i)})\in H^2\times \mathbf{K}_{A}$ be the solution of \eqref{3.10} for given $(\tilde{u}_{R}^{(i)},\tilde{\theta}_{R}^{(i)})\in \mathbf{K}_{A}$ respectively. We denote
$$
\begin{aligned}
&\bar{\rho}_{j}=\rho_{j}^{(1)}-\rho_{j}^{(2)},j=1,2,3,\quad \bar{\rho}_{R}=\rho_{R}^{(1)}-\rho_{R}^{(2)},\quad \bar{P}_{3}={P}_{3}^{(1)}-{P}_{3}^{(2)},\\ &\bar{u}_{2}=u_{2}^{(1)}-u_{2}^{(2)},\quad \bar{u}_{R}=u_{R}^{(1)}-u_{R}^{(2)},\quad \bar{\theta}_{1}=\theta_{1}^{(1)}-\theta_{1}^{(2)},\quad \bar{\theta}_{R}=\theta_{R}^{(1)}-\theta_{R}^{(2)}\\
&\bar{\tilde{u}}_{R}=\tilde{u}_{R}^{(1)}-\tilde{u}_{R}^{(2)},\quad \bar{\tilde{\theta}}_{R}=\tilde{\theta}_{R}^{(1)}-\tilde{\theta}_{R}^{(2)}.
\end{aligned}
$$
Then $(\bar{\rho}_{1},\bar{u}_{2},\bar{\theta}_{1},\bar{P}_{3})$ satisfies
\begin{align}\label{6.4}
	\left\{
	\begin{aligned}
		&\nabla (\theta_{0}\bar{\rho}_{1}+\rho_{0}\bar{\theta}_{1})=0,\quad \int_{\Omega}\bar{\rho}_{1}\,{\rm d}x=0,\\
		&\operatorname{div}(\rho_{0}\bar{u}_{2})=-\operatorname{div}(\bar{\rho}_{1}u_{1}),\\
		&\rho_{0}(u_{1}\cdot \nabla)\bar{u}_{2}+\rho_{0}(\bar{u}_{2}\cdot \nabla)u_{1}+\nabla \bar{P}_{3}=-\bar{\rho}_{1}(u_{1}\cdot \nabla)u_{1}+\mu\Delta \bar{u}_{2}+\zeta\nabla\operatorname{div}\bar{u}_{2},\\
		&\kappa\Delta\bar{\theta}_{1}=-\theta_{0}(\bar{u}_{2}\cdot \nabla)\rho_{0}+\bar{\rho}_{1}(u_{1}\cdot \nabla)\theta_{0}+\rho_{0}(u_{1}\cdot \nabla)\bar{\theta}_{1}+(\rho_{0}\bar{\theta}_{1}+\bar{\rho}_{1}\theta_{0})\operatorname{div}u_{1}\\
		&\qquad\quad\,\, +(\rho_{0}\theta_{0})\operatorname{div}\bar{u}_{2}-2\theta_{0}(\bar{\tilde{u}}_{R}\cdot \nabla)\rho_{0},\\
		&\bar{u}_{2}=0,\quad \bar{\theta}_{1}=0,\quad \text{on }\partial\Omega.
	\end{aligned}
	\right.
\end{align}
By $\eqref{6.4}_{1}$, one has
\begin{align}\label{6.4-1}
\bar{\rho}_1=\frac{1}{\theta_{0}}(\bar{P}_{1}-\rho_{0}\bar{\theta}_{1})\quad \text{with }\bar{P}_{1}=\frac{\int_{\Omega}\theta_{0}^{-1}\rho_{0}\bar{\theta}_{1}\,{\rm d}x}{\int_{\Omega}\theta_{0}^{-1}\,{\rm d}x}.
\end{align}

Let $\bar{v}_2=\rho_{0}\bar{u}_{2}$, then by similar calculations as in \eqref{3.2-1}, we can rewrite \eqref{6.4} as
\begin{align}\label{6.4-2}
	\left\{
\begin{aligned}
	&\operatorname{div}\bar{v}_{2}=\operatorname{div}\big[\frac{\rho_{0}\bar{\theta}_{1}}{\theta_{0}}u_{1}\big]-\bar{P}_{1}\operatorname{div}\big[\frac{u_{1}}{\theta_{0}}\big],\\
	%\frac{\rho_{0}\bar{\theta}_{1}}{\theta_{0}}\operatorname{div}u_{1}+u_{1}\cdot \nabla(\frac{\rho_{0}\bar{\theta}_{1}}{\theta_{0}})-\frac{\bar{P}_{1}}{\theta_{0}}\operatorname{div}u_{1}-\bar{P}_{1}u_{1}\cdot\nabla(\frac{1}{\theta_{0}}),\\
	&-\mu\Delta \bar{v}_{2}+\nabla (P_{0}\bar{P}_{3})=(\theta_{0}-1)\mu\Delta \bar{v}_{2}-\rho_{0}(u_{1}\cdot \nabla \theta_{0})\bar{v}_{2}-P_{0}(u_{1}\cdot \nabla)\bar{v}_{2}-P_{0}(\bar{v}_{2}\cdot \nabla)u_{1}\\
	&\qquad\qquad\qquad\qquad\quad\,\, +\mu\bar{v}_{2}\Delta\theta_{0}+2\mu\nabla\theta_{0}\cdot \nabla (\bar{v}_{2})^{t}-\rho_{0}(\bar{P}_{1}-\rho_{0}\bar{\theta}_{1})(u_{1}\cdot \nabla)u_{1}\\
	&\qquad\qquad\qquad\qquad\quad\,\,+\zeta \nabla \big[\rho_{0}\bar{\theta}_{1}\operatorname{div}u_{1}\theta_{0}u_{1}\cdot \nabla (\frac{\rho_{0}\bar{\theta}_{1}}{\theta_{0}}{\theta_{0}})-\bar{P}_{1}\operatorname{div}u_{1}-\bar{P}_{1}\theta_{0}u_{1}\cdot \nabla (\frac{1}{\theta_{0}})\big]\\
	&\qquad\qquad\qquad\qquad\quad\,\,+\zeta\nabla (\bar{v}_{2}\cdot \nabla \theta_{0}),\\
	& \kappa \Delta\bar{\theta}_{1}=-\frac{\theta_{0}}{\rho_{0}}\bar{v}_{2}\cdot \nabla\rho_{0}+(\bar{P}_{1}-\rho_{0}\bar{\theta}_{1})(u_{1}\cdot \nabla)\theta_{0}+\rho_{0}(u_{1}\cdot \nabla)\bar{\theta}_{1}+\bar{v}_{2}\cdot \nabla \theta_{0}\\
	&\qquad \qquad +\rho_{0}\bar{\theta}_{1}\operatorname{div}u_{1}+\theta_{0}u_{1}\cdot \nabla(\frac{\rho_{0}\bar{\theta}_{1}}{\theta_{0}})-\bar{P}_{1}\theta_{0}u_{1}\cdot (\frac{1}{\theta_{0}})-2\theta_{0}(\bar{\tilde{u}}_{R}\cdot \nabla )\rho_{0},\\
	&\bar{v}_{2}=0,\quad \bar{\theta}_{2}=0\quad \text{on }\partial \Omega.
\end{aligned}
\right.
\end{align}
Noting the smallness of $\delta_{0}$, by similar arguments as in \eqref{3.7-1}, one has
\begin{align}\label{6.5}
\|\bar{v}_{2}\|_{H^3}+\|\bar{P}_{3}\|_{H^2}+\|\bar{\theta}_{1}\|_{H^3}\leq C\|\bar{\tilde{u}}_{R}\|_{H^1}\|\nabla \theta_{0}\|_{H^2},
\end{align}
which yields that
\begin{align}\label{6.5-2}
\|\bar{u}_{2}\|_{H^3}+\|\bar{\rho}_{1}\|_{H^3}\leq C\|\bar{\tilde{u}}_{R}\|_{H^1}\|\nabla \theta_{0}\|_{H^2}.
\end{align}
Noting \eqref{3.7-4}--\eqref{3.7-3}, we have
\begin{align}\label{6.5-1}
\|\bar{\rho}_{2}\|_{H^3}+\|\bar{\rho}_{3}\|_{H^3}\leq C\|\bar{\tilde{u}}_{R}\|_{H^1}\|\nabla \theta_{0}\|_{H^2},
\end{align}
where we point out that $P_{2}$, which is determined by \eqref{ge2}, remains unchanged for $(\tilde{u}_{R}^{(i)},\tilde{\theta}_{R}^{(i)})\,\,(i=1,2)$.

For $(\bar{\rho}_{R},\bar{u}_{R},\bar{\theta}_{R})$, it satisfies
\begin{align*}
\left\{
\begin{aligned}
	&\operatorname{div}[\bar{\rho}_{R}(u_{1}+\v(u_{2}^{(1)}+u_{R}^{(1)}))]+\v\operatorname{div}[\rho_{R}^{(2)}(\bar{u}_{2}+\bar{\tilde{u}}_{R})]=-\frac{1}{\v}\operatorname{div}(\rho_{0}\bar{u}_{R})+(\mathfrak{R}_{1}^{(1)}-\mathfrak{R}_{1}^{(2)}),\\
	&\mu\Delta \bar{u}_{R}+\zeta\nabla \operatorname{div}\bar{u}_{R}=\frac{1}{\varepsilon}\operatorname{div}(\rho_{0}\bar{\theta}_{R}+\bar{\rho}_{R}\theta_{0})+\rho_{0}(u_{1}\cdot \nabla \bar{u}_{R}+\bar{u}_{R}\cdot \nabla u_{1})+\nabla (\bar{\rho_{R}}\theta_{1}^{(1)}+\rho_{R}^{(2)}\bar{\theta}_{1})\\
	&\qquad\qquad\qquad\qquad\quad+\nabla (\bar{\rho}_{1}\theta_{R}^{(1)}+\rho_{1}^{(2)}\bar{\theta}_{R}) +(\tilde{F}^{\varepsilon,(1)}-\tilde{F}^{\varepsilon,(2)})+(\mathfrak{R}_{2}^{(1)}-\mathfrak{R}_{2}^{(2)}),\\
	&\frac{\kappa}{\theta_{0}}\Delta \bar{\theta}_{R}=\frac{1}{\varepsilon}\operatorname{div}(\rho_{0}\bar{u}_{R})+\frac{\bar{P}_{1}}{\theta_{0}}\operatorname{div}u_{R}^{(1)}+\frac{P_{1}^{(2)}}{\theta_{0}}\operatorname{div}\bar{u}_{R}+\frac{1}{\theta_{0}}(\rho_{0}\bar{\theta}_{R}+\bar{\rho}_{R}\theta_{0})\operatorname{div}u_{1}\\
	&\qquad\qquad+\frac{1}{\theta_{0}}\big\{\rho_{0}(\bar{u}_{R}\cdot \nabla \theta_{1}^{(1)}+u_{R}^{(2)}\cdot\nabla\bar{\theta}_{1}+u_{1}\cdot \nabla \bar{\theta}_{R})+\bar{\rho_{1}}{u}_{R}^{(1)}\cdot \nabla\theta_{0}\\
	&\qquad\qquad\qquad\quad+\rho_{1}^{(2)}\bar{u}_{R}\cdot \nabla \theta_{0}+\bar{\rho}_{R}u_{1}\cdot \nabla \theta_{0}\big\}+\frac{1}{\theta_{0}}(\tilde{G}^{\v,(1)}-\tilde{G}^{\v,(2)})+(\mathfrak{R}_{3}^{(1)}-\mathfrak{R}_{3}^{(2)}),\\
	&\bar{u}_{R}\vert_{\partial\Omega}=\bar{\theta}_{R}\vert_{\Omega}=0,\qquad \int_{\Omega}\bar{\rho}_{R}\,{\rm d}x=0,
\end{aligned}
\right.
\end{align*}
where $F^{\varepsilon,(i)}$, $G^{\varepsilon,(i)}$, $\mathfrak{R}_{1}^{(i)}$, $\mathfrak{R}_{2}^{(i)}$ and $\mathfrak{R}_{3}^{(i)}$ represent the terms in \eqref{4.5}--\eqref{4.8} and \eqref{3.10-1} with $(\rho_{1},\rho_{2},\rho_{3},\rho_{R},u_{2},u_{R},\theta_{1},\theta_{R})$ replacing $(\rho_{1}^{(i)},\rho_{2}^{(i)},\rho_{3}^{(i)},\rho_{R}^{(i)},u_{2}^{(i)},u_{R}^{(i)},\theta_{1}^{(i)},\theta_{R}^{(i)})$ respectively.
%and
%\begin{align}\label{6.7}
%\bar{\Psi}&=2\mu D(\varepsilon\bar{u}_{2}+\varepsilon\bar{\tilde{u}}_{R}):[D(u_{1}+\varepsilon u_{2}^{(1)}+\varepsilon \tilde{u}_{R}^{(1)})+ D(u_{1}+\varepsilon u_{2}^{(2)}+\varepsilon \tilde{u}_{R}^{(2)})]\nonumber\\
%&\quad +\lambda\operatorname{div}(\varepsilon\bar{u}_{2}+\varepsilon\bar{\tilde{u}}_{R})[\operatorname{div}(u_{1}+\varepsilon u_{2}^{(1)}+\varepsilon \tilde{u}_{R}^{(1)})+\operatorname{div}(u_{1}+\varepsilon u_{2}^{(2)}+\varepsilon \tilde{u}_{R}^{(2)})],
%\end{align}
%and
%\begin{align*}%\label{r1-2}
%&\bar{r}_{1}=-\operatorname{div}(\bar{\rho}_{2}u_{1})-\operatorname{div}(\bar{\rho}_{1}u_{2}^{(1)}+\rho_{1}^{(2)}\bar{u}_{2})-\varepsilon\operatorname{div}(\bar{\rho}_{2}u_{2}^{(1)}+\rho_{2}^{(2)}\bar{u}_{2}+\bar{\rho}_{3}u_{1})-\varepsilon^2\operatorname{div}(\bar{\rho}_{3}u_{2}^{(1)}+\rho_{3}^{(2)}\bar{u}_{2}).
%\end{align*}

As in \eqref{I0}--\eqref{I1}, we denote
\begin{align}\label{I0-1}
\tilde{\mathcal{I}}_{0}=\|(\rho_{0}-P_{0},u_{1},\theta_{0}-1)\|_{H^3}+\max_{i=1,2}\{\|\rho_{1}^{(i)},u_{2}^{(i)},P_{3}^{(i)},\theta_{1}^{(i)}\|_{H^3}+\|(\rho_{2}^{(i)},\rho_{3}^{(i)})\|_{H^3}\}\lesssim \delta_{0},
\end{align}
\begin{align}\label{I1-1}
\tilde{\mathcal{I}}_{1}=\tilde{\mathcal{I}}_{0}+\max_{i=1,2}\{\|\tilde{u}_{R}^{(i)}\|_{\mathcal{K}}+\|\tilde{\theta}_{R}^{(i)}\|_{H^3}\}\lesssim \delta_{0}+\max_{i=1,2}\{\|\tilde{u}_{R}^{(i)}\|_{\mathcal{K}}+\|\tilde{\theta}_{R}^{(i)}\|_{H^3}\},
\end{align}
and $C(\tilde{\mathcal{I}}_{0})$ and $C(\tilde{\mathcal{I}}_{1})$ as two small positive constants depending on $\tilde{\mathcal{I}}_{0}$ and $\tilde{\mathcal{I}}_{1}$ respectively satisfying $C(\tilde{\mathcal{I}}_{0}), C(\tilde{\mathcal{I}}_{1})\to 0$ as $\tilde{\mathcal{I}}_{0},\tilde{\mathcal{I}}_{1}\to 0$.

Using similar arguments as in the proof of Lemma \ref{lem2.2}, we have
\begin{align}\label{6.8}
\|\bar{\rho}_{R}\|_{L^2}^2+\|(\bar{u}_{R},\bar{\theta}_{R})\|_{H^1}^2
&\leq C(\tilde{\mathcal{I}}_{1})\big[\|(\bar{\rho}_{1},\bar{u}_{2},\bar{\theta}_{1})\|_{H^3}^2+\|(\bar{\rho}_{2},\bar{\rho}_{3})\|_{H^3}^2+\|(\bar{\tilde{u}}_{R},\bar{\tilde{\theta}}_{R})\|_{H^1}^2\big]\nonumber
\nonumber\\
&
\quad+C(\tilde{\mathcal{I}}_{1})\big[\|\bar{\rho}_{R}\|_{L^2}+\|(\bar{u}_{R},\bar{\theta}_{R})\|_{H^1}^2\big] +C\sum\limits_{j=1}^3\|\mathfrak{R}_{j}^{(1)}-\mathfrak{R}_{j}^{(2)}\|_{L^2}^2.
%&\leq
%C[(\|\rho_{1}^{(2)}\|_{H^2}+\varepsilon\|\rho_{2}^{(2)}\|_{H^2}+\varepsilon^2\|\rho_{3}^{(2)}\|_{H^2})\|\bar{\tilde{u}}_{R}\|_{H^1}+\|\bar{\tilde{u}}_{R}^{(1)}\|_{H^3}(\|\bar{\rho}_{1}\|_{H^2}+\varepsilon\|\bar{\rho}_{2}\|_{H^2}+\varepsilon^2\|\bar{\rho}_{3}\|_{H^2})]^2\nonumber\\
%&\quad +C(\|\bar{r}_{1}\|_{L^2}^2+\|r_{2}^{(1)}-r_{2}^{(2)}\|_{L^2}^2+\|r_{3}^{(1)}-r_{3}^{(2)}\|_{L^2}^2+\|\bar{\Psi}\|_{L^2}^2)\nonumber\\
%&\quad +C\varepsilon^{2} \sum\limits_{i,j=1}^2 \|\rho_{R}^{(i)}\|_{H^2}\big[\|\bar{\tilde{u}}_{R}\|_{H^1}^2+\|\bar{u}_{2}\|_{H^1}^2+(\|\tilde{u}_{R}^{(j)}\|_{H^2}^2+\|\tilde{u}_{2}^{(j)}\|_{H^2}^2)(\|\bar{u}_{2}\|_{H^1}^2+\|\bar{u}_{R}\|_{H^1}^2)+\|\bar{\theta}_{R}\|_{H^1}^2\big]\nonumber\\
%&\quad +C\varepsilon^2 \sum\limits_{i,j=1}^2 \|\rho_{R}^{(i)}\|_{H^2}\big[\|\bar{\theta}_{1}\|_{H^1}^2+\|\bar{\theta}_{R}\|_{H^1}^2+(\|\tilde{u}_{R}^{(j)}\|_{H^2}^2+\|\tilde{u}_{2}^{(j)}\|_{H^2}^2)(\|\bar{\theta}_{1}\|_{H^1}^2+\|\bar{\theta}_{R}\|_{H^1}^2)\nonumber\\
%&\qquad\qquad \qquad \qquad \quad\quad  +(\|\bar{\tilde{u}}_{R}\|_{H^1}^2+\|\bar{\tilde{u}}_{2}\|_{H^1}^2)(\|{\theta}_{1}^{(j)}\|_{H^2}^2+\|{\theta}_{R}^{(j)}\|_{H^2}^2)\big],
\end{align}
%where we use the smallness of $\varepsilon$ and $A$.
%It follows from \eqref{6.5}--\eqref{6.5-1} and a
A direct calculation shows that
\begin{align}\label{6.9}
&\sum\limits_{j=1}^3\|\mathfrak{R}_{j}^{(1)}-\mathfrak{R}_{j}^{(2)}\|_{L^2}^2\leq C(\tilde{\mathcal{I}}_{1})\big[\|(\bar{\rho}_{1},\bar{u}_{2},\bar{\theta}_{1})\|_{H^3}^2+\|(\bar{\rho}_{2},\bar{\rho}_{3})\|_{H^3}^2+\|\bar{\tilde{u}}_{R}\|_{H^1}^2\big].
\end{align}
Substituting \eqref{6.9} into \eqref{6.8}, and using \eqref{6.5}--\eqref{6.5-1} and the smallness of $\tilde{\mathcal{I}}_{1}$, we obtain
$$
\|\bar{\rho}_{R}\|_{L^2}^2+\|\bar{u}_{R}\|_{H^1}^2+\|\bar{\theta}_{R}\|_{H^1}^2\leq C(\tilde{\mathcal{I}}_{1})(\|\bar{\tilde{u}}_{R}\|_{H^1}^2+\|\bar{\tilde{\theta}}_{R}\|_{H^1}^2) %C(\varepsilon^2+\|\nabla\theta_{0}\|_{H^2}^2)(\|\bar{\tilde{u}}_{R}\|_{H^1}^2+\|\bar{\tilde{\theta}}_{R}\|_{H^1}^2),
$$
which implies $\mathcal{Q}$ is continuous on $\mathbf{K}_{A}$ with respect to the norm of $X$ provided $\delta_{0}$ and $A$ are small enough. Therefore the proof of Lemma \ref{lem6.1} is complete. $\hfill\square$

\medskip

\textbf{Proof of Theorem \ref{thm1.1}}. Using Lemmas \ref{thm6.2} and \ref{lem6.1}, we get a fixed point of $\mathcal{Q}$ in $\mathbf{K}_{A}$, which implies the existence of solution
\begin{align*}
&(\rho_{1}, u_{2},\theta_{1},P_{3},\rho_{2},\rho_{3})\in H^4\times (H_{0}^1\cap H^4)\times (H_{0}^1\cap H^4)\times H^3\times H^3\times H^3,\\
&(\rho_{R},u_{R},\theta_{R})\in
H^2\times \mathbf{K}_{A},
\end{align*}
 to the coupled systems \eqref{lge2}--\eqref{r2}. Furthermore, applying the same arguments in \eqref{6.4}--\eqref{6.9} again, and using the smallness of $\delta_{0}$ and $A$, one has
\begin{align*}
&\|(\bar{\rho}_{1},\bar{\rho}_{2},\bar{\rho}_{3})\|_{H^3}+\|(\bar{u}_{2},\bar{\theta}_{1})\|_{H^3}+\|\bar{P}_{3}\|_{H^3}\leq C\|\bar{u}_{R}\|_{H^1}\|\nabla \theta_{0}\|_{H^2},\\
&\|\bar{\rho}_{R}\|_{L^2}^2+\|\bar{u}_{R}\|_{H^1}^2+\|\bar{\theta}_{R}\|_{H^1}^2\leq \frac{1}{2}(\|\bar{{u}}_{R}\|_{H^1}^2+\|\bar{{\theta}}_{R}\|_{H^1}^2),
\end{align*}
which yields the uniqueness of solution $(\rho_{1}, u_{2},\theta_{1},P_{3},\rho_{2},\rho_{3}, \rho_{R},u_{R},\theta_{R})$ to the coupled systems \eqref{lge2}--\eqref{r2}. Moreover, we obtain from \eqref{5.21} that
\begin{align}\label{5.21-0}
	&\|\rho_{R}\|_{H^2}+\|u_{R}\|_{\mathcal{K}}+\|\theta_{R}\|_{H^3}+\|\frac{\nabla (\rho_{0}\theta_{R}+\theta_{0}\rho_{R})}{\varepsilon}\|_{L^2}\leq C\|u_{1}\|_{H^3}^2
\end{align}
where the constant $C>0$ is independent of $\varepsilon$.

Finally, let
$$
\rho^{\v}=\rho_{0}+\varepsilon\rho_{1}+\varepsilon^2\rho_{2}+\varepsilon^3\rho_{3}+\varepsilon^2\rho_{R},\quad \mathfrak{u}^{\v}=\varepsilon u_{1}+\varepsilon^2 u_{2}+\varepsilon^2u_{R},\quad \theta^{\v}=\theta_{0}+\varepsilon\theta_{1}+\varepsilon^2\theta_{R},
$$
we obtain a unique solution $(\rho,\mathfrak{u},\theta)$ of \eqref{sNS} satisfying the expansion \eqref{expansion} and \eqref{c}. This completes the proof of Theorem \ref{thm1.1}. $\hfill\square$

\bigskip

\noindent{\bf Acknowledgments.}
Feimin Huang’s research is partially supported by National Key R\&D Program of China, grant No. 2021YFA1000800, and National Natural Sciences Foundation of China, grant No. 12288201. Yong Wang's research is partially supported by the National Natural Science Foundation of China, grants No. 12022114 and No. 12288201,
and CAS Project for Young Scientists in Basic Research, grant No. YSBR-031.
%The paper is partially supported by National Key R\&D Program of China No. 2021YFA1000800, and National Natural Sciences Foundation of China No. 12288201, 12022114, and CAS Project for Young Scientists in Basic Research, Grant No. YSBR-031.
%Yong Wang's research is partially supported by
%the National Natural Science Foundation of China

\appendix
\renewcommand{\appendixname}{Appendix~\Alph{section}}


\begin{thebibliography}{99}
	%\bibitem{Alazard-2006} T. Alazard, Low Mach number limit of the full Navier-Stokes equations, Arch. Ration. Mach. Anal, 180 (2006), 1--73.
	
	\bibitem{Alazard-2008} T. Alazard, A minicourse on the low Mach number, Discrete Contin. Dyn. Syst. Ser. S, 1 (2008), 365--404.
	
	\bibitem{Bardos-Levermore-Ukai-Yang} C. Bardos, C. D. Levermore, S. Ukai and T. Yang, Kinetic equations: fluid dynamical limits and viscous heating, Bull. Inst. Math. Acad. Sin. (N.S.), 3 (2008), pp. 1--49.
	
	\bibitem{Veiga-1987} H. Beir${\rm \tilde{a}}$o de Veiga, An $L^p$-theory for the n-dimensional, stationary,
	compressible Navier–Stokes equations, and incompressible limit for compressible
	fluids. The equilibrium solutions. Comm. Math. Phys., 109 (1987), 229--248.
	
	\bibitem{Boyer-Fabrie} F. Boyer and P. Fabrie, {\it Mathematical Tools for the Study of the Incompessible Navier-Stokes Equations and related models}, Springer, New York, 2013.
	
    \bibitem{Chen-Chen-Liu-Sone} C.-C. Chen, I.-K. Chen, T.-P. Liu and Y. Sone, Thermal transpiration for the linearized Boltzmann equation, Comm. Pure Appl. Math., 60 (2007), pp. 147–163.
	
	\bibitem{Choe-Jin} H. Choe and B. Jin, Existence of solutions of stationary compressible Navier-Stokes equations with large force, J. Funct. Anal., 177 (2000), 54--88.
	
	\bibitem{Dou-Jiang-Jiang-Yang} C.-S. Dou, F. Jiang, S. Jiang and Y.-F. Yang, Existence of strong solutions to the steady Navier-Stokes equations for a compressible heat-conductive fluid with large forces. J. Math. Pures Appl. 103  (2015), no.5, 1163--1197.	
	
	\bibitem{Esposito-Guo-Marra-Wu-2023-1} R. Esposito, Y. Guo, R. Marra and L. Wu, Ghost effect from the Boltzmnn theory, arXiv: 2301.09427v3.
	
	\bibitem{Esposito-Guo-Marra-Wu-2023-2} R. Esposito, Y. Guo, R. Marra and L. Wu, Ghost effect from the Boltzmnn theory: Expansion with remainder term. arXiv: 2301.09560v3.
	
	\bibitem{Feireisl-2018} E. Feireisl, Singular limits for models of compressible, viscous, heat conducting, and/or rotating fluids. In Handbook
	of mathematical analysis in mechanics of viscous fluids, pages 2771--2825, Springer, Cham, 2018.
	
	%\bibitem{Galdi} G.-P. Galdi, {\it An Introduction to the Mathematical Theory of the Navier-Stokes Equations: Steady-State Problems}, Second Edition, Springer, 2011.
	
	\bibitem{Huang-2015} F.-M. Huang, Thermal creep flow for the Boltzmann equation, Chin. Ann. Math. Ser. B, 36 (2015), 855--870.
	
	\bibitem{Huang-Tan} F.-M. Huang and W. Tan, On the strong solution of the ghost effect system, SIAM J. Math. Anal., 49 (2017), pp. 3496--3526.
	
	\bibitem{Huang-Wang-Wang-Yang} F.-M. Huang, Y. Wang, Y. Wang, and T. Yang, Justification of limit for the Boltzmann equation related to Korteweg theory, Quart. Appl. Math., 74 (2016), pp. 719--764.
	
	\bibitem{Jiang-Masmoudi} N. Jiang and N. Masmoudi, Low Mach number limits and acoustic waves. In Handbook of mathematical analysis in mechanics of viscous fluids, pages 2721--2770. Springer, Cham, 2018.
	
	\bibitem{Ju-Ou-2022} Q.-C. Ju and Y.-B. Ou, Low mach number limit of navier-stokes equations with large temperature variations in
	bounded domains. J. Math. Pures Appl, 164 (2022) 131--157.
	
	\bibitem{Levermore-Sun-Trivisa} C. D. Levermore, W. Sun and K. Trivisa, Local well-posedness of a ghost effect system, Indiana Univ. Math. J., 60 (2011), 517--576.
	
	%\bibitem{GT} D. Gilbarg, N.-S. Trudinger, {\it Elliptic Partial Differential Equations of Second Order}, Springer, 1998.
	
	\bibitem{Li-Liao-2019} Y.-P. Li and J. Liao, Existence of strong solutions to the stationary compressible Navier-Stokes-Korteweg equations with large external force. Nonlinear Anal. Real World Appl. 47 (2019), 204--223.
	
	\bibitem{Lions} P.-L. Lions, {\it Mathematical Topics in Fluid Dynamics, vol. 2, Compressible Models}, Oxford Science Publications, Oxford, 1998.
	
	\bibitem{Novotny-Padula-1994} A. Novotn\'{y} and M. Padula,
	$L^p$-approach to steady flows of viscous compressible fluids in exterior domains. Arch. Rational Mech. Anal.126 (1994), 243--297.
	
	\bibitem{Novotny-Straskraba} A. Novotn\'{y} and I. Stra{s}kraba, {\it Introduction to the Theory of Compressible Flow}, Oxford University Press, Oxford, 2004.
	
	%\bibitem{Sone-2002} Y. Sone, Kinetic theory and fluid dynamics, Birkhauser Boston, Inc., Boston, MA, 2002.
	
	\bibitem{Sone-2007} Y. Sone, Molecular gas dynamics. Theory, techniques, and applications., Birkhauser Boston, Inc., Boston, MA, 2007.
	
	\bibitem{Sun-2022} C.-Z. Sun, Uniform regularity in the low Mach number and inviscid limits for the full Navier-Stokes system in domains with boundaries, arXiv: 2204.09799v2.
	
	\bibitem{Valli-1987} A. Valli,
	On the existence of stationary solutions to compressible Navier-Stokes equations. Ann. Inst. H. Poincaré Anal. Non Linéaire 4 (1987), no.1, 99--113.
	
	%\bibitem{Wu-Ouyang} L. Wu and Z.-M. Ouyang, Asymptotic analysis of Boltzmann equation in bounded domains, arXiv:2008.10507.
\end{thebibliography}
\end{document}